\documentstyle[12pt,leqno,german,epsfig,twoside]{mystyle}

\makeindex

\title{Komplexe elliptische Geschlechter und $S^1$-"aquivariante
Kobordismustheorie}

\author{Diplomarbeit\\ von\\ Gerald H"ohn}

\date{Bonn und Vallendar, August 1991}

\begin{document}




\parindent=0cm


\multiply\baselineskip by 4
\divide\baselineskip by 3
\parskip=2mm


\hyphenation{
Man-nig-fal-tig-keit Ko-dimen-sion Ge-schlecht Man-nig-fal-tig-keit-en
ell-ip-tisch Auf-bla-sen Ab-bil-dung kom-plex kom-pakt Ein-bettung
Chern-klas-se Ho-mo-mor-phis-mus Grup-pe For-mel Fah-nen-man-nig-fal-tig-keit
Iso-mor-phis-mus Vor-les-ungs-mit-schrift Ein-bett-ung-en
Tan-gen-tial-bun-del uber-ein Ko-bor-dis-mus Ide-al dif-fe-ren-zier-bar
dif-fe-ren-zier-bare Lift-ung-en}

\hsize=16cm
\vsize=30cm

\addtolength{\topmargin}{-65pt}
\addtolength{\textheight}{+140pt}

\tolerance=1000     
\hyphenpenalty=10
\pretolerance=-1

\oddsidemargin -9pt
\evensidemargin -9pt


\def\R{\relax\ifmmode{\bf R}\else {$\bf R\ $}\fi}
\def\C{\relax\ifmmode{\bf C}\else {$\bf C\ $}\fi}
\def\H{\relax\ifmmode{\bf H}\else {$\bf H\ $}\fi}
\def\N{\relax\ifmmode{\bf N}\else {$\bf N\ $}\fi}
\def\Z{\relax\ifmmode{\bf Z}\else {$\bf Z\ $}\fi}
\def\Q{\relax\ifmmode{\bf Q}\else {$\bf Q\ $}\fi}
\def\L{\relax\ifmmode{\cal L}\else {$\cal L\ $}\fi}
\def\A{\relax\ifmmode{\hat A}\else {$\hat A\ $}\fi}
\def\CP{\relax\ifmmode{\bf CP}\else {$\bf CP\ $}\fi}
\def\rad{\mathop{\rm Rad}\nolimits}
\def\res{\mathop{\rm res}\nolimits}
\def\ker{\mathop{\rm Ker}\nolimits}
\def\sign{\mathop{\rm sign}\nolimits}
\def\rg{\mathop{\rm rg}\nolimits}
\def\rang{\mathop{\rm rang}\nolimits}
\def\kdim{\mathop{\rm K\hbox{-}dim}\nolimits}
\def\Td{\mathop{\rm Td}\nolimits}
\def\ch{\mathop{\rm ch}\nolimits}
\def\cod{\mathop{\rm codim\bf_C}\nolimits}
\def\dimc{\mathop{\rm dim\bf_C}\nolimits}
\def\grad{\mathop{\rm grad}\nolimits}
\def\id{\hbox{id}}
\def\iso{\cong}
\def\homotop{\simeq}
\def\cP{{\cal P}}
\def\cQ{{\cal Q}}
\def\cR{{\cal R}}
\def\Frage#1{$$\quad \hbox{\bf ???????? --- Da hab ich noch ein Problem --- ????????} \quad $$\write16{#1}}
\def\Luecke#1{$$\quad \hbox{\bf !!!!!!!! --- hier fehlt noch was --- !!!!!!!!} \quad $$\write16{#1}}

\def\SZDOLLAR#1{\if$#1\ $\else {\ #1}\fi}
\def\SZS#1{\if;#1;\ \else \SZDOLLAR{#1}\fi}
\def\SZD#1{\if:#1:\ \else \SZS{#1}\fi}
\def\SZK#1{\if,#1,\ \else \SZD{#1}\fi}
\def\SZ#1{\if.#1.\ \else \SZK{#1}\fi}

\def\MGF#1{Mannigfaltigkeit\SZ{#1}}
\def\MGFS#1{Mannigfaltigkeiten\SZ{#1}}
\def\UMGF#1{Untermannigfaltigkeit\SZ{#1}}
\def\VB#1{Vektorb\"undel\SZ{#1}}   
\def\VBS#1{Vektorb\"undels\SZ{#1}}
\def\GE#1{Geschlecht\SZ{#1}}
\def\SO#1{$S^1$-Operation\SZ{#1}}
\def\SOS#1{$S^1$-Operationen\SZ{#1}}
\def\SB#1{$S^1$-B"undel\SZ{#1}}
\def\NM#1{$N$-Mannigfaltigkeit\SZ{#1}}
\def\NMS#1{$N$-Mannigfaltigkeiten\SZ{#1}}
\def\NS#1{$N$-Struktur\SZ{#1}}
\def\NSS#1{$N$-Strukturen\SZ{#1}}
\def\ZSHGD#1{zusammenh"angend\SZ{#1}}
\def\ZSHGDE#1{zusammenh"angende\SZ{#1}}

\def\lto{\longrightarrow}
\def\ST{\tilde{ S^1}}
\def\TWIST{\widetilde{\CP}(E\oplus F)}
\def\GU{\varphi_{ell}}
\def\GN{\varphi_N}
\def\OU{\Omega_*^U}
\def\OSU{{\Omega_*^{SU}}}
\def\OUQ{\Omega_*^U\otimes\Q}
\def\OUN{\Omega_*^{U,N}}
\def\OUNC{\Omega_*^{U,N}\otimes{\bf C}}
\def\OUNQ{\Omega_*^{U,N}\otimes{\bf Q}}
\def\OSUQ{{\Omega_*^{SU}\otimes{\bf Q}}}
\def\OSOQ{{\Omega_*^{SO}\otimes{\bf Q}}}
\def\OUC{\Omega_*^U\otimes{\bf C}}
\def\OSUC{{\Omega_*^{SU}\otimes{\bf C}}}
\def\JSU{J_*^{SU}}
\def\ISU{I_*^{SU}}
\def\X1#1{\overline{{\bf H}/\Gamma_1(#1)}}
\def\G1#1{\Gamma_1(#1)}
\def\ZN{\Z/N\Z}
\def\DCUP{\mathop{\dot{\hbox{$\bigcup$}}}}
\def\Dcup{\mathop{\dot{\hbox{$\cup$}}}}
\def\c1{\bar c_1}
\def\JC{J_*^{\bf C}(PSL_2({\bf Z}))}
\def\JZ{J_*^{\bf Z}(PSL_2({\bf Z}))}
\def\CHIYL{\chi_y(q,{\cal L}X_d)}
\def\AN{{\cal A}_N}
\def\FC{\overline{F}}

\def\SN#1{\bf Satz\ #1: \sl\leftskip=0mm}
\def\LN#1{\bf Lemma\ #1: \sl\leftskip=0mm}
\def\DN#1{{\bf Definition\ #1: \sl\leftskip=0mm}}
\def\SE{\rm\leftskip=0mm}

\newbox\nitmit
\def\notmid{%
\mathrel{\mathchoice%
{\setbox\nitmit\hbox{$\textstyle/$}%
\hbox to 0pt{\hbox to \wd\nitmit{\hss$\textstyle/$\hss}\hss}%
\hbox to \wd\nitmit{\hfil$\textstyle\mid$\hfil}}%
{\setbox\nitmit\hbox{$\textstyle/$}%
\hbox to 0pt{\hbox to \wd\nitmit{\hss$\textstyle/$\hss}\hss}%
\hbox to \wd\nitmit{\hfil$\textstyle\mid$\hfil}}%
{\setbox\nitmit\hbox{$\scriptstyle/$}%
\hbox to 0pt{\hbox to \wd\nitmit{\hss$\scriptstyle/$\hss}\hss}%
\hbox to \wd\nitmit{\hfil$\scriptstyle\mid$\hfil}}%
{\setbox\nitmit\hbox{$\scriptscriptstyle/$}%
\hbox to 0pt{\hbox to \wd\nitmit{\hss$\scriptscriptstyle/$\hss}\hss}%
\hbox to \wd\nitmit{\hfil$\scriptscriptstyle\mid$\hfil}}%
}}
\def\mapr#1{\smash{\mathop{\lto}\limits^{#1}}}
\def\mapd#1{\Big\downarrow\rlap{$\vcenter{\hbox{$\scriptstyle#1$}}$}}
\def\klein#1{\hbox{$\scriptstyle{#1}$}}

\def\kleinn#1{\hbox{$\footnotesize#1$}}

\def\BOX{\quad\phantom.\hfill\hbox{\vrule height 5pt width 5pt}\relax}

\maketitle

\newpage

{\small\tt {\bf Anmerkung zu der auf den Archivserver abgelegten Version der
Diplomarbeit:} 

Dies ist eine Rekonstruktion der Originalversion mittels
alter Sicherungs\-disket\-ten. Einige Seitenumbr"uche sind wegen der
eingeschr"ankten R"uckw"arts\-kompatibilit"at von LaTeX2e
leicht verschieden.

In den folgenden Formeln bleiben einige Pfeile
zu erg"anzen:

\begin{itemize}
\item[-] S.~8, Def.~1.1.3 (2): von $Q$ nach $B$ und von $P$ nach $B$; 
\item[-] S.~9, unten: von $B$ nach $\widetilde{BU}(n)$ und von $B$
                    nach $BU(n)$;
\item[-] S.~19, Gleichung (1.2): von $E \oplus F$ nach $B$;
\item[-] S.~45, nach Def.~2.3.3.: von $\Omega_*^U\otimes {\bf C}$ nach $K(C_N)$;
\item[-] S.~73, gro"ses Diagramm: von $F(\nu)$ nach $\tilde Y$ und von
                $\tilde Y$ nach $Y$.
%
\end{itemize} }
\newpage

\pagenumbering{roman}
\setcounter{page}{3}
\addcontentsline{toc}{chapter}{\protect\numberline{Vorwort}}

\vspace*{30pt}
{\Large\bf Vorwort\hfill}

\medskip

In der Diplomarbeit werden kobordismustheoretische Fragestellungen
mit Hilfe von komplexen elliptischen Geschlechtern behandelt sowie
Eigenschaften von diesen untersucht.

Die Theorie der elliptischen Geschlechter lernte ich in einer Vorlesung
kennen, die Prof.~F.~Hirzebruch im WS 1987/88 an der Universit"at Bonn hielt.
In der Vorlesung wurde einerseits eine Einf"uhrung
in das damals neue Gebiet der elliptischen Geschlechter gegeben, anderseits
wurde eine Verallgemeinerung auf komplexe \MGFS vorgenommen. Die so erhaltenen
komplexen
Geschlechter der Stufe $N$ ordnen einer komplexen \MGF eine Modulform
zur Modulgruppe $\G1{N}$ zu. F"ur $N=2$ erh"alt man das urspr"ungliche
f"ur differenzierbare \MGFS erkl"arte
elliptische Geschlecht, wie es von S.~Ochanine, E.~Witten und anderen
eingef"uhrt wurde. Die wesentlichen Eigenschaften --- Ganzzahligkeit und
Starrheit ---, die im Stufe-$2$-Fall f"ur $Spin$-\MGFS erf"ullt sind, gelten
im Stufe-$N$-Fall, wenn die komplexe Mannigfaltigkeit --- oder allgemeiner die
stabil fastkomplexe Mannigfaltigkeit --- eine durch $N$ teilbare erste Chernklasse besitzt.

Ziel der Diplomarbeit ist es, m"oglichst viele der im
differenzierbaren Fall sonst noch
bekannten Ergebnisse in geeigneter Weise auf den
komplexen Fall zu "ubertragen.

Ich m"ochte R.~Jung f"ur die Diskussionen, die viel zur Entstehung der Arbeit
beigetragen haben, sowie f"ur sonstige Hilfen danken. Weiter danke ich
Herrn Prof.~Kreck f"ur die Beantwortung vieler Fragen,
Herrn Prof.~W.~Vogel aus Halle f"ur seine Hilfe beim Beweis von Satz 2.3.5,
H.~von Eitzen f"ur die "Ubersetzung von \cite{kona}, H.~Kastenholz f"ur seine
ausf"uhrlichen Korrekturvorschl"age, sowie allen anderen,
mit denen ich "uber meine Arbeit sprechen konnte.

Ganz besonders gilt mein Dank
Herrn Prof.~Hirzebruch f"ur die Anregung zu dieser Arbeit und seine begleitende
Hilfe. Aus seinen Vorlesungen und Arbeiten stammen auch
die meisten meiner Kenntnisse "uber Kobordismustheorie, charakteristische
Klassen und Geschlechter.

Schlie"slich halfen mir beim Aufschreiben und Korrekturlesen der Arbeit
R.~Kleinrensing, J.~Kettner und mein Vater. Ihnen sei ebenfalls gedankt.

\hfill {\sc Gerald H"ohn}

\newpage

\setcounter{page}{4}
\phantom{xxx}
\newpage

\setcounter{page}{5}
\tableofcontents

\newpage
\phantom{xxx}
\newpage
\pagenumbering{arabic}
\setcounter{page}{1}
\addcontentsline{toc}{chapter}{\protect\numberline{Einleitung}}

\vspace*{40pt}
{\Large\bf Einleitung\hfill}

\bigskip
\bigskip
Die Diplomarbeit besch"aftigt sich mit der Verallgemeinerung der Theorie
der elliptischen Geschlechter
von orientierten differenzierbaren \MGFS auf $U$-Man\-nig\-fal\-tig\-keiten
(\MGFS mit komplexer Struktur auf dem stabilen Tangentialb"undel).

Das universelle elliptische Geschlecht der Stufe zwei, $\varphi_2:\Omega^{SO}_*\otimes\C \to M_*
(\Gamma_1(2))\iso\C[\delta,\epsilon]$
(f"ur eine Einf"uhrung vgl.~ersten Aufsatz in \cite{land} sowie
\cite{botau,last,och1,och2,tau,witt1,witt2}),
l"a"st sich nach Ochanine \cite{och1}
wie folgt charakterisieren: sei $J^{Spin}_*\subset \Omega^{Spin}_*\otimes\C
\iso\Omega^{SO}_*\otimes\C$ das Ideal, das von projektiven B"undeln
$\CP(E)$ erzeugt wird, wobei $E$ ein komplexes \VB von
geradem Rang "uber einer orientierten \MGF ist. Dann gilt
$\ker \varphi_2=J^{Spin}_*$. Das Geschlecht $\varphi_2$ geh"ort zur charakteristischen Potenzreihe
$Q(x)=\frac{x}{f(x)}$, wobei die Umkehrfunktion
$f^{-1}(y)=g(y)$ von $f(x)$ das elliptische Integral
$g(y)=\int_0^y{\frac{dt}{\sqrt{1-2\delta t^2+\epsilon t^4}}}$ ist.
Witten \cite{witt1} zeigte, da"s sich $\varphi_2(X)$ als
die "aquivariante Signatur des freien Schleifenraumes $\sign(q,\L X)$
interpretieren l"a"st, welche man formal als eine Reihe $\sum_{n\ge 0} \sign(X,R_n)q^n$
schreiben kann. Hierbei sind die $R_n$ eine Folge
zu $TX$ assoziierter \VB.
Witten vermutete, da"s das elliptische Geschlecht starr bei \SOS
auf Spin-\MGFS
ist, d.h.~der virtuelle $S^1$-Charakter $\sign(\lambda,X,R_n)\in\Z[\lambda,\lambda^{-1}]$
ist trivial f"ur alle $n$. Dies wurde von Ochanine \cite{och2} im Falle von
semifreien \SOS und im allgemeinen Fall von Taubes \cite{tau,botau} bewiesen.
Weiter folgt, da"s der Kern von $\varphi_2$ das Ideal $I^{Spin,odd}_*$ ist, welches von
zusammenh"angenden Spin-\MGFS mit ungerader \SO erzeugt wird.
F"ur semifreie \SOS vgl.~dazu auch \cite{bos}.
Die Starrheit impliziert Multiplikativit"at in allen \hbox{orientierten} Faserb"undeln
mit Spin-\MGFS als Faser: $\varphi_2(E)=\varphi_2(B)\cdot\varphi_2(F)$
(Genauer sei $E$ ein lokal triviales orientiertes Faserb"undel mit orientierter \MGF
$B$ als Basis, orientierter \MGF $F$ als Faser und kompakter
zusammenh"angender Liegruppe als Strukturgruppe.) Falls die Faser eine
beliebige orientierte \MGF sein darf,
ist die Signatur das einzige Geschlecht mit dieser Eigenschaft
(vgl.~\cite{chhise,bohi}). Operiert auf $E$ eine endliche Gruppe, so ist auch die
"aquivariante Signatur multiplikativ \cite{ossa}, so da"s man vermuten k"onnte,
da"s dies auch f"ur $\sign(q,\L X)$ (kleine Schleifen verwenden) der Fall ist.
Da der freie Schleifenraum $\L X$ nur f"ur Spin-\MGFS $X$ orientierbar ist
\cite{segal}, k"onnte dies eine Erkl"arung f"ur die Spin-Bedingung liefern.

Das universelle elliptische Geschlecht $\varphi_2(X)$ ist eine Modulform zu $\Gamma_1(2)$,
so da"s man die Werte in den beiden Spitzen von $\G1{2}$ betrachten kann. Diese
Werte sind aber gerade die Signatur und das $\A$-Geschlecht. Die Starrheit des
$\A$-Geschlechtes bei \SOS auf Spin-\MGFS  wurde von Atiyah-Hirzebruch
\cite{athis1} gezeigt; au"serdem: $\ker \A=I_*^{Spin}$, wobei $I_*^{Spin}$
das Ideal ist, welches von
zusammenh"angenden Spin-\MGFS mit effektiver \SO erzeugt wird.
\medskip

Hirzebruch folgend \cite{himod,hiell} werden in meiner Diplomarbeit die
obigen Charakterisierungen von $\varphi_2$ in geeigneter Weise
auf komplexe elliptische Geschlechter "ubertragen.
Dazu definieren wir das universelle
elliptische Geschlecht $\GU$ als das Geschlecht, das zur
charakteristischen Potenzreihe $Q(x)=\frac{x}{f(x)}$ geh"ort, wobei $f$
die normierte L"osung der Differentialgleichung
$$\left(\left(\frac{f'}{f}\right)'\right)^2=S\left(\frac{f'}{f}\right)$$
mit
$$S(y)=\left(y+\frac{A}{2}\right)^4-\frac{1}{4}B\,\left(y+\frac{A}{2}\right)^2+4C\,\left(y+\frac{A}{2}\right)+\frac{1}{64}B^2-2D $$
ist.
Das universelle elliptische Geschlecht einer komplex $n$-di\-men\-sio\-nalen
$U$-Man\-nig\-fal\-tig\-keit ist dann ein
gewichtet homogenes Polynom vom Gewicht $n$ in den Unbestimmten $A$, $B$, $C$ und $D$,
wenn diese mit den Gewichten $1$ bis $4$ versehen werden. F"ur $SU$-\MGFS ist es
sogar ein Polynom, das nur von $B$, $C$ und $D$ abh"angt.
Obiges Polynom $S(y)$ ist gerade so
gew"ahlt worden, da"s $\GU$ auf bestimmten ausgezeichneten $U$-\MGFS $W_1$ bis
$W_4$ die Werte $A$ bis $D$ annimmt.

Das elliptische Geschlecht der Stufe $N$, wie von Hirzebruch eingef"uhrt, geh"ort
zu solchen Polynomen $S(y)$, die einer Differentialgleichung
$$\hbox{\quad(*)\quad} P'(y)^2S(y)=N^2y^2(P(y)^2-c^2)$$
gen"ugen mit einem Polynom
$P(y)=y^N+\dots+d_{N-1}y+d_N$ und einer Konstante $c\not=0$. Die Gleichung (*) liefert zwei homogene
Bedingungen $R_{N-1}$ und $R_{N+1}$ vom Gewicht $N-1$ und $N+1$ an die
Unbestimmten $A$ bis $D$.
Sie definieren eine Kurve $C_N$ im gewichtet projektiven Raum $\CP^{1,2,3,4}$.
Wie in \cite{jung} gezeigt, sind die irreduziblen Komponenten dieser Kurve gerade
die Bilder von Abbildungen der Modulkurven $\X1{n}$ in den $\CP^{1,2,3,4}$,
wobei $n>1$ die Teiler von $N$ durchl"auft.
Mit $\tilde\varphi_N$ bezeichnen wir dann den wie folgt definierten
Algebrenhomomorphismus von $\OUC$ in die Koordinatenalgebra $\C[A,B,C,D]/I(C_N)$
von $C_N$:\hfill\break
Wir konstruieren eine Basisfolge $W_1,W_2,\dots$  von $\OUC$ mit $\GU(W_1)=A$,
$\GU(W_2)=B$, $\GU(W_3)=C$, $\GU(W_4)=D$ und $\GU(W_i)=0$ f"ur $i\ge 5$. Sei
$$\pi_N: \C[A,B,C,D]\to\C[A,B,C,D]/I(C_N)$$ die Projektion
auf die Koordinatenalgebra von $C_N$.
Dann definiere $\tilde\GN$ als $\pi_N\circ\GU$.
Das in \cite{hiell} definierte Geschlecht $\GN$ ist dann die Abbildung $\pi_N'\circ\GU$,
$\pi_N': \C[A,B,C,D]\to\C[A,B,C,D]/I(C_N')$, wobei $C_N'$ die irreduzible Komponente
von $C_N$ ist, die der Modulkuve $\X1{N}$ entspricht. Alternativ kann
$\GN$ auch als Polynom in Modulformen $A$ bis $D$ vom Gewicht 1 bis 4 zu $\Gamma_1(N)$ aufgefa"st
werden. Die Werte von $\varphi_N$ in den Spitzen von $\G1{N}$ sind die Geschlechter $\chi_y(X)$
mit $y=e^{2\pi i \frac{l}{N}}$ und  $l=1,\dots,N-1$, $(l,N)=1$,
sowie die Geschlechter $\chi(X,K^{k/N})$ mit $k=1,\dots,N-1$. Die Geschlechter
$\chi(X,K^{k/N})$ fassen wir zu einem Geschlecht $\tilde A_N$ zusammen:\hfill\break
Es ist das Geschlecht $\mu_N\circ\GU$, wobei $\mu_N:\C[A,B,C,D]\to\C[A,B,C,D]/I(M_N)$
die Projektion auf die Koordinatenalgebra von $M_N=\bigcup_{k=1}^{[\frac{N}{2}]}P_{N,k}
\subset\CP^{1,2}\subset\CP^{1,2,3,4}$ ist. Die Variet"at $M_N$ besteht aus den
Punkten $P_{N,k}=(2(1/2-k/N)t:2t^2:0:0)$, deren homogene Koordinaten die Werte der
Modulformen $A$ bis $D$ in den Spitzen von $\G1{N}$ sind,
die zu $\chi(X,K^{k/N})$ geh"oren.

Sei $\OUN$ der Kobordismusring von stabil fastkomplexen \MGFS mit fest ge\-w"ahl\-ter
$N$-ten Wur\-zel des De\-ter\-mi\-nan\-ten\-b"un\-dels, d.h.~die Objekte sind
$U$-Man\-nig\-fal\-tig\-keit\-en $X$ zusammen mit einer $N$-Struktur:
das ist ein Geradenb"undel $L$, so da"s $L^N\iso\det(TX)$.
Besitzt eine $U$-\MGF $X$ eine \NS, so gilt $c_1(X)\equiv 0\pmod{N}$. Wir bezeichnen
$U$-\MGFS mit durch $N$ teilbarer erster Chernklasse daher als \NMS.
Umgekehrt besitzen \NMS stets
eine (nicht notwendig eindeutig bestimmte) $N$-Struktur.
Mit dem Pontrjagin-Thom-Theorem berechnet man
$\OUNQ\iso\OUQ\iso\Q[W_1,W_2,\dots]$.

Der Typ $t$ einer effektiven, mit der $U$-Struktur kompatiblen
\SO  auf einer zusammenh"angenden
$N$-\MGF $X$ mit $[X]\not=0$ in $\OUNQ$ ist die Summe der Drehzahlen im Normalenb"undel
einer Fixpunktkomponente modulo $N$. Er ist unabh"angig von der Auswahl der
Fixpunktkomponente.
Sei $I_*^{N,t}$ das Ideal in $\OUNC$, welches von zusammenh"angenden $N$-\MGFS
mit effektiven \SOS vom Typ $t$ erzeugt wird. Da $I_*^{N,t}=I_*^{N,ggT(N,t)}$ gilt,
betrachten wir nur $t=0$ und Teiler $t$ von $N$. Das Ideal $I_*^{N}=I_*^{N,0}$
sei das Ideal der \NMS mit effektiver \SO. Weiter definiere noch ein Ideal
$J_*^{N}$ wie folgt: F"ur zwei komplexe \VB $E$ und $F$ vom Rang $p$ und $q$ "uber einer
$U$-\MGF $B$ sei das getwistet projektive B"undel $\TWIST$ als differenzierbare
\MGF das gew"ohnliche projektive B"undel $\CP(E\oplus \FC)$. Es wird aber mit einer anderen,
"`getwisteten"' stabil fastkomplexen Struktur versehen, so da"s eine Faser
$\widetilde\CP_{p,q}$ dieses B"undels die totale Chernklasse
$c(\widetilde\CP_{p,q})=(1+g)^p(1-g)^q$ besitzt. Das Ideal $J_*^N$ von
$\OUNC$ wird nun von solchen B"undeln $\TWIST$ erzeugt, bei denen
der Totalraum und die Faser \NMS sind. (Die zweite Bedingung
ist "aquivalent zu $N\mid p-q$.) Schlie"slich definieren wir noch analog die Ideale
$I_*^{SU}$, $I_*^{SU,t}$ und $J_*^{SU}$ f"ur $SU$-\MGFS. Der Typ $t$ ist
hierbei eine ganze Zahl.

Ein Hauptresultat der Arbeit ist der

{\bf Satz:\sl}
\vskip -15mm
\begin{eqnarray}
(1)&&\qquad \ker \GU\mid_{\OSUC}=\JSU=I_*^{SU,t},\ \hbox{f"ur}\ t\not=0\cr
(2)&&\qquad \ker \hat A\mid_{\OSUC}=\ISU\cr
(3)&&\qquad \ker \tilde\GN=J_*^N=I_*^{N,1}\cr
(4)&&\qquad \ker \tilde A_N=I_*^N\nonumber
\end{eqnarray}
\SE

Die Inklusion $\ker \tilde A_N\supset I_*^N$ findet sich bei Hattori \cite{hat}
und Kri\v cever \cite{kri}, vgl.~auch \cite{may}, die Inklusion $\ker \tilde\GN\supset I_*^{N,1}$ bei
Hirzebruch \cite{hiell}.
Wir haben damit also eine geometrische Charakterisierung der komplexen
elliptischen Geschlechter $\GU$ und $\tilde\GN$.

Da $\hbox{Bild}(\varphi_{ell}\!\mid_{\OSUC})\,\,\subset\C[B,C,D]$, k"onnen wir $\GU$ als ${\bf C}$-Algebrenhomomorphismus
von $\OSUC$ in die Koordinatenalgebra der gewichtet projektiven Ebene $\CP^{2,3,4}$
auffassen. Die Ebene wiederum k"onnen wir mit dem geeignet kompaktifizierten Modulraum
von "Aquivalenzklassen elliptischer Kurven "uber \C mit ausgezeichnetem Punkt auf der Kurve identifizieren.
Hierbei entspricht dem Gitter $L=2\pi i (\Z\tau+\Z)$ mit ausgezeichnetem
Punkt $z\in\C/L$, $z\not=0$ der Punkt
$(B:C:D)=(24\,\wp_L(z):\wp_L'(z):6\,\wp_L^2(z)-\frac{1}{2}g_2(L))\in\CP^{2,3,4}$.
Dabei steht $\wp_L(x)$ f"ur die Weierstra"ssche $\wp$-Funktion zum Gitter $L$,
und  $g_2(L)$, $g_3(L)$ seien die Gitterkonstanten.
Unter dieser Identifikation ist das universelle elliptische Geschlecht f"ur
$SU$-\MGFS
bis auf einen Homogenit"atsfaktor $\lambda^{-1}\in\C^*$, der dem
"Ubergang von $(L,z)$ zu $(\lambda L,\lambda z)$ entspricht,
das Geschlecht zur charakteristischen Potenzreihe
$$Q(x)=
\frac{x}{1-e^{-x}}(1+ye^{-x})\prod_{n=1}^{\infty}\frac{1+yq^ne^{-x}}
{1-q^ne^{-x}}\cdot\frac{1+y^{-1}q^ne^{x}}{1-q^ne^{x}}
\cdot\Bigl(\!(1+y)\!\prod_{n=1}^{\infty}\frac{(1+yq^n)(1+y^{-1}q^n)}
{(1-q^n)^2}\Bigr)^{-1}\!\!,$$
wobei $q=e^{2\pi i \tau}$ und $-y=e^{z}$.
So geschrieben ist $\GU(X_d)$ f"ur eine $d$-dimensionale $SU$-\MGF $X_d$
eine Jacobiform vom Gewicht $d$ und Index $0$ (holomorph bis auf einen Pol in $0$)
zur vollen Modulgruppe $PSL_2(\Z)$ (vgl.~\cite{eizag}). Das unnormierte universelle
elliptische Geschlecht $\CHIYL$\index{$\CHIYL$}
l"a"st sich wegen obiger Produktentwicklung auch als Potenzreihe in
$q$ und $y$ schreiben mit Indizes getwisteter
Dolbeault-Komplexe als Koeffizienten:
\begin{eqnarray}\CHIYL& =&\prod_{i=1}^d
\frac{x_i}{1-e^{-x_i}}(1+ye^{-x_i})\prod_{n=1}^{\infty}\frac{1+yq^ne^{-x_i}}
{1-q^ne^{-x_i}}\cdot\frac{1+y^{-1}q^ne^{x_i}}{1-q^ne^{x_i}}[X_d]\cr
\noalign{\vskip1mm}
&=& \chi_y\big(X_d,\bigotimes_{n=1}^{\infty}\Lambda_{yq^n}T^*
\otimes\bigotimes_{n=1}^{\infty}\Lambda_{y^{-1}q^n}T\otimes
\bigotimes_{n=1}^{\infty}S_{q^n}(T\oplus T^*)\big)
\cr
\noalign{\vskip2mm}
&=&\sum_{n\ge 0}\sum_{\mid i \mid\le c_n}\chi(X,R_{n,i})y^iq^n.\nonumber
\end{eqnarray}
Die $R_{N,i}$ sind dabei zum stabil fastkomplexen Tangentialb"undel $T$ assoziierte
komplexe \VB, und die $x_i$ sind die formalen Wurzeln der Chernklasse von $X_d$.
Man kann formal $\CHIYL$ als das $S^1$-"aquivariante $\chi_y$-Geschlecht
des freien Schleifenraumes ${\cal L}X_d$ interpretieren.

In Verallgemeinerung der Geschlechter der Stufe $N$ gilt der

{\bf Satz:}\sl~Das universelle komplexe elliptische Geschlecht ist starr bei \SOS
auf $SU$-\MGFS, d.h.~der virtuelle $S^1$-Charakter $\chi(\lambda,X,R_{n,i})\in\Z[\lambda
,\lambda^{-1}]$ ist f"ur alle $n$ und $i$ trivial.
\SE

Die Starrheit bei \SOS ist "aquivalent zu der Multiplikativit"at in stabil fastkomplexen
Faserb"undeln mit $SU$-\MGFS als Faser und kompakter zusammenh"angender
Liegruppe von $U$-Automorphismen der Faser als Strukturgruppe. Der Koeffizient
von $q^0$ in der Potenzreihenentwicklung von $\CHIYL$ ist das
$\chi_y$-Geschlecht, das multiplikativ in allen $U$-B"undeln
ist (vgl.~\cite{hihab}).
Die Bedingung $c_1=0$ an die Faser k"onnte man als das Hindernis f"ur eine
$U$-Struktur auf ${\cal L}X$ zu erkl"aren versuchen. Die Interpretation von
$\CHIYL$
als das "aquivariante $\chi_y$-Geschlecht des freien Schleifenraumes liefert
Zusammenh"ange mit der konformen Quantenfeldtheorie (Charaktere von
Quantengruppen (vgl.~\cite{roan})?).
Es schlie"st sich die Frage an, ob es ein komplexes Analogon zur elliptischen
Kohomologie im Stufe-$2$-Fall gibt.

Schlie"slich beweisen wir noch eine Aussage "uber das Verhalten von $\GN$
beim Aufblasen:

{\bf Satz:}\sl~Sei $Y$ eine kompakte komplexe Untermannigfaltigkeit der kompakten komplexen
\MGF $X$, sei $\tilde X$ die Aufblasung von $X$ entlang $Y$. Falls $\cod Y
\equiv 1 \pmod{N}$, gilt $\GN(\tilde X)=\GN(X)$.
\SE



\newpage

\chapter{Der Kobordismusring $ \Omega_*^{U,N} $ der N-Mannigfaltigkeiten}

\section{Geometrische Definition von $ \Omega_*^{U,N} $}
In diesem Abschnitt referieren wir kurz "uber die Kobordismusringe
$\OU$ und $\OSU$. Es wird dann der Kobordismusring $\OUN$ der \NMS mit \NS eingef"uhrt.
Dies sind \MGFS mit einer $U$-Struktur auf dem stabilen Tangentialb"undel
und fest gew"ahlter "`$N$-ten Wurzel"' des zum stabilen Tangentialb"undel geh"origen
komplexen Determinantenb"undels.

Wir erinnern an die Definition von $U$-\MGFS (vgl.~z.B.~\cite{cofl2}).
Sei $X$ eine kompakte differenzierbare (d.h.~$C^{\infty}$) \MGF.
Eine komplexe Struktur  auf dem stabilen Tangentialb"undel von $X$ ist ein
Vektorb"undelisomorphismus $J:TX\oplus\epsilon_{\R}^k\to TX\oplus\epsilon_{\R}^k$
mit $J^2=-\id$, wobei $\epsilon_{\R}^k$
das triviale reellen B"undels vom Rang $k$ "uber $X$ ist.
Durch eine komplexe Struktur wird $TX\oplus\epsilon_{\R}^k$ zu einem komplexen
\VB.
Ist $J$ eine komplexe Struktur auf $TX\oplus\epsilon_{\R}^k$ und
$J'$ eine komplexe Struktur auf $TX\oplus\epsilon_{\R}^{k'}$, so geh"oren
$J$ und $J'$ der gleichen stabilen "Aquivalenzklasse $\Phi$ an, falls es $s$
und $s'$ mit $k+2s=k'+2s'=:t$ gibt, so da"s die komplexen Strukturen auf
$TX\oplus\epsilon_{\R}^t$, die zu den komplexen \VB
$TX\oplus\epsilon_{\R}^k\oplus\epsilon_{\C}^s$ und  $TX\oplus\epsilon_{\R}^{k'}\oplus\epsilon_{\C}^{s'}$
geh"oren ($\epsilon_\C^k$ sei das triviale reelle \VB vom Rang $2k$ mit der
standardkomplexen Struktur), stetig ineinander deformiert werden k"onnen.
Eine $U$-Struktur auf $X$ ist eine "Aquivalenzklasse $\Phi$ von
komplexen Strukturen auf dem stabilen Tangentialb"undel. Eine \MGF zusammen mit
einer $U$-Struktur bezeichnen wir als $U$-\MGF.
Die Menge der $U$-Strukturen auf $X$ kann mit der Menge der Homotopieklassen
von Liftungen der klassifizierenden Abbildung $X\to BO$ des stabilen
Tangentialb"undels bzgl.~der Faserung $BU\to BO$ identifiziert werden. F"ur eine
reell $n$-dimensionale \MGF $X_n$ k"onnen anstelle der stabilen R"aume $BO$ und
$BU$ die klassifizierenden R"aume $BO(2k)$ und $BU(k)$ verwendet werden,
falls $2k\ge n+2$.

Der Rand einer $U$-\MGF ist wieder eine $U$-\MGF.
Das Negative einer $U$-\MGF $X$, deren $U$-Struktur durch eine komplexe Struktur
auf $TX\oplus\epsilon_{\R}^k$ repr"asentiert wird, ist $X$ zusammen mit der
$U$-Struktur, die durch die komplexe Struktur
$TX\oplus\epsilon_{\R}^k\oplus\overline{\epsilon}_{\C}^1$ rep"asentiert wird.
Zwei $U$-\MGF hei"sen isomorph, falls
das Differential eines Diffeomorphismus die $U$-Struktur respektiert.
Disjunkte Vereinigung und direktes Produkt zweier $U$-\MGFS sind wieder $U$-\MGFS.
(Beim direkten Produkt zweier \MGFS mit Rand sind die "`Ecken zu gl"atten"'.
Dies ist bis auf $U$-Isomorphie eindeutig m"oglich.)

Ein komplexes \VB $E$ "uber einer $n$-dimensionalen \MGF $X_n$ besitzt die totale
Chernklasse $c(E)=1+c_1+c_2+\cdots+c_{[\frac{n}{2}]}\in H^*(X,\Z)$, $c_i\in
H^{2i}(X,\Z)$. Wegen der Multiplikativit"at der Chernklasse (d.h.~$c(E\oplus F)
=c(E)\cdot c(F)$) h"angt die Chernklasse einer komplexen Struktur auf dem stabilen
Tangentialb"undel nur von der stabilen "Aquivalenzklasse $\Phi$ ab, d.h.~die
Chernklasse $c(X)=c(TX\oplus\epsilon_{\R}^k)$ einer $U$-\MGF ist wohldefiniert.
Eine $U$-Struktur auf $X$ definiert
eine Orientierung f"ur $X$. Wir k"onnen
daher f"ur eine $U$-\MGF $X_n$ kanonisch einen Erzeuger $[X_n]$ von
$H_n(X_n,\Z)\iso\Z$ (Fundamentalzyklus der Mannigfaltigkeit) auszeichnen. Sei
$I=(i_1,\dots,i_n)$, $\sum\nu\cdot i_{\nu}=n$ eine Partition von $n$. Die
Chernzahl $c_I[X_{2n}]$ einer $2n$-dimensionalen $U$-\MGF $X_{2n}$ ist
dann durch $c_I[X_{2n}]=(c_1^{i_1}(X)\cdot c_2^{i_2}(X)\cdot\dots\cdot
c_n^{i_n}(X))[X_{2n}]$ definiert.

Zwei $n$-dimensionale $U$-\MGFS $X_n$ und $X'_n$ hei"sen kobordant, falls es
eine $(n+1)$-dimensionale $U$-\MGF $Y_{n+1}$ gibt mit $\partial Y_{n+1}=X_n-X_n'$.
Dies ist eine "Aquivalenzrelation und sie ist mit der disjunkten Vereingung und dem
direkten Produkt vertr"aglich. Die "Aquivalenzklassen bilden einen graduiert
kommutativen Ring $\OU$. (Addition $\hat =$ disjunkte Vereingung, Multiplikation
$\hat =$ Produkt.)

"Uber die Struktur von $\OU$ informiert

\SN{1.1.1 (Milnor)}\index{$\OU$}\index{$s(X_n)$} Sei $M_2$, $M_4$, $M_6$, $\dots$
eine Folge von $U$-\MGFS mit \hfill
$$s(M_{2n})=\cases{\pm 1 & falls $n+1$ keine Primzahlpotenz ist,\cr \pm p &
falls $n+1$ eine Potenz der Primzahl $p$ ist,}$$
wobei die Milnorzahl
$s(M_{2n})$ f"ur eine $2n$-dimensionale $U$-\MGF der Wert des Geschlechtes ist,
das zur charakteristischen Potenzreihe $Q(x)=1+x^n$ geh"ort (siehe Abschnitt 2.1).
Dann gilt $\OU\iso\Z[M_2,M_4,M_6,\dots]$.
Die Gruppen $\Omega_n^U$ sind insbesondere torsionsfrei.
\SE

{\bf Korollar:} \sl $\OUQ\iso\Q[M_2,M_4,M_6,\dots]$, falls $s(M_{2n})\not=0$.\SE

Eine solche Folge $M_{2n}$ von \MGFS mit $s(M_{2n})\not=0$ wird als Basisfolge
f"ur $\OUQ$ bezeichnet.

{\bf Beispiel:} $\OUQ\iso\Q[\CP_1,\CP_2,\CP_3,\dots]$, da
$s(\CP_n)=(1+g^n)^{n+1}[\CP_n]=n+1\not=0$, falls $g$ der Erzeuger von
$H^2(\CP_n,\Z)$ ist, der poincar\'e-dual zu einem Hyperebenenschnitt mit
der durch die komplexe Struktur induzierten Orientierung ist.

Analog zu $U$-\MGFS k"onnen $SU$-\MGFS definiert werden: Eine $SU$-Struktur
auf $X$ ist eine Homotopieklasse von Liftungen der klassifizierenden Abbildung
$X\to BO$ bez"uglich der Faserung $BSU\to BO$. Jede $SU$-Struktur auf $X$ induziert
eine $U$-Struktur. Eine $U$-\MGF l"a"st genau dann eine (nicht notwendig
eindeutige) $SU$-Struktur zu,
die die gegebene $U$-Struktur induziert, falls $c_1(X)=0$.
Als Strukturresultat f"ur den zugeh"origen Kobordismusring hat man

\SN{1.1.2}\index{$\OSU$} Es gilt $\OSUQ\iso\Q[M_4,M_6,M_8,\dots]$, falls $M_{2n}$ f"ur $n\ge2$
eine Folge von $SU$-\MGFS mit $s(M_{2n})\not=0$ ist.

{\bf Korollar:} \sl Sind Chernzahlen $c_I$ ($I=(i_1,\dots,i_n)$, $\sum \nu\cdot i_{\nu}=n$) vorgegeben,
so existiert genau dann ein $[X]\in\OSUQ$ mit diesen Chernzahlen, wenn alle
Chernzahlen mit $c_1$ als Faktor verschwinden.\SE

{\bf Beweis:} Die Notwendigkeit ist klar, die Umkehrung folgt aus Dimensionsgr"unden.
\BOX

Die Torsion von $\Omega_n^{SU}$ ist komplizierter zu beschreiben und zu
berechnen. Sie ist genau dann von Null verschieden, falls $n$ von der Gestalt
$n=8k+1$ oder $n=8k+2$ ist. Auch die Struktur des freien Anteils
von $\OSU$ ist wesentlich komplizierter als die von $\OU$.

Kommen wir nun zur Definition von $N$-\MGFS und $N$-Strukturen.
Dabei sei $N$ stets eine feste nat"urliche Zahl gr"o"ser oder gleich 2.\index{$N$}

\DN{1.1.3}\index{$N$-Struktur} Sei $B$ ein endlicher zusammenh"angender $CW$-Komplex und
$\xi$ ein $U(n)$-Prinzipalb"undel "uber $B$, sei $P$ das assoziierte
$S^1$-Prinzipalb"undel,
das zur Darstellung $\Lambda^n:U(n)\to U(1)=S^1$ geh"ort.
Eine \NS auf $\xi$ ist ein Paar $(Q,\pi)$, bestehend aus \hfill\break
\quad (1)\quad einem $S^1$-Prinzipalb"undel $Q$ "uber $B$; und\hfil\break
\quad (2)\quad einer Abbildung $\pi:Q\to P$, so da"s das folgende Diagramm
kommutativ ist:\typeout{2 Pfeile einf"ugen bei Def. Nstruktur}
$$\begin{array}{ccccc}
Q\times S^1              &\mapr{}& Q         &           &   \\
\mapd{\pi\times\lambda_N}&       &\mapd{\pi} &\qquad\quad& B \\
P\times S^1              &\mapr{}& P         &           &   \\
\end{array}.$$
Hierbei ist $\lambda_N$ die Abbildung von $S^1$ nach $S^1$, die $\lambda\in S^1$
auf $\lambda^N$ abbildet, und die horizontalen Pfeile sind die nat"urlichen
Rechtsoperationen von $S^1$ auf $P$ und $Q$. Eine \NS auf einem komplexen \VB
ist eine \NS auf einem zugeh"origen $U(n)$-Prinzipalb"undel.

Die Abbildung $\pi: Q\to P$ ist eine $N$-bl"attrige "Uberlagerung.
Sind $K$ und $L$ die zu $P$ und $Q$ assozierten Geradenb"undel, so ist $K$ die $N$-te
Tensorpotenz von $L$.

Eine zweite \NS $(Q',\pi')$ auf $\xi$ hei"st "aquivalent zu $(Q,\pi)$, falls
es einen $S^1$-Prinzipalb"undelisomorphismus $f:Q'\to Q$ gibt, so da"s
$\pi \circ f=\pi'$ gilt.

{\it Beachte:} Es kann vorkommen, da"s zwei $N$-Strukturen $(Q,\pi)$ und $(Q',\pi')$ auf
$\xi$ nicht "aquivalent sind, obwohl die $S^1$-B"undel $Q$ und $Q'$ isomorph
sind (d.h.~gleiche erste Chernklasse besitzen).

Sei $\phi:B\to BU(n)$ die bis auf Homotopie wohlbestimmte klassifizierende
Abbildung des $U(n)$-b"undels $\xi$. Die Determinantenabbildung $U(n)\to S^1$
induziert eine Abbildung $\det:BU(n)\to BS^1$ der klassifizierenden
R"aume; die zur Bildung der $N$-ten Tensorpotenz eines komplexen Linienb"undels geh"orende
Abbildung $\lambda_N:S^1\to S^1$ induziert eine Faserung $\mu_N:BS^1\to BS^1$
mit Faserhomotopietyp $K(\ZN,1)$ (lange exakte Homotopiesequenzen zu $\lambda_N$,
$\mu_N$ und $ES^1\to BS^1$ betrachten). Der Homotopietyp der betrachteten R"aume
und die Homotopieklasse der Abbildungen $\det$ und $\mu_N$ sind dabei eindeutig
festgelegt. "Aquivalenzklassen von \NSS entsprechen
dann Faserhomotopieklassen von Liftungen der Abbildung $\det\circ\,\phi:B\to BS^1$
bez"uglich der Faserung $\mu_N:BS^1 \to BS^1$.

Die Kohomologiesequenz zur Koeffizientensequenz $0\mapr{}\Z\mapr{\cdot N}\Z\mapr{}
\ZN\mapr{}0$ liefert die nat"urliche Transformation $\c1:H^2(\,\, . \,\, ,\Z)\to
H^2(\,\, . \,\, ,\ZN)$, welche von einer ebenfalls mit $\c1$ bezeichneten
Abbildung $K(2,\Z)\to K(2,\ZN)$ der klassifizierenden R"aume der Kohomologiefunktoren
induziert ist. Wir k"onnen die Homotopieklasse dieser Abbildung mit
dem Element $\c1(\gamma)\in H^2(BS^1,\ZN)$, der modulo-$N$-Reduktion der
ersten Chernklasse des universellen $S^1$-B"undels $\gamma$ "uber
$BS^1\homotop K(2,\Z)$, identifizieren. Die Faserung $\mu_N:BS^1\to BS^1$
ist dann die mit $\c1$ zur"uckgeholte (Faserprodukt) Wegeraumfaserung
$p_N:PK(\ZN,2)\to K(\ZN,2)$ mit Faserhomotopietyp $\Omega K(\ZN,2)\homotop
K(\ZN,1)$. Holt man die Faserung $\mu_N$ mit $\det$ weiter nach $BU(n)$ zur"uck,
erh"alt man eine Faserung $(\tilde\mu_N)_n:\widetilde{BU}(n)\to BU(n)$:
$$\begin{array}{ccccccc}\typeout{2 Pfeile einfuegen !!!}
 &\qquad    &\widetilde{BU}(n) &\mapr{}    &BS^1        &\mapr{}   &PK(\ZN,2)\\
 &\klein{\tilde\psi}&\mapd{(\tilde\mu_N)_n}&           &\mapd{\mu_N}&          &\mapd{p_N}\\
 &          &BU(n)             &\mapr{\det}&BS^1        &\mapr{\c1}&K(\ZN,2)\\
 &\klein{\phi}&                  &           &            &          &\\
B&          &                  &           &            &          &
\end{array}.$$
"Aquivalenzklassen von \NSS auf $\xi$ entsprechen den Faserhomotopieklassen
von Liftungen $\tilde\psi$ von $\phi$ bez"uglich dieser Faserung.

Hindernistheorie (vgl.~\cite{spa}, Kap.~8) besagt, da"s $\phi:B\to BU(n)$
genau dann eine Liftung besitzt, falls $\c1\circ \det \circ\,\phi$ nullhomotop
ist, d.h.~die Reduktion modulo $N$ der ersten Chernklasse von $\xi$
verschwindet. In diesem Fall k"onnen die Differenzen von Faserhomotopieklassen
von Liftungen $\tilde\psi:B\to\widetilde{BU}(n)$ von $\phi$ eineindeutig mit den
Elementen von $H^1(B,\ZN)$ identifiziert werden.

Ist das $U(n)$-B"undel $\xi$ durch "Ubergangsfunktionen $f_{ij}:U_i\cap U_j\to
U(n)$ gegeben, wobei ${\cal U}=\{U\}_{i \in I} $ eine $\xi$ trivialisierende
"Uberdeckung von $B$ ist, so hat man die folgende Beschreibung von
"Aquivalenzklassen von \NSS. Sei $\{g_{ij}\}$ der durch $g_{ij}(x)=\det f_{ij}(x)$,
$x\in U_i\cap U_j$ definierte ${\cal U}$-Kozyklus des zu $\xi$ geh"orenden
Determinantenb"undels. Die "Ubergangsfunktionen $h_{ij}:U_i\cap U_j\to S^1$ einer
$N$-ten Wurzel des Determinatenb"undels gen"ugen den Gleichungen $h_{ij}^N=g_{ij}$,
sind also f"ur feste $i$ und $j$ schon bis auf eine $N$-te Einheitswurzel
eindeutig bestimmt. Sind zwei feste $N$-te Wurzeln mit "Ubergangsfunktionen
$\{h_{ij}'\}$ und $\{h_{ij}''\}$ gegeben, so ist $\{h_{ij}'/h_{ij}''\}$ ein
$1$-Kozykel mit Werten in den $N$-ten Einheitswurzeln, repr"asentiert daher
ein Element von $H^1(B,\ZN)$. Man sieht so auch sofort, da"s jedes Element von
$H^1(B,\ZN)$ als Differenz von $N$-ten Wurzeln auftritt.

\LN{1.1.4} Seien $E$ und $\widetilde E$ zwei komplexe \VB{} "uber einem endlichen
$CW$-Komplex $B$. Dann bestimmt die Wahl einer $N$-Struktur auf zweien der drei
B"undel $E$, $\widetilde E$ und $E\oplus\widetilde E$ eindeutig eine \NS auf dem dritten.
\SE

{\bf Beweis:} Man versehe die B"undel $E$ und $\widetilde E$ mit einer Hermiteschen
Metrik, so da"s man die zugeh"origen $U(n)$-Prinzipalb"undel betrachten kann.

F"ur die erste Chernklasse gilt $c_1(E\oplus\widetilde E)=
c_1(E)+c_1(\widetilde E)$. Lassen daher zwei der Vektorb"undel eine \NS zu, so
auch das dritte.

Seien $\{g_{ij}\}$ bzw.~$\{\widetilde g_{ij}\}$ Kozykel, die die
Determinantenb"undel von $E$ bzw.~$\widetilde E$ repr"asentieren. Wegen
$\det(E\oplus\widetilde E)=\det(E)\otimes\det(\widetilde E)$ repr"asentiert
$\{g_{ij}\cdot \widetilde g_{ij}\}$ das Determinantenb"undel von
$E\oplus\widetilde E$. Sind $\{h_{ij}'\}$, $\{h_{ij}''\}$ bzw.~$\{\widetilde
h_{ij}'\}$, $\{\widetilde h_{ij}''\}$ Kozykel mit ${h_{ij}'}^N={h_{ij}''}^N=g_{ij}$,
$\hbox{$\widetilde h_{ij}'$}^N=\hbox{$\widetilde h_{ij}''$}^N=\widetilde g_{ij}$, die
zwei $N$-te Wurzeln von $\det E$ bzw.~$\det\widetilde E$ repr"asentieren, so
repr"asentieren $\{h_{ij}'\cdot \widetilde h_{ij}'\}$ und
$\{h_{ij}''\cdot \widetilde h_{ij}''\}$ zwei $N$-te Wurzeln von $\det(E\oplus
\widetilde E)$. Aus der Beziehung
$$\frac{h_{ij}'}{h_{ij}''}\cdot\frac{\widetilde h_{ij}'}{\widetilde h_{ij}''}=
\frac{h_{ij}'\cdot\widetilde h_{ij}'}{h_{ij}''\cdot\widetilde h_{ij}''}$$
zwischen den $1$-Kozykeln $\{h_{ij}'/h_{ij}''\}$, $\{\widetilde h_{ij}'/\widetilde h_{ij}''\}$
und $\{h_{ij}'\cdot\widetilde h_{ij}'/(h_{ij}''\cdot\widetilde h_{ij}'')\}$
mit Werten in den $N$-ten Einheitswuzeln ergibt sich die Behauptung.
\BOX

Ein die $U$-Struktur $\Phi$ einer $U$-\MGF $X_n$ repr"asentieren\-des
$U(n+k)$-Prin\-zi\-pal\-b"un\-del l"a"st ge\-nau dann eine $N$-Struktur zu,
falls $c_1(X_n)$ durch $N$ teilbar ist. Wir machen daher folgende

\DN{1.1.5}\index{$N$-Mannigfaltigkeit}\index{$N$-Struktur} Eine {\bf $N$-Mannigfaltigkeit} ist eine $U$-\MGF, deren erste Chernklasse durch $N$
teilbar ist.\hfill\break
Eine {\bf $N$-Struktur} auf einer \NM ist eine "Aquivalenzklasse von $N$-Struk\-tur\-en auf einem
Prinzipalb"undel, das zu einem stabilen Tangentialb"undel geh"ort, das die
$U$-Struktur repr"asentiert.
Wenn die $N$-\MGF nicht zusammenh"angend ist, ist dabei auf jeder Komponente
eine \NS zu w"ahlen.

Die Definition ist unabh"angig von dem die $U$-Struktur $\Phi$
repr"asentierenden Prinzipalb"undel (vgl.~\cite{law}, S.~207).

Ist eine \NM einfach zusammenh"angend, so ist die \NS eindeutig bestimmt.


\LN{1.1.6} Das kartesische Produkt zweier \NMS mit \NS ist kanonisch eine
\NM mit \NS. Eine \NS auf einer \NM $X$ mit Rand $\partial X$ induziert
kanonisch eine \NS auf der \NM $\partial X$.\SE

{\bf Beweis:} F"ur die zugrunde liegenden $U$-Strukturen vgl.~\cite{cofl2}.
Wir betrachten hier nur noch die zus"atzliche \NS.

Das Tangentialb"undel des kartesischen Produktes $X\times Y$ zweier
stabil fastkomplexer \MGFS $X$, $Y$ ist verm"oge eines kanonisch gew"ahlten
Isomorphismus stabil isomorph zu $\pi_X^*(TX\oplus\epsilon_{\R}^k)\oplus
\pi_Y^*(TY\oplus\epsilon_{\R}^{k'})$. Hierbei sind $\pi_X$ und $\pi_Y$ die
Projektionen von $X\times Y$ auf die Faktoren. Holt man die \NSS mit
$\pi_X$ und $\pi_Y$ von $X$ und $Y$ auf das Produkt zur"uck, so erh"alt
man mit Lemma 1.1.4 eine eindeutige \NS auf $X\times Y$.

Das Tangentialb"undel von $\partial X$ ist mittels eines kanonisch gew"ahlten
Isomorphismus stabil isomorph zu $TX\oplus\epsilon_{\R}^k\mid\partial X$.
Die Einschr"ankung $(Q\mid \partial X,\pi\mid\partial X)$ des die \NS definierenden
$S^1$-Prinzipalb"undels $Q$ und der Projektionsabbildung $\pi$ definiert damit
eine \NS auf $\partial X$.\BOX

Zwei \NMS $X$ und $Y$ mit \NS hei"sen differenzierbar "aquivalent, falls
ein $U$-Diffeomorphismus $f:X\to Y$ existiert, so da"s die mit $f$ nach $X$
zur"uckgeholte \NS von $Y$ die von $X$ ist.

Sei $I=[0,1]$ mit der $U$-Struktur versehen, da"s f"ur jede geschlossene
$U$-\MGF $X$, die auf der Komponente $X\times\{0\}$ von $\partial (X\times I)$
induzierte $U$-Struktur die von $X$ ist. Die \NM $I$ besitzt dann eine eindeutige
\NS f"ur obige $U$-Struktur. F"ur eine geschlosse \NM $X$ mit \NS
bezeichnen wir mit $-X$ die auf die Komponente $X\times\{1\}$ von $\partial
(X\times I)$ induzierte \NS (zur induzierten $U$-Struktur).

\DN{1.1.7}\index{$\OUN$} Zwei geschlossene $n$-dimensionale \NMS $X_n$ und $Y_n$ mit \NS
hei"sen {\bf $N$-kobordant}, falls es eine (n+1)-dimensionale $N$-Man\-nig\-fal\-tig\-keit
 $W_{n+1}$ mit \NS
gibt, so da"s $\partial W_{n+1}=X_n\Dcup-Y_n$ gilt (im Sinne von \NMS mit $N$-Strukturen).
$N$-Kobordanz ist eine "Aquivalenzrelation, die Menge der
"Aquivalenzklassen von $n$-dimensionalen \MGFS sei mit $\Omega_n^{U,N}$
bezeichnet. Da die "Aquivalenzrelation mit der disjunkten Vereinigung, der Bildung
des Negativen und der Produktbildung vertr"aglich ist, erh"alt
$\OUN=\bigoplus_{n=0}^{\infty}\Omega_n^{U,N}$ die Struktur eines graduiert
kommutativen Ringes, welchen wir als den {\bf $N$-Kobordismusring} bezeichnen.

$N$-kobordante \MGFS sind auch $U$-kobordant, da ein $N$-Kobordismus auch ein
$U$-Kobordismus ist. Damit h"angen Chernzahlen einer $N$-\MGF nur von der
$N$-Kobordismusklasse ab. Sie sind offensichtlich unabh"angig von der gew"ahlten
\NS. Das "`Vergessen"' der $N$-Struktur definiert einen nat"urlichen
Ringhomomorphismus\index{$\varrho^N_*$}
$$\varrho^N_*:\OUN\to\OU,$$
der weder injektiv noch surjektiv ist. Im n"achsten
Abschnitt werden wir sehen, da"s zumindest $\varrho_*^N\otimes\Q$ ein
Isomorphismus ist.

Die Konstruktion von $\OUN$ ist durch die analoge Definitonen von $Spin$-\MGFS
und $Spin$-Strukturen motiviert worden. Eine $U$-\MGF ist genau dann eine
$Spin$-\MGF, wenn sie eine $2$-\MGF ist, da die Stiefel-Whitney Klasse
$w_2$ die modulo-$2$-Reduktion der Chernklasse $c_1$ ist.
Die Kobordismustheorie $\Omega_*^{2,N}$
ist "aquivalent zum komplex-spin Kobordismusring $\Omega_*^{c-s}$, wie er in
\cite{stcs} definiert wird.


\section{Berechnung von $\OUNQ$ mit Pontrjagin-Thom Theorem}
Wir berechnen $\OUNQ$ mit Hilfe des Pontrjagin-Thom Theorems. Dazu zeigen
wir, wie sich die $N$-Kobordismustheorie in die in \cite{stong} definierten
$(B,f)$-Kobordismustheorien einordnet.

Die im letzten Abschnitt betrachteten Faserungen $(\tilde\mu_N)_n:
\widetilde{BU}(n)\to BU(n)$ sind durch die Faserung $\tilde\mu_N:
\widetilde{BU}\to BU$ "uber dem stabilen klassifizierenden Raum induziert:
\vskip1mm
$$\begin{array}{ccccccc}
\widetilde{BU}(n)&\mapr{k_n}&\widetilde{BU}(n+1)&\mapr{}&\widetilde{BU}&\mapr{}&BS^1\quad\cr
\mapd{(\tilde\mu_N)_n}& &\mapd{(\tilde\mu_N)_{n+1}}& &\mapd{\tilde\mu_N}& &\mapd{\mu_N}\quad\cr
BU(n)&\mapr{i_n}&BU(n+1)&\mapr{}&BU&\mapr{det}&BS^1\quad.
\end{array}$$

Setzen wir $B_{2n}=B_{2n+1}:=\widetilde{BU}(n)$ so, kommutiert das Diagramm
\vskip1mm
\begin{equation}\label{faserungen}
\begin{array}{ccccccccc}
\dots&\mapr{}&B_{2n}&\mapr{g_{2n}\,:=id}&B_{2n+1}&\mapr{g_{2n+1}\,:=k_n}&B_{2n+2}&\mapr{}&\dots\quad\cr
&&\mapd{f_{2n}}& &\mapd{f_{2n+1}}& &\mapd{f_{2n+2}}& & \cr
\dots&\mapr{}&BO(2n)&\mapr{j_{2n}}&BO(2n+1)&\mapr{j_{2n+1}}&BO(2n+2)&\mapr{}&\dots\quad.
\end{array}
\end{equation}
\vskip2mm
Dabei sei $j_k:BO(k)\to BO(k+1)$ die nat"urliche Abbildung, die
zum Addieren eines trivialen reellen Geradenb"undels geh"ort, und die Faserung
$f_k$ sei
definiert durch $f_{2n}:=h_n\circ(\tilde\mu_N)_n$, $f_{2n+1}:=j_{2n}\circ f_{2n}$,
falls $h_{n}:BU(n)\to BO(2n)$ die nat"urliche Abbildung ist, die zum "`Vergessen"'
der komplexen Struktur geh"ort.

In \cite{stong}, Kap.~II wird zu einer Familie von Faserungen $f_k:B_k\to BO(k)$,
f"ur die (\ref{faserungen}) kommutativ ist, die Kobordismustheorie $\Omega_*(B,f)$
definiert. Es gilt nun

\SN{1.2.1} Die Kobordismusgruppen $\Omega_n(B,f)$, die zu der oben definierten
Familie von Faserungen $f_k:B_k\to BO(k)$ geh"oren, entsprechen den
Kobordismusgruppen $\Omega_n^{U,N}$.
\SE

{\bf Beweis:} (Skizze)
a) \NMS $X$ mit \NS entsprechen eineindeutig stabilen "Aquivalenzklassen von
Homotopiepieklassen von Liftungen zu $f_l:B_l\to BO(l)$ der klassifizierenden
Abbildung des Normalenb"undels von $X_n$:\hfill\break
Eine $N$-Struktur auf dem stabilem Tangentialb"undel induziert eine
$N$-Struktur auf dem stabilen Normalenb"undel $\nu$\/: Bette $X$ differenzierbar
in den $\R^{2s}$ ein, $s$ gro"s. Dann ist
$$(TX\oplus\epsilon_{\R}^k)\oplus\nu\iso\epsilon_{\R}^{2s}.$$
Versehe das Tangentialb"undel $\epsilon_{\R}^{2s}$ des $\R^{2s}$ mit der
kanonischen $U$- und $N$-Struktur und schr"anke diese dann auf $X_n$ ein.
Das Normalenb"undel $\nu$ kann dann mit einer $U$-Struktur versehen werden
(s.~\cite{cofl2}, S.~21), und mit Lemma 1.1.4 erh"alt man die $N$-Struktur.
Der Zusammenhang zwischen $U$-Strukturen und Liftungen ist klar, den
Zusammenhang zwischen $N$-Strukturen und Liftungen haben wir im
letzten Abschnitt diskutiert.

b) Die in \cite{stong} definierte Menge von "Aquivalenzklassen $\Omega_n(B,f)$
entspricht den "Aquivalenzklassen $\Omega_n^{U,N}$.
\BOX

Man kann auf $\Omega_*(B,f)$ auch eine Ringstruktur definieren, die dann mit
der von $\OUN$ "ubereinstimmt. Wir werden dies aber nicht weiter ausf"uhren,
da wir uns nur f"ur $\OUNQ$ interessieren.

\SN{1.2.2} F"ur die Dimension des rationalen $N$-Kobordismusringes gilt\hfil\break
$\dim_{\Q}\Omega_n^{U,N}\otimes\Q=\cases{p(n/2),& falls $n$ gerade \cr 0,&falls
$n$ ungerade}$\,. Dabei bezeichnet $p(k)$ die Anzahl der Partitionen von $k$.
\SE

{\bf Beweis:} Die Kobordismusgruppen $\Omega_n^{U,N}\otimes\Q$ sind nach Satz 1.2.1 vom Typ
$\Omega_n(B,f)$ wie in \cite{stong}, Kap.~II definiert. Wir k"onnen daher
das verallgemeinerte Pontrjagin-Thom-Theorem (s.~\cite{stong}, S.~18) anwenden:
$$\Omega_n^{U,N}\iso\lim_{k\to\infty}\pi_{n+k}(T\widetilde{BU}_k,\infty).$$
Hierbei ist $T\widetilde{BU}_k$ der Thom-Raum des mit $f_k=(j_{k-1}\circ)h_{[k/2]}\circ
(\tilde\mu_N)_{[k/2]}:\widetilde{BU}([\frac{k}{2}])\to BO(k)$ nach
$\widetilde{BU}([\frac{k}{2}])$ zur"uckgeholten universellen Vektorb"undels
$\gamma_k$ "uber $BO(k)$. (Der Thom-Raum eines mit einer euklidischen Metrik
versehenen \VBS $E$ besteht aus allen Vektoren $v \in E$ mit
$\vert v\vert \leq 1$, wobei die Vektoren der L"ange $1$ zu einem Punkt
$\infty$ zusammengeschlagen sind. Er kann auch als die
Ein-Punkt-Kompaktifizierung von $E$ angesehen werden.)

Da der Raum $T\widetilde{BU}_k$ $(k-1)$-fach zusammenh"angend ist, ist der
mit $\Q$ tensorierte Hurewicz-Homomorphismus $\pi_{n+k}(T\widetilde{BU}_k,\infty)
\otimes\Q\to H_{n+k}(T\widetilde{BU}_k,\infty,\Z)\otimes\Q$ nach einem Theorem von
Serre (vgl.~\cite{mist}, S. \ \ ) f"ur $n+k<2k-1$, also $k>n+1$ ein Isomorphismus.
Der Thom-Isomorphismus ($f_k^*(\gamma_k)$ ist orientiert, da $f_k$
von einem komplexen \VB induziert) besagt $H_{n+k}(T\widetilde{BU}_k,\infty,\Z)
\iso H_n(\widetilde{BU}([\frac{k}{2}]),\Z)$. Nun gilt
$H_n(\widetilde{BU}([\frac{k}{2}]),\Q)\iso H_n(BU([\frac{k}{2}]),\Q)$,
wie man mit Hilfe der Leray-Spektralsequenz f"ur die Faserung
$(\tilde\mu_N)_{[k/2]}:\widetilde{BU}([\frac{k}{2}])\to BU([\frac{k}{2}])$ sieht,
da die Faser ein $K(\ZN,1)$ ist und $H_*(K(\ZN,1),\Q)\iso\Q$.

Insgesamt ergibt sich somit f"ur $k>n+1$
$$\pi_{n+k}(T\widetilde{BU}_k,\infty)\otimes \Q
\iso H_{n+k}(T\widetilde{BU}_k,\infty,\Z)\otimes \Q
\iso H_n(\widetilde{BU}([\kleinn{\frac{k}{2}}]),\Q)
\iso H_n(BU([\kleinn{\frac{k}{2}}]),\Q).$$

F"ur $2[\frac{k}{2}]\ge n$ ist $H^n(BU([\frac{k}{2}]),\Q)$ ein rationaler Vektorraum
der Dimension $p([\frac{k}{2}])$, falls $n$ gerade und der Dimension $0$, falls
$n$ ungerade. Nach dem universellen Koeffiziententheorem hat dann
$H_n(BU([\frac{k}{2}]),\Q)$ die gleiche Dimension, und es folgt die Behauptung.\BOX

Eine unmittelbare Folgerung ist

\SN{1.2.3} Die Abbildung $\varrho_*^N\otimes\Q:\OUNQ\to\OUQ$ ist ein
Ringisomorphismus.\SE

{\bf Beweis:} Nach Satz 1.2.2 und Satz 1.1.1 stimmen die Dimensionen von
$\OUNQ$ und $\OUQ$ "uberein.

Sei $X_4$, $X_6$, $X_8$, $\dots$ eine Basisfolge f"ur $\OSUQ$ und sei $X_2$
eine Riemannsche Fl"ache vom Geschlecht $g=N+1$, also $c_1(X_2)=-2N\cdot h$,
$h[X_2]=1$, und $s(X_2)\not=0$. Dann ist nach Satz 1.1.1 und 1.1.2
$X_2$, $X_4$, $X_6$, $\dots$ eine Basisfolge f"ur $\OUQ$, die aus lauter \NMS
besteht, d.h.~$\varrho_*^N\otimes\Q$ ist surjektiv und damit bijektiv. \BOX

Die rationale $N$-Kobordismusklasse einer \NM $X$ mit \NS ist somit
unabh"angig von der gew"ahlten \NS und schon durch die Chernzahlen von $X$
bestimmt. Umgekehrt kann jede Kombination von Chernzahlen durch ein Element
von $\OUNQ$ realisiert werden.

Im n"achsten Abschnitt werden wir den Typ einer \SO auf einer $N$-Man\-nig\-fal\-tig\-keit
$X$ definieren. Es wird sich zeigen, da"s dieser Typ unabh"angig von der
gew"ahlten \NS ist, falls nicht alle Chernzahlen von $X$ verschwinden. Wir
h"atten alternativ $\OUNQ$ einfach als den von $U$-\MGFS $X$ mit $N\mid c_1(X)$
erzeugten Ring {\it definieren} k"onnen, wobei zwei \MGFS kobordant sein sollen,
falls ihre Chernzahlen "ubereinstimmen. Der folgende Teil der Diplomarbeit ist
also unabh"angig von den Betrachtungen in den Abschnitten~1.1 und~1.2.

\section{"Aquivariante Kobordismustheorie und getwistet projektive B"undel:
die Ideale $I_*^N$, $I_*^{N,t}$, $J_*^N$ von $\OUNQ$ und
$I_*^{SU}$, $I_*^{SU,t}$, $J_*^{SU}$ von $\OSUQ$}
Obwohl, wie im letzten Abschnitt gezeigt, der Kobordismusring $\OUNQ$
isomorph zu $\OUQ$ ist, kann man interessante von $N$ abh"angige Ideale
definieren. Jede $U$-\MGF ist in $\OUQ$ zu einer \MGF mit effektiver
\SO kobordant, da die komplex projektiven R"aume solche zulassen und
$\OUQ$ erzeugen. Bei $N$-\MGFS ist dies anders, so mu"s f"ur $2$-\MGFS mit
effektiver \SO das $\hat A$-\GE verschwinden \cite{athis1}.
Wir definieren daher Ideale $I_*^N$ in $\OUNQ$, die von zusammenh"angenden \NMS mit
effektiver \SO erzeugt werden.
Die Definition des Typs einer \SO auf einer \NM als Element von $\ZN$
wird zu den Idealen $I_*^{N,t}$ mit $t\in\ZN$ f"uhren. Komplexe projektive B"undel
mit komplexer Faserdimension $n$ sind nur f"ur die $N$, die Teiler von
$n+1$ sind, \NMS. Wir werden daher allgemeiner das Ideal $J_*^N$ betrachten, das
von {\it getwistet\/}-projektiven B"undeln, die zugleich \NMS sind, erzeugt wird.
Analoge Definitionen sind f"ur $SU$-\MGFS m"oglich.

Sei $X_n$ eine zusammenh"angende \NM mit fest gew"ahlter \NS. Sei $\alpha:S^1\times X \to
X$ eine differenzierbare \SO, die die $U$-Struktur respektiert, d.h.~die
induzierte Abbildung $T\alpha:S^1\times TX\oplus\epsilon^k_{\R}\to TX\oplus\epsilon^k_{\R}$,
$ T\alpha(\lambda)=D\alpha(\lambda)\oplus\hbox{id}$ auf dem stabilen Tangentialb"undel
kommutiert f"ur alle $\lambda\in S^1$ mit einer die stabil fastkomplexe Struktur
definierenden Abbildung $J:TX\oplus\epsilon^k_{\R}\to TX\oplus\epsilon^k_{\R}$.

Wir wollen den Typ $t$ der
\SO $\alpha$ als Element von $\ZN$ definieren. Sei $K$ das komplexe
Determinantenb"undel des stabilen Tangentialb"undels und $L$ die festgew"ahlte
\NS, d.h.~ein komplexes Linienb"undel mit $L^N=K$. Betrachte die zu $K$ und $L$
geh"origen $S^1$-Prinzipalb"undel $P$ und $Q$.
Auf $Q$ operiert $\ZN\subset S^1$,
und der Quotient kann nach Definition der $N$-Struktur mit $P$
identifiziert werden.
Es bezeichne $\pi:Q\to P$ diese $N$-bl"attrige "Uberlagerung.
\epsfig{file=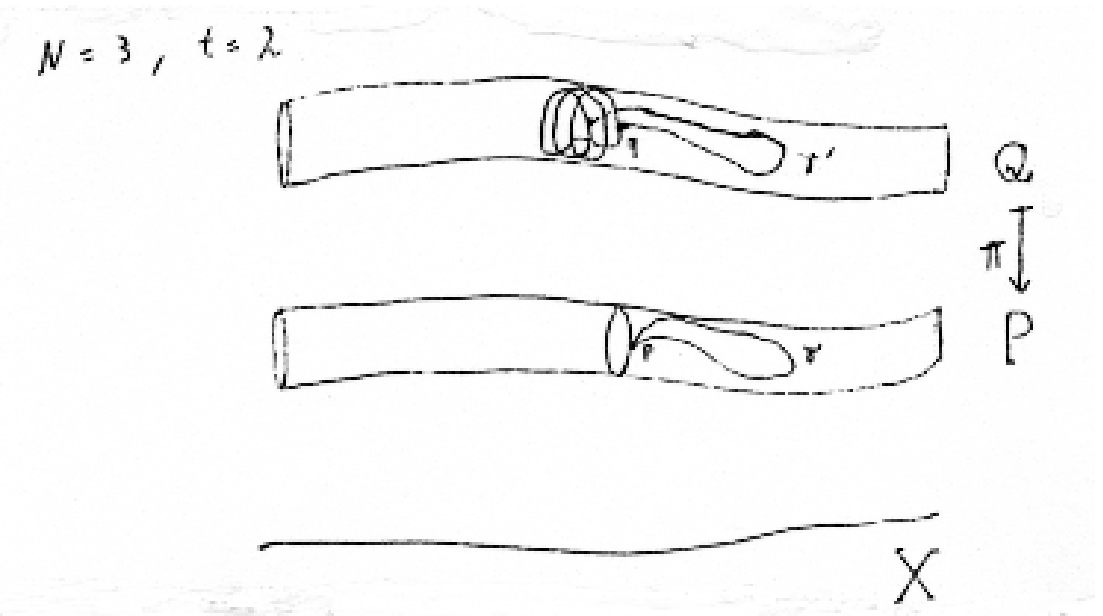,height=50mm}

Da die $S^1$ auf $X$ verm"oge $U$-Diffeomorphismen operiert, induziert sie eine
B"undelabbildung $S^1\times P\to P$, $(\lambda,p)\mapsto \det(T\alpha(\lambda))p$.
Betrachte f"ur einen beliebigen Punkt $p\in P$ den geschlossenen Weg
$\gamma: S^1\iso\R/\Z\to P$, $\lambda\mapsto \det(T\alpha(\lambda))p$, die Bahn von $p$
unter der \SO. F"ur jedes $q\in Q$ mit $\pi(q)=p$ gibt
es genau einen gelifteten Weg $\gamma':[0,1]\to Q$ mit
$\pi\circ\gamma'=\gamma$ und $\gamma'(0)=q$. Da $\pi\circ\gamma'(0)=\pi\circ\gamma'(1)=p$,
gibt es ein wohldefiniertes Element $t\in\ZN\subset S^1$ mit $t\gamma'(0)=
\gamma'(1)$. Dieses Element $t\in\ZN$ ist offenbar eine lokalkonstante Funktion
von $q \in Q$ und, da $X$ als \ZSHGD vorausgesetzt ist, also auch das
$S^1$-Prinzipalb"undel $Q$ \ZSHGD ist,
somit unabh"angig vom
gew"ahlten Bezugspunkt $q$. Damit ist die nachfolgende Definition sinnvoll:

\DN{1.3.1}\index{$t$} Der Typ einer \SO $\alpha:S^1\times X\to X$ auf einer zusammenh"angenden
\NM $X$ mit festgew"ahlter \NS ist das oben definierte Element $t\in\ZN$.

Der Typ $t$ einer \SO auf einer \NM h"angt im allgemeinen von der gew"ahlten \NS ab.

Beispiel: Die $S^1$ mit ihren $N$ verschiedenen $N$-Strukturen besitzt in
Abh"angigkeit der gew"ahlten \NS{} \SOS beliebigen Typs.
\typeout{Beispiel S1 erlaeutern}

Es wird sich zeigen, da"s dies nicht der Fall ist, falls die \SO Fixpunkte
besitzt.

Wir identifizieren im folgenden die $S^1$ mit den komplexen Zahlen
$\lambda=e^{2\pi i w}$, $w\in\R/\Z$, vom Betrag $1$.
Die Fixpunktmenge $X^{S^1}$ einer \SO $\alpha$ auf einer kompakten differenzierbaren
\MGF $X_n$ ist eine disjunkte Vereinigung von kompakten differenzierbaren
Untermannigfaltigkeiten verschiedener Dimensionen. Sei $Y\subset X^{S^1}$
eine Fixpunktkomponente mit Normalenb"undel $\nu$ und sei $p\in Y$.
Respektiert $\alpha$ eine stabil fastkomplexe Struktur, so wird der komplexe
$l:=(n+k)/2$ dimensionale Vektorraum $T_pX\oplus(\epsilon_{\R}^k)_p=T_pY\oplus\nu_p\oplus
(\epsilon_{\R}^k)_p$ zu einem komplexen $S^1$-Modul, der sich eindeutig in
Summanden $E^i_p$, $i\in\Z$, zerlegt, wobei ein Element $\lambda\in S^1$
auf $E_p^i$ duch Multiplikation mit $\lambda^i$ operiert. Bez"uglich einer
Basis von Eigenvektoren operiert $\lambda\in S^1$ auf dem stabilen Tangentialraum
verm"oge einer Diagonalmatrix $\hbox{diag}(\lambda^{m_1},\lambda^{m_2},\dots,\lambda^{m_l})$.
Dabei sind die als Drehzahlen bezeichneten $m_1$, $m_2$, $\dots$, $m_l$ ganze Zahlen,
und die Anzahl der Drehzahlen $m_{\nu}$ mit $m_{\nu}=i$ ist die komplexe
Dimension von $E^i_p$. Die Zerlegung von $T_pX\oplus(\epsilon_{\R}^k)_p$ l"a"st
sich zu einer Zerlegung von $TX\oplus\epsilon_{\R}^k\mid_Y$
in Eigenraumb"undel $E^i$ fortsetzen. Die Drehzahlen $m_1$, $m_2$, $\dots$, $m_l$
h"angen somit (bis auf die Reihenfolge) nur von der Fixpunktkomponente $Y$ ab.
Das Eigenraumb"undel $E^0$ ist gerade das stabile Tangentialb"undel $TY\oplus\epsilon_{\R}^k$
von $Y$, insbesondere ist $Y$ also selbst wieder stabil fastkomplex.

\SN{1.3.2} Sei $X_n$ eine zusammenh"angende \NM mit \SO $\alpha$. Falls
$[X_n]\not=0$ in $\OUNQ$, d.h.~nicht alle Chernzahlen verschwinden, ist der
Typ der \SO unabh"angig von der gew"ahlten \NS. Er ist gerade die Restklasse
modulo $N$ der Summe der Drehzahlen einer beliebigen Fixpunktkomponente.
\SE

{\bf Beweis:} Die Fixpunktmenge $X^{S^1}$ ist nicht leer, da sonst alle
Chernzahlen verschwinden w"urden (vgl.~\cite{hisl}, S.~15).

F"ur einen festen Punkt $p\in X^{S^1}$ operiert die Kreislinie auf dem stabilen
Tangentialraum $T_pX\oplus(\epsilon^k_{\R})_p$ verm"oge der Diagonalmatrix
$\hbox{diag}(\lambda^{m_1},\lambda^{m_2},\dots,\lambda^{m_l})$, wobei die
$m_{\nu}$ die Drehzahlen sind. Auf dem Vektorraum $K_p=\det(TX\oplus
\epsilon_{\R}^k)_p$, der Faser des Determinantenb"undels $K$ im Punkt $p$,
operiert $\lambda \in S^1$ dann verm"oge der Multiplikation mit
$\det(\hbox{diag}(\lambda^{m_1},\lambda^{m_2},\dots,\lambda^{m_l}))=\lambda^{m_1+m_2+\cdots+m_l}$.
Sei $P$ das $S^1$-Prinzipalb"undel zum Determinantenb"undel $K$.
Unabh"angig von der gew"ahlten \NS $(Q,\pi)$ ist die die \NS definierende
$N$-fache "Uberlagerung $\pi:Q\to P$ "uber dem Punkt $p\in X$ isomorph zur
$N$-fachen "Uberlagerung $\lambda_N:S^1\to S^1$, $\mu\mapsto\mu^N$.
Die Bahn eines Punktes $z\in P_p$ unter der \SO ist der Weg $\gamma:[0,1]\to P_p$,
$w\mapsto e^{2\pi i(m_1+m_2+\cdots+m_l)w}\cdot z$.
Er liftet zu einem Weg $\gamma':[0,1]\to Q_p$,
$w\mapsto e^{2\pi i\frac{1}{N}(m_1+m_2+\cdots+m_l)w}\cdot z'$, wobei $z'\in Q_p$
ein beliebiger Punkt mit $\pi(z')=z$ ist. Also gilt $\gamma'(1)=
e^{2\pi i\frac{1}{N}(m_1+m_2+\cdots+m_l)}\gamma'(0)$, d.h.~der Typ der \SO
$\alpha$ ist nach Definition $t\equiv m_1+m_2+\cdots+m_l \pmod{N}$.
\nopagebreak
\epsfig{file=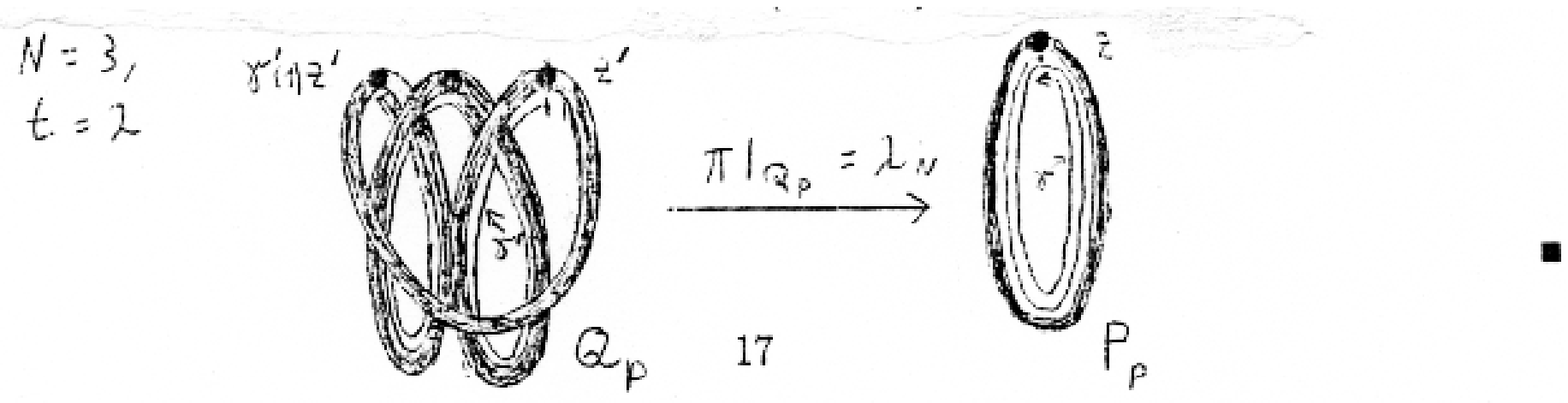,height=28.2mm,width=150mm}
\pagebreak[3]

Unter den Voraussetzungen des Satzes ist der Typ $t$ daher einerseits
unabh"angig von der gew"ahlten \NS, andererseits h"angt die Restklasse modulo $N$
der Summe der Drehzahlen nicht von der gew"ahlten Fixpunktkomponente ab.

Ist $\alpha:S^1\times X\to X$, $(\lambda,x)\mapsto\alpha(\lambda)x$ eine \SO
vom Typ $t$, so ist f"ur $k\in \Z$ die \SO $\alpha^k:S^1\times X\to X$, $(\lambda,x)\mapsto\alpha(\lambda^k)x$
vom Typ $kt$. L"a"st daher eine zusammenh"angende \NM eine effektive
\SO vom Typ $t$ zu, so l"a"st sie auch
eine vom Typ $s$ f"ur alle $s$ in der von $t$ erzeugten zyklischen Untergruppe
$\langle t\rangle$ von $\ZN$ zu, insbesondere also eine nichttriviale \SO
vom Typ $0$. Ist $\alpha$ eine \SO vom Typ $t$ auf $X$ und $\beta$ eine \SO
vom Typ $s$ auf $Y$, so ist $\alpha\times\beta:S^1\times X\times Y\to X\times Y$,
$(\lambda,x,y)\mapsto(\alpha(\lambda)x,\beta(\lambda)y)$ eine \SO vom Typ $t+s$.

Man k"onnte nun versuchen, die unbeschr"ankten Bordismusgruppen
$\Omega_*^{U,N}({\ST})$ von \SOS auf $N$-Man\-nig\-fal\-tig\-kei\-ten zu untersuchen
(d.h.~die Objekte sind $N$-Man\-nig\-fal\-tig\-kei\-ten mit einer die $U$-Struktur respektierenden \SO,
deren Isotropiegruppen beliebige Untergruppen der $S^1$ sein d"urfen
(vgl.~\cite{bos,cofl}). Es besteht die Zerlegung
$$\Omega_*^{U,N}({\ST})=\bigoplus_{t\in\ZN}\Omega_*^{U,N,t}({\ST})$$
bzgl.~des Typs der \SO. Aufgrund obiger Bemerkungen ist
$\Omega_*^{U,N,*}(\tilde S^1)$ eine $\N\times\ZN$ bigraduierte Algebra
"uber $\OUN$ (bzgl.~des ersten Grades (reelle Dimension)
graduiert kommutativ und bzgl.~zweiten Grades (Typ der $S^1$-Operation)
kommutativ). In dieser Allgemeinheit
werden wir die Kobordismusgruppen nicht weiter untersuchen (auch nicht "uber $\Q$),
sondern uns in der Diplomarbeit dem folgenden Problem zuwenden: Sei
$\epsilon_t:\Omega_*^{U,N,t}({\ST})\to\OUN$, $[X,\alpha]\mapsto[X]$
die Augmentation (Vergessen der $S^1$-Operation), ein $\OUN$-Modulhomomorphismus.
Die Einschr"ankung
von $\epsilon_t\otimes\Q$ auf $\ ^{\hbox{\scriptsize eff}}\Omega_*^{U,N,t}({\ST})\otimes\Q$, den rationalen
Kobordismusgruppen von \NMS mit effektiver \SO vom Typ $t$, liefert die
Ideale $I_*^{N,t}:=\hbox{Bild}(\epsilon_t\otimes\Q
\mid_{\ ^{\hbox{\tiny eff}}\Omega_*^{U,N,t}({\ST})\otimes\Q})$ von $\OUNQ$.
Die Berechnung dieser Ideale wird der Gegenstand der weiteren Betrachtungen
sein. Wir verwenden im folgenden aber die einfachere

\DN{1.3.3}\index{$I_*^N$}\index{$I_*^{N,t}$}
Bezeichne mit $I_*^N$ das Ideal in $\OUNQ$, das von zusammenh"angenden
\NMS $X$ mit effektiver \SO erzeugt wird. Bezeichne mit $I_*^{N,t}$, $t\in \ZN$,
das Ideal in $\OUNQ$, das von zusammenh"angenden \NMS mit effektiver \SO vom
Typ $t$ erzeugt wird.

Ist $[X]\not=0$ in $\OUNQ$, so ist wegen Satz 1.3.2 der Typ unabh"angig von der gew"ahlten
\NS. Daher mu"s man bei der Untersuchung von $I_*^N$ und $I_*^{N,t}$ auf
$X$ keine bestimmte $N$-Struktur festlegen.\pagebreak[3]

Nicht jede \NM $X$ mit $[X]\in I_*^N$($I_*^{N,t})$ l"a"st eine effektive \SO
(vom Typ $t$) zu, sondern dies bedeutet nur, da"s man eine \NM $Y$
(evtl.~aus mehre\-ren Komponenten bestehend) mit effektiver
\SO (vom Typ $t$) finden kann, so da"s $[mX]=[Y]$ f"ur eine nat"urliche
Zahl $m$ gilt.

Es ist $I_*^N=I_*^{N,0}$, denn l"a"st $X$ eine effektive \SO $\alpha$ vom Typ $t$
zu, so ist $\alpha^N$ eine effektive \SO vom Typ 0. Die umgekehrte Inklusion gilt nach
Definition. Da $I_*^{N,t}=I_*^{N,\hbox{\scriptsize ggT}(t,N)}$ (vgl.~Bemerkung weiter oben), werden
wir im folgenden nur die Ideale $I_*^{N,t}$ f"ur echte Teiler von $N$ betrachten.
F"ur $t\mid s$ gilt $I_*^{N,t}\subset I_*^{N,s}\subset I_*^{N,0}=I_*^N$.

Au"ser \NMS mit \SOS werden wir komplex projektive B"undel untersuchen.
F"ur kobordimustheoretische Betrachtungen ist die Klasse der projektiven B"undel
f"ur uns allerdings etwas zu klein, so da"s wir die allgemeinere Klasse von
getwistet-projektiven B"undeln betrachten werden.

Seien $E$ und $F$ komplexe \VB{}vom Rang $p$ und $q$ "uber einer $U$-\MGF $B$.
Die $U$-Struktur von $B$
und die komplexe Struktur von $E\oplus F$ liefern die
$U$-Struktur $\sigma^*(E\oplus F)\oplus \sigma^*TB$ auf $E\oplus F \mapr{\sigma} B$, sowie
--- nach Wahl einer Hermiteschen Metrik --- auf dem Scheibenb"undel
$D(E\oplus F)$ und dem Sph"arenb"undel $S(E\oplus F)=\partial D(E\oplus F)$.
Betrachte auf $E\oplus F$ die \SO $\alpha:S^1\times E\oplus F\to
E\oplus F$, $\alpha(\lambda,(p,e_p,f_p))=(p,\lambda e_p,\lambda^{-1} f_p)$. Diese
\SO ist eingeschr"ankt auf $S(E\oplus F)$
frei und respektiert die $U$-Struktur. Der Quotient $S(E\oplus F)/
\alpha$ ist als differenzierbare \MGF das komplex projektive B"undel
$\CP(E\oplus\FC)$, wobei $\FC$ das zu $F$ konjugierte
komplexe B"undel sei. "Uber dem komplex projektiven B"undel
$\CP(E\oplus\FC)$ haben wir die folgende
Sequenz komplexer \VB{}:

\begin{equation}\label{taut}
\begin{array}{ccccccccc}
0 &\mapr{}&S&\mapr{}& \pi^*(E\oplus\FC) &\mapr{}           & Q &\mapr{}& 0\\
  &       & &       & \mapd{}                              &         &   &       &  \\
  &       & &       & \CP(E\oplus\FC)                       &         &   &       &  \\
  &       & &       & \ \ \qquad\mapd{\pi} \qquad{}_{\sigma}&\quad E\oplus F&   &       &  \\
  &       & &       & B&                                   &   &       &.
\end{array}
\end{equation}\typeout{ein Pfeil !!!! }

Dabei sind $S$ das tautologische Linienb"undel und $Q$ das Quotientenb"undel.
Das B"undel $\CP^{\Delta}(E\oplus\FC)$ entlang der Fasern zu $\CP(E\oplus\FC)
{\buildrel \pi \over \longrightarrow}B$ ist isomorph
$Q\otimes S^*$. Aufgrund (\ref{taut}) folgt
$\CP^{\Delta}(E\oplus\FC)\oplus\epsilon_{\C}^1\iso Q\otimes S^*\oplus S\otimes S^*\iso
 S^*\otimes\pi^*E\oplus S^*\otimes \pi^*\FC$.

\DN{1.3.4}\index{$\TWIST$}
Das getwistet-projektive B"undel $\TWIST$ zu zwei komplexen
Vektorb"undeln $E$ und $F$ "uber einer $U$-\MGF $B$ ist als differenzierbare
\MGF das B"undel $\CP(E\oplus\FC)$, versehen mit der stabil fast komplexen Struktur
\begin{equation}\label{ttwist}
T\TWIST:=S^*\otimes\pi^*E\oplus\overline{S^*}\otimes\pi^*F\oplus\pi^*TB.
\end{equation}
\vskip 1mm
Die durch die "`getwistet"' stabil fastkomplexe Struktur induzierte
Orientierung von $\TWIST$ ist $(-1)^q$ mal die gew"ohnliche Orientierung
von $\CP(E\oplus\FC)$.
Den getwistet-projektiven Raum $\widetilde{\CP}(\C^p,\C^q)$ f"ur die \VB
$E=\C^p$, $F=\C^q$ "uber einem Punkt bezeichnen wir mit $\widetilde{\CP}_{p,q}$.
\index{$\widetilde{\CP}_{p,q}$}
Auf $\widetilde{\CP}_{p,q}$ operiert die Untergruppe $U(p)\times U(q)\subset
U(p+q)$ verm"oge $U$-Diffeomorphismen. Reduziert man die Strukturgruppe des
\VBS  $E \oplus \overline{F}$ "uber $B$ zu $U(p)\times U(q)$, so ist $\TWIST$ gerade
das $U$-Faserb"undel $B\times_{U(p)\times U(q)}\widetilde{\CP}_{p,q}$.

Der Kohomologiering von $\TWIST$ ist der von $\CP(E\oplus\FC)$:
\begin{equation}\label{htwist}
\quad\quad H^*(\TWIST,\Z)\iso H^*(B,\Z)[t]/\langle t^{p+q}+c_1(E\oplus\FC)t^{p+q-1}+\cdots+c_{p+q}(E\oplus\FC)\rangle,
\end{equation}
wobei $t=c_1(S^*)$ die erste Chernklasse des dualen B"undels des tautologischen
Linien\-b"undels $S$ ist.
F"ur die erste Chernklasse von $\TWIST$ erh"alt man aufgrund von (\ref{ttwist})
und (\ref{htwist})
\begin{equation}\label{ctwist}
c_1(\TWIST)=c_1(B)+c_1(E)+c_1(F)+(\rg E-\rg F)t.
\end{equation}
Die $U$-\MGF $\TWIST$ ist also genau dann eine $N$-\MGF, wenn
$N\mid c_1(B)+c_1(E)+c_1(F)$ und $\rg E\equiv \rg(F) \pmod{N}$ gilt.

Das durch $\alpha$ gegebene $S^1$-B"undel $S(E\oplus F)$ "uber $\TWIST$ hat
die erste Chernklasse $-t$. Ist $f$ eine dieses $S^1$-Prinzipalb"undel
klassifizierende Abbildung, so definiert $[\TWIST,f]$ ein Element in
$\OU(\CP_{\infty})$. Die Bordismusgruppe $\OU(\CP_{\infty})$ ist als
$\OU$-Modul isomorph zu $\Omega^U_{*+1}(\hbox{Frei})$, der Bordismusgruppe
von $U$-\MGF mit freier die $U$-Struktur respektierenden \SOS. Die
stabil fastkomplexe Struktur auf $\TWIST$ wurde gerade so gew"ahlt, da"s
$[\TWIST,f]$ unter obigem Isomorphismus der $S^1$-\MGF $[S(E\oplus F),\alpha]$
entspricht (vgl.~\cite{hata}, \S~3).

\DN{1.3.5}\index{$J_*^N$}
Bezeichne mit $J_*^N$ das Ideal von $\OUNQ$, das von
getwistet-projektiven B"undeln $\TWIST$  erzeugt wird, bei denen der Totalraum
und damit die Faser \NMS sind.

Bemerkung: Wir h"atten auch das Ideal $\tilde J_*^N$ betrachten k"onnen,
das von getwistet-projektiven B"undeln erzeugt wird, bei denen nur die Faser
eine $N$-\MGF ist. In Kapitel 2 wird $J_*^N=\ker \GN$ und
$J_*^N\subset\tilde J_*^N\subset\ker \GN$ gezeigt, so da"s $J_*^N=\tilde J_*^N$ folgt.

F"ur $SU$-\MGFS kann man zu 1.3.3 und 1.3.5 analoge Definitionen treffen.

Eine $SU$-\MGF $X$ ist eine \NM f"ur alle $N$.
Gilt $[X]\not=0$ in $\OSUQ$, so ist der Typ $t\in\ZN$ einer
\SO auf einer zusammenh"angenden \MGF $X$
--- aufgefa"st als $N$-\MGF{} --- nach Satz 1.3.2 die Restklasse modulo $N$ der Summe
der Drehzahlen einer Fixpunktkomponente. Da $N$ beliebig gew"ahlt werden kann,
ist die Summe der Drehzahlen f"ur alle Fixpunktkomponenten die gleiche. Diese
ganze Zahl $t$ sei der Typ der \SO auf einer $SU$-\MGF $X$, mit $[X]\not=0$
in $\OSUQ$.
Analog zu Def.~1.3.3 seien die Ideale $I_*^{SU}$ und $I_*^{SU,t}$, $t\in \Z$
definiert.\index{$I_*^{SU}$}\index{$I_*^{SU,t}$}
Es gelten die Beziehungen $I_*^{SU,t}\subset I_*^{SU,s}$ f"ur $t\mid s$, $s\not=0$,
(\SO $\alpha$ durch $\alpha^{s/t}$ ersetzen) sowie $I_*^{SU,t}\subset I_*^{SU},
t\in\Z$. Die Inklusion  $I_*^{SU}\subset I_*^{SU,0}$ ist im Gegensatz zu
der analogen bei \NMS nicht ohne weiteres ersichtlich.

Das Ideal $J_*^{SU}$\index{$J_*^{SU}$}
von $\OSUQ$ sei von getwistet-projektiven B"undeln
$\TWIST$ erzeugt, bei denen Totalraum und Faser $SU$-\MGFS sind.
F"ur $\tilde{J}_*^{SU}$ gilt eine zu $\tilde{J}_*^{N}$ analoge Bemerkung.
Wegen Lemma 1.4.8 aus dem n"achsten Abschnitt haben wir bei der Definition
von ${J}_*^{SU}$ nicht gefordert, da"s auch die Basis eine $SU$-\MGF sein
mu"s.

Da die rationale $SU$-($N$-)Kobordismusklasse einer $SU$-($N$-)\MGF durch
ihre Chernzahlen bestimmt ist, k"onnen wir $\OSUQ$ als Unterring von $\OUNQ$
auffassen. In $\OUNQ$ gelten dann die Inklusionen
$\langle I_*^{SU}\rangle\subset I_*^{N}$,
$\langle I_*^{SU,t}\rangle\subset I_*^{N,t\,\bmod{\,N}}$
und $\langle J_*^{SU}\rangle\subset J_*^{N}$.

Einen ersten Zusammenhang zwischen getwistet-projektiven B"undeln und \SOS
liefert

\SN{1.3.6} Sei $E$ ein komplexes Vektorb"undel "uber einer $U$-\MGF $B$,
das als direkten Summanden ein Geradenb"undel abspaltet, sei $F$ ein beliebiges
komplexes Vektorb"undel "uber $B$. Falls das
ge\-twis\-tet projektive B"undel $\TWIST$ eine $N$-\MGF bzw.~$SU$-\MGF ist,
l"a"st es \SOS von beliebigem Typ $t$ zu.
\SE

{\bf Beweis:} Die \VB $E$ und $F$ m"ogen den Rang $p$ und $q$ besitzen, wobei
$E=L\oplus K$  mit einem Geradenb"undel $L$.
Betrachte die \SO $\tilde\psi$ auf dem Sph"arenb"undel $S(E\oplus F)$:
$$\tilde\psi:S^1\times S(L\oplus K\oplus F)\to S(L\oplus K \oplus F).
\quad (\mu,(p,l_p,k_p,f_p))\mapsto (p,\mu^t\cdot l_p,k_p,f_p)$$
Sie kommutiert mit der $S^1$-Operation \hbox{$\alpha:
S^1\times S(L\oplus K\oplus F)\to S(L\oplus K \oplus F)$,}\break
$\quad (\lambda,(p,l_p,k_p,f_p))\mapsto (p,\lambda\cdot l_p,\lambda\cdot k_p,\lambda^{-1}\cdot f_p)$,
definiert daher eine \SO $\psi$ auf $\TWIST$, die die stabil fastkomplexe Struktur
$T\TWIST\iso \pi^*L\otimes S^*\oplus \pi^*K\otimes S^*\oplus \pi^* F\otimes
\overline{S^*}\oplus \pi^*TB$ respektiert.
Die \SO $\psi$ besitzt die Fixpunktmengen $\widetilde{
\CP}(K\oplus F)$ und $\CP(L)\iso B$ mit den komplexen Normalenb"undeln
$\pi^*L\otimes S^*\mid_{\widetilde{\CP}(K\oplus F)}$ und
$\pi^*K\otimes S^*\oplus \pi^*F\otimes\overline{S^*}\mid_{\widetilde{\CP}(L)}$.
Die jeweiligen Drehzahlen lauten (betrachte eine Faser $\widetilde{\CP}_{p,q}$)
$t$ f"ur den $\widetilde{\CP}(K\oplus F)$ sowie $(p-1)$ mal $(-t)$ und $q$ mal
$t$ f"ur den $\CP(L)$. Ist $\TWIST$ eine $N$-\MGF bzw.~$SU$-\MGF, so gilt
nach (\ref{ctwist}) $N\mid p-q$ bzw.~$p=q$. F"ur den Typ der \SO $\psi$ erh"alt man
also $t\equiv(p-1)\cdot(-t)+q\cdot t \pmod{p-q}$ bzw.~$t=(p-1)\cdot(-t)+q\cdot t$,
falls $p=q$.\BOX

In Kapitel 2 werden wir die in diesem Abschnitt definierten Ideale als Ideale
in den komplexen Kobordismusringen $\OUNC$ und $\OSUC$ auffassen, ohne dies
in den Bezeichnungen besonders zu vermerken.

\section{Konstruktion einer speziellen Basisfolge}
In diesem Abschnitt wird eine spezielle Basisfolge $W_1$, $W_2$, $W_3$, $\dots$
von $\OUQ$ kon\-struiert. Dabei sind $W_2$, $W_3$, $W_4$, $\dots$ $SU$-\MGFS, und
es liegen $W_3$ und $W_4$ in $\ISU$ (und damit in $I_*^N$) und
$W_5$, $W_6$, $\dots$ in $\JSU$, $I_*^{SU,t}$ (und damit in $\ISU$, $J_*^N$,
$I_*^{N,t}$).\footnote{Genauer m"u"ste es hei"sen: Die rationale Kobordismusklasse
$[W_i]$ von $W_i$ liegt in den jeweiligen Idealen.}
Au"serdem betrachten wir noch die \NMS $\CP_{N-1}$ und
$\widetilde\CP_{N+1,1}$.

{\it Vorbemerkung:} Ab diesem Abschnitt (bis auf den Beweis von Lemma 1.4.4)
bedeutet der Index $n$ bei einer $N$-
($U$-, $SU$- )\MGF $X_n$ die komplexe Dimension der \MGF (falls es eine w"are),
also die halbe reelle Dimension. Da wir nur die rationalen Kobordismusringe
betrachten ($\OUNQ$, $\OUQ$ und $\OSUQ$), und diese von reell gerade
dimensionalen \MGFS erzeugt werden, kann ohne Einschr"ankung die komplexe
Dimension als ganz vorausgesetzt werden.

Wir stellen zu Beginn einige Ergebnisse zusammen, die wir in diesem Abschnitt
dann ohne weitere Erw"ahnung verwenden werden.

Der $n$-dimensionale komplex projektive Raum \quad $\CP_n$ \quad besitzt den
Kohomologiering\break
$H^*(\CP_n,\Z)\iso \Z[g]/\langle g^{n+1}\rangle$. Hierbei sei der Erzeuger
$g\in H^2(\CP_n,\Z)\iso \Z$ stets
die Poincar\'e-duale Kohomologieklasse eines
Hyperebenenschnittes, der durch die komplexe Struktur ebenso wie der $\CP_n$
selbst kanonisch orientiert ist.
Ist $H$ eine glatte Hyperfl"ache vom Grad $d$ im $\CP_n$,
z.B.~$H=\{(z_0:\dots:z_n)\mid z_0^d+\cdots+z_n^d=0\}$, so hat das komplexe
Normalenb"undel $\nu$ von $H$ im $\CP_n$ die Chernklasse $c(\nu)=1+d\cdot i^*(g)$,
wobei $i:H\to\CP_n$ die Inklusionsabbildung bezeichne.
Die totale Chernklasse von $H$ lautet daher
$$c(H)=i^*\bigl((1+g)^{n+1}\cdot (1+ d\, g)^{-1}\bigr).$$
Ist  $u \in H^{2n-2}(\CP_n,\Z)$ eine Kohomologieklasse vom Grad $2n-2$, so gilt
(vgl.~\cite{himod}, S.~37 f.)
$$i^*(u)[H]=i_*(1)\cdot u[\CP_n]=d\cdot g\cdot u[\CP_n].$$

Die Milnorzahlen von $1$- bis $4$-dimensionalen $U$-\MGFS sind
die folgenden Linear\-kombinationen von Chernzahlen (vgl.~\cite{mist}, S.~188 f.):
$$s(X_1)=c_1,\qquad s(X_2)=c_1^2-2c_2,\qquad s(X_3)=c_1^3-3c_1c_2+3c_3,$$
$$s(X_4)=c_1^4-4c_1^2c_2+2c_2^2+4c_1c_3-4c_4.\qquad\qquad\qquad$$

Das $\chi_y$-Geschlecht einer $U$-\MGF ist ein Polynom in $y$ mit ganzzahligen
Koeffizienten, da der $U$-Kobordismusring "uber $\Z$ von algebraischen \MGFS
erzeugt wird (vgl.~\cite{hikm}). F"ur $2$- bis $4$-dimensionale
$U$-\MGFS gelten die folgenden Formeln:
\begin{eqnarray}
\chi_y(X_2)&=&\frac{1}{12}(1+y)^2 c_1^2+\frac{1}{12}(1-10 y+y^2) c_2\cr
\noalign{\vskip2mm}
\chi_y(X_3)&=&\frac{1}{24}(1+y-y^2-y^3)c_1c_2+\frac{1}{2}(y^2-y)c_3,\cr
\noalign{\vskip2mm}
\chi_y(X_4)&=& -\frac{1}{720}(1+y)^4 c_1^4+ \frac{1}{180}(1+y)^4 c_1^2 c_2
+\frac{1}{720}(1-56y-114y^2-56y^3+y^4)c_1c_3\cr
\noalign{\vskip1.5mm}
&&\qquad+\frac{1}{240}(1+y)^4 c_2^2-\frac{1}{720}(1+124-474 y^2 +124 y^3+y^4) c_4.
\nonumber
\end{eqnarray}
Die Signatur erh"alt man aus dem $\chi_y$-Geschlecht f"ur $y=1$, das
$\hat A$-Geschlecht von $SU$-\MGFS f"ur $y=0$.

\vskip 2mm
Definiere die $U$-\MGF $W_1$ als den $\CP_1$. F"ur die totale Chernklasse von $W_1$ gilt
$c_1(W_1)=(1+g)^2=1+2g$, falls $g$ der Erzeuger von $H^*(W_1,\Z)$ ist.

Wir halten fest:

\LN{1.4.1}\index{$W_1$} Die Chernzahlen der $U$-\MGF $W_1$ lauten
$$ c_1[W_1]=2. $$
Die Milnorzahl von $W_1$ ist $s(W_1)=2$.
\SE

Die $U$-Kobordismusklasse von $W_1$ ist bis auf ein Vorzeichen dadurch
ausgezeichnet, da"s $[W_1]$ einer der beiden Erzeuger des $\Z$-Moduls
$\Omega_2^U\iso\Z$ ist.

Eine glatte Hyperfl"ache vom Grad $4$ im $\CP_3$ --- eine sogenannte $K3$-Fl"ache ---
sei die Basismannigfaltigkeit $W_2$. Die totale Chernklasse von $W_2$ ist
$$c(W_2)=i^*((1+g)^4/(1+4g))=i^*(1+6g^2),$$
wobei $i$ die Inklusionsabbildung von $W_2$ nach $\CP_3$ ist und $g$ der Erzeuger
von $H^*(\CP_3,\Z)$. Damit erhalten wir

\LN{1.4.2}\index{$W_2$} Die Chernzahlen der $SU$-\MGF $W_2$ lauten
$$c_1^2[W_2]=0,\quad c_2[W_2]=24. $$
Die Milnorzahl von $W_2$ ist $s(W_2)=-48$.
\SE

Die $SU$-Kobordismusklasse $[W_2]$ ist als einer der beiden Erzeuger von
$\Omega_4^{SU}\iso \Z$ ausgezeichnet: F"ur die
Signatur von $[X_2]\in\Omega_4^{SU}$ gilt $\sign(X_2)=-\frac{2}{3}c_2[X_2]$. Da
die Signatur einer 4-dimensionalen
Spin-\MGF nach dem Satz von Rohlin (vgl.~\cite{ara}) durch
$16$ teilbar ist, und eine $SU$-\MGF eine Spin-\MGF ist, folgt $24 \mid c_2[X_2]$.

Die Basismannigfaltigkeit $W_3$ sei der homogene Raum $G_2/SU(3)\iso S^6$.
Er l"a"st eine homogen fastkomplexe Struktur auf dem Tangentialb"undel zu
(s.~\cite{bohi} Teil I, S.~500). Da $G_2$ eine Liegruppe vom Rang $2$ ist, l"a"st
$W_3$ auch effektive \SOS zu.
Wegen $H^2(S^6,\Z)=H^4(S^6,\Z)=0$ ist $c_3$ die
einzige nichtverschwindende Chernklasse von $W_3$ bzgl.~dieser fastkomplexen
Struktur. Die h"ochste Chernklasse einer fastkomplexen Struktur auf dem
Tangentialb"undel einer \MGF ist aber gerade die Eulerklasse.
Es gilt daher

\LN{1.4.3}\index{$W_3$} Die Chernzahlen von $W_3$ lauten
$$c_1^3[W_3]=0,\quad c_1c_2[W_3]=0,\quad c_3[W_3]=2. $$
Die Milnorzahl von $W_3$ ist $s(W_3)=6$. Die $SU$-\MGF $W_3$ ist in $\ISU$
(und damit in $I_*^N$) enthalten.
\SE

Auch $[W_3]\in \Omega_6^{SU}$ ist wieder bis auf das Vorzeichen ausgezeichnet, denn
$\Omega_6^{SU}\iso\Z$, und aus $[X_3]\in\Omega_6^{SU}$ folgt
$\chi_y(X_3)=\frac{-y+y^2}{2}c_3[X_3]$, d.h.~$2\mid c_3[X_3]$.

%
%
%
\LN{1.4.4}\index{$W_4$} Es existiert ein Element $[W_4] \in \OSUQ$ mit den
folgenden Eigenschaften:
Die Chernzahlen von $[W_4]$ lauten
$$c_1^4[W_4]=0,\quad c_1^2c_2[W_4]=0,\quad c_2^2[W_4]=2,\quad c_1c_3[W_4]=0,\quad c_4[W_4]=6. $$
Die Milnorzahl von $[W_4]$ ist $s(W_4)=-20$, und $[W_4]$ ist in $\ISU$
(und damit in $I_*^N$) enthalten.
\SE

Obwohl nicht sicher ist, da"s es eine \MGF $W_4$ mit den verlangten
Eigenschaften gibt (es gibt Probleme mit der $2$-Torsion, s.~Beweis), werden
wir "`so tun als ob"', da in unseren Betrachtungen nur die rationale
Kobordismusklasse eine Rolle spielt. Die Normierung der Chernzahlen
wurde gerade so gew"ahlt, da"s jedes Element $[X_4]\in \Omega_8^{SU}$ mit
$[X_4]\in I_*^{SU}$ sicher ein ganzzahliges Vielfaches von $[W_4]$ ist:
Nach Satz 2.4.7 verschwindet das $\hat A$-\GE von $X_4$, woraus $c_4=3c_2^2$ und damit
$\sign(X_4)=c_2^2$ f"ur die Signatur folgt. Nach \cite{cofl}, S.~70 (19.5) ist die Signatur
einer komplex 4-dimensionalen $SU$-\MGF gerade. Mit einem $k\in\Z$ gilt also
$c_2^2=2k$, $c_4=6k$, d.h.~$[X_4]=k\cdot[W_4]$. Es gibt auch $SU$-\MGFS
mit den Chernzahlen von $[W_4]$, z.B.~die Quadrik im $\CP_5$ mit einer
geeigneten stabil fastkomplexen Struktur, nur ist nicht klar, ob
eine diese Struktur respektierende \SO existiert.

{\bf Beweis von Lemma 1.4.4:} Wir werden $W_4$ mit Hilfe "aquivarianter Kobordismustheorie
konstruieren, wozu wir die folgenden geometrischen "Uberlegungen ben"otigen:

Sei $X$ eine $U$-\MGF mit semifreier die $U$-Struktur respektierender \SO.
Seien $F_j$ die Fixpunktkomponenten und $v(F_j)=v^+_j\oplus v^-_j$ die zugeh"origen
Normalenb"undel, wobei $\lambda \in S^1\subset \C^*$ auf den komplexen B"undeln
$v_j^+$ und $v_j^-$ durch Multiplikation mit $\lambda$ bzw.~$\lambda^{-1}$
operiere. Die $S^1$ operiert dann auf dem mit der induzierten $U$-Struktur
versehenen Scheibenb"undel $D(v_j^+\oplus v_j^-)$ und frei auf dem Sph"arenb"undel
$S(v_i^+\oplus v_i^-)$. Der Quotient $S(v^+_j\oplus v^-_j)/S^1$ ist das getwistet
projektive B"undel $\widetilde{\CP}(v_j^+\oplus v_j^-)$ (s.~Abschnitt 1.3), das
zusammen mit der klassifizierenden Abbildung $f_j$ des $S^1$-B"undels
$S(v^+_j\oplus v_j^-)$ das Element $[\widetilde{\CP}(v_j^+\oplus v_j^-),f_j]$
in $\OU(\CP_{\infty})$ ergibt. Da $\bigcup_j S(v_j^+\oplus v_j^-)$ von einer
$U$-Man\-nig\-fal\-tig\-keit mit freier \SO berandet wird, n"amlich dem Komplement
der Schei\-ben\-b"un\-del
$\bigcup_j D(v_j^+\oplus v_j^-)$ in $X$, gilt wegen des Isomorphismuses
$\Omega^U_{*+1}(\hbox{Frei})\iso \OU(\CP_{\infty})$ die Beziehung
\begin{equation}\label{nullbordant}
\sum_{j}[\widetilde{\CP}(v_j^+\oplus v_j^-),f_j]=0\quad\hbox{in $\OU(\CP_{\infty})$.}
\end{equation}
Gibt man sich umgekehrt $U$-\MGFS $F_j$ und komplexe \VB $v_j^+$, $v_j^-$ vor,
so ist (\ref{nullbordant}) hinreichend daf"ur, da"s sich die Vereinigung der
Scheibenb"undel $D(v_j^+\oplus v_j^-)$ zu einer geschlossenen $U$-\MGF $X$
vervollst"andigen l"a"st, die eine mit der $U$-Struktur kompatible semifreie \SO
mit Fixpunktkomponenten $F_j$ und vorgegebener Operation in den Normalenb"undeln
besitzt. Auf diese Weise werden wir nun die Basismannigfaltigkeit $W_4$
konstruieren.

Dazu w"ahlen wir als $4$-dimensionale Fixpunktkomponente $F_1$ den $\CP_2$ mit dem
Normalenb"undel $v(F_1)=v_1^+\oplus v_1^-=L_1\oplus L_2$, wobei $L_1$ und $L_2$
komplexe Geradenb"undel mit den Chernklassen $c_1(L_1)=-g$ und $c_1(L_2)=-2g$
seien und $g$ der Erzeuger von $H^*(\CP_2,\Z)$ ist.
Weiterhin verwenden wir $3$mal eine
$0$-dimensionale Fixpunktkomponente $F_2$ mit Normalenb"undel $v_2=v_2^+\oplus
v_2^-$, wobei $v_2^+=v_2^-=\C^2$. Wegen (\ref{ctwist}) sind die getwistet
projektiven B"undel $\widetilde{\CP}(v_j^+\oplus v_j^-)$, $j=1$, $2$ $SU$-\MGFS,
so da"s die Summe
\begin{equation}\label{disver}
[\widetilde{\CP}(L_1\oplus L_2),f_1]+ 3\cdot [\widetilde{\CP}
(\C^2\oplus\C^2),f_2]
\end{equation}
ein Element von $\Omega_6^{SU}(\CP_{\infty})$ definiert (auf jeder Komponente
dabei eine $SU$-Struktur w"ahlen).

Ein Element $[X,f]\otimes 1\in \Omega_6^{SU}(\CP_{\infty})\otimes\Q$ ist
genau dann Null, wenn die verallgemeinerten Chernzahlen $c_3[X]$, $c_2d[X]$
und $d^3[X]$ verschwinden, wobei $d=f^*(g)$ ein nach $H^2(X,\Z)$ zur"uckgeholter
Erzeuger von $H^*(\CP_{\infty},\Z)\iso\Z[g]$ und $c_2$, $c_3$ die Chernklassen
der $SU$-\MGF $X$ seien (vgl.~\cite{cofl2}, S.~25). In diesem Falle ist $[X,f]\in\Omega_6^{SU}(\CP_{\infty})$
ein Torsionselement, d.h.~es existiert ein $n\in \N$, so da"s $n\cdot[X,f]=0$.
Schaut man sich die Atiyah-Hirzebruch Spektralsequenz zu $\Omega_6^{SU}
(\CP_{\infty})$ an, sieht man, da"s nur $2$-Torsion vorkommen kann.

Berechnen wir nun die verallgemeinerten Chernzahlen von (\ref{disver}):\hfill\break
Man hat
$H^*(\widetilde{\CP}(L_1\oplus L_2),\Z)\iso\Z[g,t]/\langle g^3,t^2+c_1(L_1\oplus
\overline{L_2})t+c_2(L_1\oplus\overline{L_2})\rangle$
f"ur die Kohomologie von $\widetilde{\CP}(L_1\oplus L_2)$,
also $t^2=-(-g+2g)t-(-g)(2g)=-gt+2g^2$. Damit folgt f"ur die Chernklasse
$c(\widetilde{\CP}(L_1\oplus L_2))=(1+g)^3(1-g+t)(1-2g-t)=(1+3g+3g^2)(1-3g)=
1-6g^2$. Da $f_1:\widetilde{\CP}(L_1\oplus L_2)\to\CP_{\infty}$ das
Sph"arenb"undel klassifiziert, gilt $d=f_1^*(g)=-t$. Somit ergibt sich
f"ur die verallgemeinerten Chernzahlen (beachte $[\widetilde{\CP}(L_1\oplus L_2)]=
(-1)\cdot [\CP(L_1\oplus \overline{L_2})]$):
\begin{equation}\label{ck1}
\quad\quad c_3[\widetilde{\CP}(L_1\oplus L_2)]=0,\qquad
c_2d[\widetilde{\CP}(L_1\oplus L_2)]=-6,\qquad
d^3[\widetilde{\CP}(L_1\oplus L_2)]=3.
\end{equation}
F"ur den $\widetilde{\CP}(\C^2\oplus \C^2)=\widetilde{\CP}_{2,2}$ gilt
$H^*(\widetilde{\CP}_{2,2},\Z)=\Z[t]/\langle t^4\rangle$, $c(\widetilde{\CP}_{2,2})
=(1+t)^2(1-t)^2=1-2t^2$, $d=f_2^*(g)=-t$ und $[\widetilde{\CP}_{2,2}]=[\CP_3]$.
Die verallgemeinerten Chernzahlen lauten also
\begin{equation}\label{ck2}
c_3[\widetilde{\CP}_{2,2}]=0,\qquad c_2d[\widetilde{\CP}_{2,2}]=2,\qquad d^3[\widetilde{\CP}_{2,2}]=-1.
\end{equation}
Aus (\ref{ck1}) und (\ref{ck2}) folgt das Verschwinden der verallgmeinerten
Chernzahlen von (\ref{disver}).

Es gibt also eine singul"are $SU$-\MGF $g:X_7\to\CP_{\infty}$ mit
Rand, so da"s
$$\partial (X_7,g)=n\cdot ((\widetilde{\CP}(L_1\oplus L_2),f_1)\cup 3\cdot(\widetilde{\CP}
(\C^2\oplus\C^2),f_2)). $$
Der Isomorphismus $\Theta:\OU(\CP_{\infty})\to\Omega_{*+1}^U(\hbox{Frei})$ ist
durch die folgende Konstruktion gegeben: Man ordnet einer singul"aren \MGF
$f:X\to\CP_{\infty}$ das durch $f$ induzierte $S^1$-Prinzipalb"undel
$P\mapr{\pi}X$ zu, wobei $P$ mit der der $U$-Struktur
\begin{equation}\label{ufrei}
TP=\pi^*TX\oplus\pi^*\gamma
\end{equation}
versehen wird. Hierbei ist $\gamma$ das zu $P$ assoziierte Geradenb"undel,
so da"s $\pi^*\gamma$ ein triviales B"undel ist. Mit Hilfe dieser Konstruktion
erh"alt man zu $[X_7,g]$ eine $8$-dimensionale $U$-\MGF $X_8$ mit freier \SO
und dem Rand
$$\partial X_8=n\cdot(S(L_1\oplus L_2)\cup 3\cdot S(\C^2\oplus\C^2)).$$
Da wegen (\ref{ufrei}) $c_1(X_8)=\pi^*c_1(X_7)=0$ ist, l"a"st $X_8$ ebenso wie
die Scheibenb"undel $D(L_1\oplus L_2)$ und $D(\C^2\oplus \C^2)$ eine
$SU$-Struktur zu. Diese stimmen auf dem Rand $\partial X_8$ "uberein:
Das Hindernis gegen eine Liftung einer Abbildung $f:X\to BU$ bez"uglich
der Faserung $BSU\to BU$ ist die zur"uckgeholte erste Chernklasse $f^*c_1$,
und die Menge der Faserhomotopieklassen von Liftungen wird durch die Menge
$[X,\Omega BS^1]\iso H^1(X,\Z)$ beschrieben. Nun ist $H^1(\partial X_8,\Z)=0$,
da die Komponenten von $\partial X_8$ aus $S^3$- bzw.~$S^7$-B"undeln mit
einfach zusammenh"angender Basis bestehen. Die \MGFS{}~\hbox{$n\cdot ( D(L_1\oplus L_2)
\cup 3\cdot D(\C^2\oplus \C^2))$} und $X_8$ lassen sich also entlang des Randes
zu einer geschlossenen $SU$-\MGF
$W_4'=n\cdot ( D(L_1\oplus L_2)\cup 3\cdot D(\C^2\oplus \C^2))\cup_{\partial
X_8}(-X_8)$ mit effektiver die $U$-Struktur respektierender semifreier \SO
zusammenkleben. Der Typ $t$ dieser \SO ist als Summe der Drehzahlen einer
Fixpunktkomponente $t=2-2=1-1=0$. Wir setzen $[W_4]=\frac{1}{n}[W_4']$.

Es bleiben noch die Chernzahlen von $[W_4]$ zu berechnen.\footnote{
\hbox{Die Chernzahlen lassen sich auch direkt aus der Gleichung $[W_4]=
[\widetilde{\CP}((\epsilon_{\C}^1\oplus L_1)\oplus L_2)]$} $+3\cdot[\widetilde{\CP}
(\C^3\oplus \C^2)]$ in $\Omega_8^{U}\otimes\Q$ berechnen.}
Die Involution
$T\in S^1$ hat wie die \SO auf $W_4'$ als Fixpunktkomponenten $n$ Kopien
von $\CP_2$ und $3n$ isolierte Fixpunkte. Aufgrund der Starrheit der Signatur
bei \SO und nach \cite{hiinv} gilt $\sign(W_4')=\sign(T,W_4')=\sign(W_4'^T\circ
W_4'^T)$. F"ur die Signatur des Selbstschnittes $W_4'^T\circ W_4'^T$ der
Fixpunktmenge $W_4'^T$ liefern nur die $\CP_2$ einen Beitrag. Er ist gerade
die Eulerzahl des Normalenb"undels: $\sign(W_4'^T\circ W_4'^T)=n\cdot e(L_1
\oplus L_2)[\CP_2]=2n\cdot g^2[\CP_2]=2n$. Nach Satz 2.4.7 ist $\hat A(W_4')=0$.
Da $c_2^2$ und $c_4$ die einzigen von Null verschiedenen Chernzahlen sind,
gilt $\sign(W_4)=\frac{1}{45}(14c_4+3c_2^2)=2$, $\hat A(W_4)=\frac{1}{45\cdot 128}
(-8c_4+24c_2^2)=0$ und somit  $c_2^2[W_4]=2$,
$c_4[W_4]=6$ sowie $s(W_4)=-20$ f"ur die Milnorzahl.\BOX

Die weiteren Basismannigfaltigkeiten $W_5$, $W_6$, $\dots$ werden getwistet-projektive
B"undel sein. Dazu ben"otigen wir zwei Lemmata "uber die Milnorzahlen
von solchen B"undeln "uber einer $2$- oder $3$-dimensionalen komplexen Basis.

\LN{1.4.5} Sei $B$ eine 2-dimensionale komplexe \MGF, seien $E$ und $F$ komplexe
\VB {}"uber $B$ mit $\rg E=p\ge 1$, $\rg F=q\ge 1$ und $V:=E\oplus\FC$. Dann gilt
f"ur die Milnorzahl der $U$-\MGF $\TWIST$, falls ihre Dimension $p+q+1$ ungerade
ist,
$$ s_{p+q+1}(\TWIST)=(-1)^q\bigl((p-q)(c_1^2(V)-c_2(V))-(p+q+1)c_1(V)(c_1(E)+c_1(F))\quad$$
$$\qquad\qquad\qquad\qquad\qquad\qquad\quad+{p+q+1\choose2}(c_1^2(E)-c_1^2(F)-2c_2(E)+2c_2(F))\bigr)[B].$$
\SE

{\bf Beweis:} Nach Definition von $\TWIST$ gilt
$$T\TWIST\oplus\epsilon_{\C}^1\iso S^*\otimes \pi^*E\oplus\overline{S^*}\otimes\pi^* F\oplus\pi^*TB.$$
Dabei sind $\pi$ die Projektion von $\TWIST$ auf die Basis und $S$ das
tautologische Linienb"undel (vgl.~Def.~1.3.4). Falls $c(E)=\prod_{i=1}^{p}(1+x_i)$, $c(F)=\prod_{j=1}^{q}(1+y_j)$
und $c(TB)=\prod_{k=1}^{2}(1+v_k)$ die Zerlegungen der Chernklassen in formale
Wurzeln sind, und $t=c_1(S^*)$ gesetzt wird, erh"alt man f"ur die totale Chernklasse
von $\TWIST$ das Element
$$c(\TWIST)=\prod_{i=1}^{p}(1+t+x_i)\cdot\prod_{j=1}^{q}(1-t+y_j)\cdot\prod_{k=1}^{2}(1+v_k)$$
in $H^*(\TWIST,\Z)\iso H^*(B,\Z)[t]/\langle t^{p+q}+c_1(V)t^{p+q-1}+\cdots+c_{p+q}(V)\rangle$
(siehe (\ref{htwist})). F"ur die Milnorzahl ergibt sich
\begin{equation} \label{NR}
 s_{p+q+1}(\TWIST)=\left(\sum_{i=1}^p(x_i+t)^{p+q+1}+\sum_{j=1}^q(y_j-t)^{p+q+1}+\sum_{k=1}^2v_k^{p+q+1}\right)[\TWIST].
\end{equation}
Da $\dim_{\C}(B)=2$, ist $c_k(V)=0$ f"ur $k\ge 3$, und man erh"alt die
Relationen
\begin{eqnarray} \label{*}
  t^{p+q}&=&-c_1(V)t^{p+q-1}-c_2(V)t^{p+q-2}\qquad\qquad\\
\noalign{\hbox{sowie}}
 t^{p+q+1}&=&-c_1(V)t^{p+q}-c_2(V)t^{p+q-1}\qquad\qquad\label{**}\nonumber\\
 &=&(c_1^2(V)-c_2(V))t^{p+q-1}.\end{eqnarray}
Nach der binomischen Formel gilt
\begin{eqnarray}  \label{NR1}
\sum_{i=1}^p(x_i+t)^{p+q+1}&=&\qquad\qquad\qquad\qquad\qquad\qquad\qquad\qquad\qquad
\end{eqnarray}
$$ \qquad\qquad\qquad\qquad\qquad
   p\cdot t^{p+q+1}+(p+q+1)c_1(E)t^{p+q}+{p+q+1\choose 2}
   (c_1^2(E)-2c_2(E))t^{p+q-1},$$
mit (\ref{*}) und (\ref{**})
$$ = p(c_1^2(V)-c_2(V))t^{p+q-1}-(p+q+1)c_1(E)c_1(V)t^{p+q-1}
+{p+q+1\choose 2}(c_1^2(E)-2c_2(E))t^{p+q-1}.$$
Analog ergibt sich wegen $p+q+1$ ungerade
\begin{eqnarray}\label{NR2}
\sum_{j=1}^q(y_j-t)^{p+q+1}&=&\qquad\qquad\qquad\qquad\qquad\qquad\qquad\qquad\qquad\\
\noalign{\vskip -8mm $$
-q(c_1^2(V)-c_2(V))t^{p+q-1}-(p+q+1)c_1(F)c_1(V)t^{p+q-1}-{p+q+1\choose 2}
(c_1^2(F)-2c_2(F))t^{p+q-1}.$$}\nonumber
\end{eqnarray}
\vskip -10mm
Weiter gilt $\sum_{k=1}^{2}v_k^{p+q+1}=0$, da $p,q\ge 1$.
F"ur $\alpha\in H^*(B)$ gilt $t^{p+q-1}\cdot\alpha[\TWIST]=(-1)^q\pi_*(t^{p+q-1}\cdot
\pi^*(\alpha))[B]=(-1)^q\alpha[B]$. Wenden wir dies auf (\ref{NR}) an, erh"alt man mit
(\ref{NR1}) und (\ref{NR2})
$$ s_{p+q+1}(\TWIST)=(-1)^q\bigl((p-q)(c_1^2(V)-c_2(V))-(p+q+1)c_1(V)(c_1(E)+c_1(F))\quad$$
$$\qquad\qquad\qquad\qquad\qquad\qquad\qquad+{p+q+1\choose2}(c_1^2(E)-c_1^2(F)-2c_2(E)+2c_2(F))\bigr)
[B].\hbox{\BOX}$$

\LN{1.4.6} Sei $B$ eine 3-dimensionale komplexe \MGF, seien $E$ und $F$
komplexe \VB {}"uber $B$ mit $\rg E=p\ge 1$, $\rg F=q\ge 1$ und $V:=E\oplus\FC$.
Dann gilt f"ur die Milnorzahl der $U$-\MGF $\TWIST$, falls ihre Dimension
$p+q+2$ gerade ist,
$$s(\TWIST)=(-1)^q\bigl((p+q)(-c_3(V)+2c_1(V)c_2(V)-c_1^3(V))+\qquad\qquad\qquad\qquad$$
$$(p+q+2)c_1(V)(c_1^2(V)-c_2(V))+{p+q+2\choose 2}c_1(V)(-c_1^2(V)+2c_2(V))+$$
$$\qquad\qquad\qquad\qquad{p+q+2\choose 3}(c_1^3(V)-3c_1(V)c_2(V)+3c_3(V))\bigr)[B].$$
\SE

{\bf Beweis:} Der Beweis verl"auft analog zum Beweis von Lemma 1.4.5.
Mit den Bezeichnungen wie dort hat man f"ur die Milnorzahl die Beziehung
\begin{eqnarray*}
s(\TWIST)&=&\left(\sum_{i=1}^p(x_i+t)^{p+q+2}+\sum_{j=1}^q(y_j-t)^{p+q+2}+\sum_{k=1}^3v_k^{p+q+2}\right)[\TWIST]\nonumber\\
&=&\left(\sum_{i=1}^{p+q}(z_i+t)^{p+q+2}\right)[\TWIST].
\end{eqnarray*}
Hierbei sei $\prod_{i=1}^{p+q}(1+z_i)=c(E)c(\FC)=c(V)$. Unter Verwendung der
Abk"urzungen $c_i=c_i(V)$, $i=1,2,3$ gelten die Relationen
\begin{eqnarray}
t^{p+q}&=&-c_1t^{p+q-1}-c_2t^{p+q-2}-c_3t^{p+q-3},     \label{Q*}\\
t^{p+q+1}&=&-c_1(-c_1t^{p+q-1}-c_2t^{p+q-2}-c_3t^{p+q-3})-c_2t^{p+q-1}-c_3t^{p+q-2} \label{Q**}\\
&=&(c_1^2-c_2)t^{p+q-1}+(c_1c_2-c_3)t^{p+q-2},\nonumber\\
t^{p+q+2}&=&-c_1(c_1^2-c_2)t^{p+q-1}+(c_1c_2-c_3)t^{p+q-2}\label{Q***}\\
&=&(-c_1^3+2c_1c_2-c_3)t^{p+q-1}.\nonumber
\end{eqnarray}
Mit (\ref{Q*}), (\ref{Q**}) und (\ref{Q***}) ergibt sich
\begin{eqnarray*}
\hbox{$\sum_{i=1}^{p+q}(z_i+t)^{p+q+2}$}&=&
\hbox{$\sum_{i=1}^{p+q}(t^{p+q+2}+{p+q+2\choose 1}z_it^{p+q+1}
+{p+q+2\choose 2}z_i^2t^{p+q}+{p+q+2\choose 3}z_i^3t^{p+q-1})$}\\
&=&\hbox{$(p+q)t^{p+q+2}+(p+q+2)c_1t^{p+q+1}+$}\\
\noalign{\vskip -1mm\hfill\hbox{$+{p+q+2\choose 2}(c_1^2-2c_2)t^{p+q}+{p+q+2\choose 3}(c_1^3-3c_1c_2+3c_3)t^{p+q-1}$}\vskip 1mm}
&=&(p+q)(-c_3+2c_1c_2-c_1^3)t^{p+q-1}+(p+q+2)c_1(c_1^2-c_2)t^{p+q-1}+\\
\noalign{\hfill\hbox{$+{p+q+2\choose 2}(c_1^2-2c_2)(-c_1)t^{p+q-1}+{p+q+2\choose 3}(c_1^3-3c_1c_2+3c_3)t^{p+q-1}.$}}
\end{eqnarray*}
\vskip -0.8cm
Es folgt die Behauptung wie bei Lemma 1.4.5.\BOX

{\bf Definition:} a) Sei $\nu$ das komplexe Normalenb"undel der Hyperfl"ache
$W_2$ vom Grad $4$ im $\CP_3$. F"ur $n\ge 2$ setze $E:=\underbrace{\epsilon_{\C}^1
\oplus\dots\oplus\epsilon_{\C}^1}_{n-2 \hbox{ mal}}\oplus\epsilon_{\C}^1\oplus\nu^2$ und $F:=
\underbrace{\epsilon_{\C}^1\oplus\dots\oplus\epsilon_{\C}^1}_{n-2 \hbox{ mal}}\oplus\nu^{-1}\oplus\nu^{-1}$.
Die Basismannigfaltigkeit $W_{2n+1}$ sei dann das getwistet-projektive B"undel
$\TWIST$.

b) Sei $K$ das komplexe Determinantenb"undel "uber dem $\CP_3$. F"ur $n\ge 2$
setze \hbox{$E:=\underbrace{\epsilon_{\C}^1\oplus\dots\oplus\epsilon_{\C}^1}_{n-1\hbox{ mal}}\oplus
K$} und $F:=\underbrace{\epsilon_{\C}^1\oplus\dots\oplus\epsilon_{\C}^1}_{n-1\hbox{ mal}}\oplus
K^{-2}$. Die Basismannigfaltigkeit $W_{2n+2}$ sei dann das getwistet-projektive B"undel
$\TWIST$.

\LN{1.4.7}\index{$W_i$, $i\ge 5$} Die oben definierten $U$-\MGFS $W_5,W_6,W_7,\dots$ haben die
folgenden von Null verschiedenen Milnorzahlen:
\begin{eqnarray}
 s(W_{2n+1})&=&(-1)^n\cdot 128\, n(2n+1),\cr
 s(W_{2n+2})&=&(-1)^n\cdot 192\, (n-1)(2n-3)(2n+3)\not=0 \hbox{\ f"ur\ }n\ge 2.\nonumber
\end{eqnarray}
Sie sind getwistet-projektive B"undel, bei denen sowohl Faser als auch Totalraum
$SU$-\MGFS sind, d.h.~sie liegen in $J_*^{SU}$ (und damit in $J_*^N$). Sie lassen
\SOS von beliebigem Typ $t\in\Z$ zu, d.h.~sie liegen in $I_*^{SU,t}$, $t\in\Z$
beliebig (und damit in $I_*^{SU}$, $I_*^{N,t}$ und $I_*^N$).
Die keinen Faktor $c_1$ enthaltenden Chernzahlen von
$W_5$ sind:
$$c_2c_3[W_5]=-256,\quad c_5[W_5]=0,$$
sowie von $W_6$:
$$c_2^3[W_6]=192,\quad c_2c_4[W_6]=192,\quad c_3^2[W_6]=192,\quad c_6[W_6]=0.$$
\SE

{\bf Beweis:} Zu a) \ Die \MGFS $W_{2n+1}$, $n\ge2$: Die Chernklasse des
Normalen\-b"undels~$\nu$ der Hyperfl"ache $W_2$ im $\CP_3$ ist $c(\nu)=i^*(1+4g)$,
wobei wieder $i$ die Inklusionsabbildung von $W_2$ nach $\CP_3$ und $g$ der
Erzeuger von $H^*(\CP_3,\Z)$ seien. F"ur die Chernklassen der B"undel $E$ und $F$ vom
Rang $n$ sowie f"ur $V=E\oplus \FC$ gilt daher
$$c(E)=i^*(1+8g),\qquad c(F)=i^*(1-8g+16g^2),\qquad c(V)=i^*(1+16g+80g^2).$$
Lemma 1.4.5 liefert f"ur die Milnorzahlen
$$s(W_{2n+1})=(-1)^n{2n+1\choose2}(2\cdot16)i^*(g^2)[W_2]=128\cdot n(2n+1).$$
Nach (\ref{ctwist}) ist $c_1(W_{2n+1})=c_1(W_2)+c_1(E)+c_1(F)+(\rg E-\rg F)\cdot t=0$,
d.h.~$W_{2n+1}$ ist eine $SU$-\MGF. Die Faser des B"undels $W_{2n+1}$ ist der
$\widetilde\CP_{n,n}$, also ebenfalls eine $SU$-\MGF. Nach Satz 1.3.6 l"a"st $W_{2n+1}$
Operationen der $S^1$ von beliebigem Typ zu, liegt also in den angegeben Idealen. \hfil\break
Die totale Chernklasse von $W_5$ ist nach (\ref{ttwist})
$$(1+t)(1+t+8\cdot i^*(g))(1-t-4\cdot i^*(g))^2\cdot i^*((1+g)^4/(1+4g)).$$
Unter Ber"ucksichtigung der Relation $t^4=-c_1(V)t^3-c_2(V)t^2=
-16\,i^*(g)t^3-80\,i^*(g^2)t^2$ (siehe (\ref{htwist})) ergeben sich die angegebenen
Chernzahlen.

Zu b) \ Die \MGFS $W_{2n+2}$, $n\ge2$: Die Chernklasse des
Determinantenb"undels $K$ "uber $\CP_3$ ist $c(K)=1+4g$,
wobei $g$ der Erzeuger von $H^*(\CP_3,\Z)$ sei.
F"ur die Chernklassen der B"undel $E$ und $F$ vom
Rang $n$ sowie f"ur $V=E\oplus \FC$ gilt daher
$$c(E)=1+4g,\qquad c(F)=1-8g,\qquad c(V)=1+\underbrace{12g}_{c_1:=\quad}+\underbrace{32g^2}_{c_2:=\quad}.$$
Lemma 1.4.6 liefert f"ur die Milnorzahlen
\begin{eqnarray*}
s(W_{2n+2})&=&(-1)^n\bigl(2n(2c_1c_2-c_1^3)+(2n+2)(c_1^3-c_1c_2)+\\
\noalign{\vskip -8mm $$\qquad\qquad\qquad\qquad\qquad\qquad{2n+2\choose 2}(2c_1c_2-c_1^3)+{2n+2\choose 3}(c_1^3-3c_1c_2)\bigr)[\CP_3]$$\vskip 2mm}
&=&(-1)^n(4g)^3\bigl(2n(12-27)+(2n+2)(27-6)+\\
\noalign{\vskip -8mm$$\qquad\qquad\qquad\qquad(2n+2)(2n+1)(12-27)/2+(2n+2)(2n+1)n(27-18)/3\bigr)[\CP_3]$$\vskip 2mm}
&=&(-1)^n192\,(n-1)(2n-3)(2n+3).
\end{eqnarray*}
Nach (\ref{ctwist}) ist $c_1(W_{2n+2})=c_1(\CP_3)+c_1(E)+c_1(F)+(\rg E-\rg F)\cdot t=0$,
d.h.~$W_{2n+2}$ ist eine $SU$-\MGF. Die Faser des B"undels $W_{2n+2}$ ist der
$\widetilde\CP_{n,n}$, also ebenfalls eine $SU$-\MGF. Nach Satz 1.3.6 l"a"st $W_{2n+2}$
Operationen der $S^1$ von beliebigem Typ zu, liegt also in den angegeben Idealen. \hfil\break
Die totale Chernklasse von $W_6$ ist nach (\ref{ttwist})
$$(1+t)(1-t)(1+t+4g)(1-t-8g)(1+g)^4.$$
Unter Ber"ucksichtigung der Relation $t^4=-c_1(V)t^3-c_2(V)t^2=
-12gt^3-32g^2t^2$ ergeben sich die angegebenen Chernzahlen.
\BOX

Die Basismannigfaltigkeit $W_5$ ist ein getwistet-projektives B"undel, bei dem
au"ser der Faser und dem Totalraum sogar die Basis des B"undels eine $SU$-\MGF
ist. Im Gegensatz dazu gilt

\LN{1.4.8} Sind bei einem getwistet-projektiven B"undel der komplexen
Dimension~$6$ Basis, Faser und Totalraum $SU$-\MGFS, so verschwinden
alle Chernzahlen.
\SE

{\bf Beweis:} Sei $X=\TWIST$ ein solches B"undel mit Basis $B$. Aufgrund von
Gleichung (\ref{ctwist}), $c_1(B)+c_1(E)+c_1(F)+(\rg E-\rg F)t=0$,
sind die folgenden $3$ F"alle m"oglich:
$$\hbox{a)\quad}\dimc B=1\quad\hbox{und}\quad \rg E=\rg F=3,$$
$$\hbox{b)\quad}\dimc B=3\quad\hbox{und}\quad \rg E=\rg F=2,$$
$$\hbox{c)\quad}\dimc B=5\quad\hbox{und}\quad \rg E=\rg F=1.$$
Weiter mu"s $c_1(E)=c_1(\FC)=:a$ gelten. Wir verwenden im folgenden wieder
die gleichen Bezeichnungen wie beim Beweis von Lemma 1.4.5.

Zu a) \ Wegen $\dimc B=1$ ist $c(E)=c(\FC)=1+a$ und $c(V)=1+2a$. Damit gilt
f"ur die Chernklasse von $X$ in $H^*(X,\Z)\iso H^*(B,\Z)[t]/\langle t^6+2at^5\rangle$
die Beziehung
\begin{eqnarray*}
c(X)&=&\prod_{i=1}^3(1+x_i+t)\cdot\prod_{j=1}^3(1+y_j-t) \\
&=&(1+a+3t+2at+3t^2+at^2+t^3)(1-a-3t+2at+3t^2-at^2-t^3) \\
&=& 1\underbrace{-3t^2-2at}_{=\,c_2}+\underbrace{3t^4+4at^3}_{=\,c_4}.
\end{eqnarray*}
Es folgt $c_6=c_4c_2=c_3^3=c_2^3=0$. Wir haben dabei die Beziehung
$t^6=-2at^5$ verwendet.

Zu b) \ Mit den Abk"urzungen $b:=c_2(E)$, $c:=c_2(\FC)$, $u:=c_2(B)$ und $v:=c_3(B)$
gilt $c(E)=1+a+b$, $c(\FC)=1+a+c$ und $c(V)=1+2a+b+c+a^2+a(b+c)$.
Da $H^*(X,\Z)\iso H^*(B,\Z)[t]/\langle t^4+2at^3+(b+c+a^2)t^2+a(b+c)t\rangle$,
also $t^4=-2at^3-(b+c+a^2)t^2-a(b+c)t$, hat man f"ur die Chernklasse von $X$
\begin{eqnarray*}
c(X)&=&\prod_{i=1}^2(1+x_i+t)\cdot\prod_{j=1}^2(1+y_j-t)\cdot(1+v+w)\\
&=&(1+a+2t+t^2+at+b)(1-a-2t+t^2+at+c)(1+v+w)\\
&=&1+\underbrace{-a^2-2at+b+c-2t^2+v}_{=\,c_2}+\underbrace{-a(b+c)-2bt+2ct+w}_{=\,c_3}+\\
\noalign{\vskip -5mm $$\qquad\qquad\qquad\qquad\qquad\qquad\qquad
+\underbrace{-2atv-2t^2v}_{=\,c_4}+\underbrace{2t^2w}_{=\,c_5}+\underbrace{(-a^2-b-c)t^2v-(b+c)atv}_{=\,c_6}.$$}
\end{eqnarray*}
\vskip -10mm
Es folgt $c_6=c_4c_2=c_3^3=c_2^3=0$.

Zu c) \ In $H^*(X,\Z)\iso H^*(B,\Z)[t]/\langle t^2+2at \rangle$ gilt
$$c(X)=(1+a+t)(1-a-t)c(B)=c(B),$$
denn $(1+a+t)(1-a-t)=1-2at-t^2=1$. Da $H^{12}(B,\Z)=0$, verschwinden die Chernzahlen.
\BOX

\index{$\CP_{N-1}$, $\widetilde\CP_{N+1,1}$}
In Abschnitt 2.4 ben"otigen wir f"ur $N\ge 2$ noch den projektiven Raum
$\CP_{N-1}$ und den getwistet-projektiven Raum $\widetilde\CP_{N+1,1}$, welcher diffeomorph
zum $\CP_{N+1}$ ist. Da $c(\CP_{N-1})=(1+g)^N$ sowie $c(\widetilde\CP_{N+1,1})
=(1+g)^{N+1}(1-g)$, sind beide R"aume \hbox{$N$-Man\-nig\-faltig\-keiten}.
Sie liegen nach Konstruktion
in $J_*^N$ und nach Satz 1.3.6 in $I_*^{N,t}$ f"ur $t\in\ZN$ beliebig.



%
%
%

\chapter{Elliptische Geschlechter und Modulkurven}

\section{Komplexe Geschlechter}

Wir stellen in diesem Abschnitt den notwendigen Kalk"ul f"ur das Arbeiten mit
Geschlechtern bereit.

Sei $\OUQ \cong \Q[\CP_1,\CP_2,\CP_3,\dots]$ der in Abschnitt 1.1 definierte
rationale Kobordismusring der stabil fastkomplexen \MGFS.

\DN{2.1.1} Ein komplexes Geschlecht $\varphi$\index{$\varphi$} ist ein graduierter
Al\-ge\-bren-Ho\-mo\-mor\-phis\-mus
von $\OUQ$ in eine graduierte kommutative $\Q$-Algebra $\Lambda$ mit Einselement,
wobei $\varphi(1)=1$.

Einem Geschlecht sind zugeordnet (vgl.~\cite{hihab,mist}):

1) Die Potenzreihe \index{$g(y)$}$$g(y)=\sum_{n=0}^{\infty}\frac{\varphi(\CP_n)}{n+1}y^{n+1}
\in \Lambda[[y]],$$ die auch Logarithmus des Geschlechts hei"st.

2) Die Potenzreihe \index{$Q(x)$}\index{$f(x)$}$Q(x)=\frac{x}{f(x)}=1+a_1x+a_2x^2+\dots \in \Lambda[[x]]$,
wobei $f(x)$ die Umkehrfunktion von $g(y)$ ist, d.h $f(g(y))=y$, $g(f(x))=x$.
Man bezeichnet $Q(x)$ als die zu $\varphi$ geh"orige charakteristische Potenzreihe.

3) Das formale Gruppengesetz \index{$F(u,v)$}$\qquad F(u,v)=g^{-1}(g(u)+g(v))$.

\index{$K_n(c_1,c_2,\dots,c_n)$}
4) Eine multiplikative Folge von Polynomen $K_0=1$, $K_1(c_1)$, $\dots$,
$K_n(c_1,c_2,\dots,c_n)$, $\cdots$ mit $K_n\in \Lambda[c_1,c_2,c_3,\dots,c_n]$,  wobei die
$K_n$ homogen vom Gewicht $n$ in den $c_1$, $\dots$, $c_n$ sind ($c_i$ mit dem
Gewicht $i$ versehen). Dabei hei"st eine Folge multiplikativ, wenn aus
 $(1+c_1z+c_2z^2+\dots)=(1+c'_1z+c'_2z^2+\dots)
(1+c''_1z+c''_2z^2+\dots)$ die Gleichung
$$ \sum_{i=0}^{\infty}K_i(c_1,\dots,c_i)z^i=
\sum_{i=0}^{\infty}K_i(c'_1,\dots,c'_i)z^i\cdot
\sum_{i=0}^{\infty}K_i(c''_1,\dots,c''_i)z^i$$ \nopagebreak
folgt. Aus der charakteristischen Potenzreihe bestimmen sich die Polynome
verm"oge $Q(x_1)\cdot\dots\cdot Q(x_n)=\sum_{i=0}^{n}K_i(c_1,\dots,c_i)+\sum_{i=n+1}^{\infty}K_i(c_1,\dots,c_n,0,\dots,0)$,
falls $c_i$ die $i$-te elementarsymmetrische Funktion in den $x_\nu$
ist. Umgekehrt erh"alt man mit $Q(x)=\sum_{i=0}^{\infty}K_i(x,0,\dots,0)$
wieder die charakteristische Potenzreihe.

Jedes der Objekte in 1) bis 4) bestimmt umgekehrt eindeutig ein Geschlecht
$\varphi$, falls die Koeffizienten aus $\Lambda$ geeignet graduiert sind. Das
Geschlecht einer $n$-dimensionalen \hbox{$U$-\MGF $X$} mit totaler Chernklasse $c(X)=
1+c_1+c_2+\dots+c_n \in H^*(X,\Z)$ ist $\varphi(X)=K_n(c_1,c_2,\dots,c_n)[X]
\in \Lambda$, also eine feste $\Lambda$-Linearkombination von Chernzahlen von $X$. Ist die charakteristische Potenzreihe eine gerade Potenzreihe,
so ist das Geschlecht schon f"ur orientierte \MGFS erkl"art, da man es dann
schon mit den Pontrjaginklassen alleine berechnen kann.
Anstelle von Geschlechtern mit Werten in einer $\Q$-Algebra k"onnen
wir auch Geschlechter von $\OUC$ in eine $\C$-Algebra betrachten, z.B. k"onnen wir
ein Geschlecht "uber $\Q$ mit $\C$ tensorieren. Der Kalk"ul bleibt der gleiche.
Auch ist ein Geschlecht schon auf $\OU$ definiert.

{\bf Beispiele:}
1)\index{$\varrho^U_{SO}$} Sei $\varrho^U_{SO}: \OUQ \to \OSOQ$,
$\varrho^U_{SO}([\CP_{2n}])=[\CP_{2n}]$,
$\varrho^U_{SO}([\CP_{2n+1}])=0$, das komplexe Geschlecht, das einer
$U$-Kobordismusklasse die durch die Abbildung $BU\to BSO$ induzierte
$SO$-Kobordismusklasse zuordnet.
"Uber $\varrho^U_{SO}$ l"a"st sich jedes andere $SO$-Geschlecht faktorisieren.

2) Sei \index{$\sign$}$\sign$: $\OUQ\to \Q[t]$, $t$ vom Gewicht 2, das Geschlecht zur Potenzreihe
$Q(x)=\frac{\sqrt{t}x}{tanh(\sqrt{t}x)}$ oder $g(y)=\sum_{n=0}^{\infty}
\frac{t^n}{n+1}y^{n+1}$. Setzen wir $t=1$, so ist $\sign(X)$ die Signatur,
d.h.~die Signatur der Schnittform auf der mittleren Kohomologie von $X$.
Die Unbestimmte $t$ haben wir eingef"ugt, um zu erreichen, da"s das Geschlecht
die Graduierung von $\OUQ$ respektiert. Analog kann man so mit jedem nicht
graduierten Geschlecht mit Werten in $\Q$ verfahren, indem man $Q(x) \in \Q[[x]]$
durch \hbox{$Q(sx)\in \Q[s][[x]]$}, grad $s=1$, ersetzt. Diese Vorgehensweise hat
z.B. den Vorteil, da"s man in der Aussage --- \em Die Signatur ist bis auf Skalierung
das einzige Geschlecht (f"ur $\Omega_*^{SO}$), welches multiplikativ in Faserb"undeln
ist \rm --- auf den Zusatz \em "`bis auf Skalierung"' \rm verzichten kann. Ist
im folgenden ein Geschlecht $\varphi:\OUQ\to\Q$ gegeben, werden wir --- falls
erforderlich --- annehmen, da"s es homogen geschrieben sei, d.h.~wir ersetzen
$\varphi$ durch $\varphi_{\rm hom}:\OUQ\to \Q[s]$, $[X_n]\mapsto \varphi
(X_n)\cdot s^n$, grad~$s=1$.

Man kann auch Geschlechter zu Potenzreihen $Q(x)$ betrachten, die einen
konstanten Term $a_0\not=1$ besitzen. Ist $a_0$ eine Einheit, so kann man zur normierten
Potenzreihe $a_0^{-1}Q(a_0x)$ "ubergehen, welche das gleiche \GE liefert.

\nopagebreak
Wir wollen nun komplexen projektiven Variet"aten Geschlechter von $\OUC$
in \nopagebreak $\C$-Algebren zuordnen. W"ahle dazu eine feste
Basisfolge $X_1$, $X_2$, $\dots$ von $\OUC$. Sei \index{$\CP^{a_1,\dots,a_n}$}$\CP^{a_1,\dots,a_n}$ der
gewichtet projektive Raum mit Gewichten $a_1$,$\dots$,$a_n$ $\in \N$,
d.h.~$(\C^n\setminus\{0\})/\sim$, wobei genau dann $(x_1,\dots,x_n)\sim
(y_1,\dots,y_n)$ ist, wenn es ein $\lambda \in \C^*$ gibt mit $y_i=\lambda
^{a_i}x_i$ f"ur alle $i=1,\dots,n$ (vgl.~\cite{dol}). Sei \index{$V$} $V$ eine (nicht
notwendig irreduzible) projektive Variet"at in $\CP^{a_1,\dots,a_n}$ mit
gewichtet homogenem Verschwindungsideal \index{$I(V)$}$I(V)\subset\C[x_1,\dots,x_n]$, wobei $x_i$ das
Gewicht $a_i$ haben m"oge. Die graduierte Koordinatenalgebra von $V$ ist dann
$K(V)\cong\C[x_1,\dots,x_n]/I(V)$, und $\pi:\C[x_1,\dots,x_n]\to K(V)$\index{$K(V)$} m"oge
die Projektionsabbildung sein.
Setzt man nun $a_1<a_2<\dots<a_n $ voraus, so ist die Abbildung
$$\lambda_{a_1,\dots,a_n}:\OUC \to \C[x_1,\dots,x_n],\quad X_i\mapsto\cases{
x_\nu, \mbox{\enspace falls\enspace} i=a_\nu \cr 0,\ \,\mbox{\enspace sonst}}$$
wohldefiniert. Das zur Variet"at $V$ geh"orende
Geschlecht $\varphi_V:\OUC\to K(V)$ ist dann definiert durch $\varphi_V:=
\pi\circ\lambda_{a_1,\dots,a_n}$. Das Geschlecht h"angt von der Auswahl der
Basisfolge und der Wahl des Koordinatensystems f"ur den $\CP^{a_1,\dots,a_n}$
ab. Der Kern von \index{$\varphi_V$}$\varphi_V$ ist das Ideal
$ \lambda_{a_1,\dots,a_n}^{-1}(I(V))$. Es wird sich sp"ater zeigen, da"s
das elliptische Geschlecht $\GN$ der Stufe $N$ sich in dieser Weise beschreiben
l"a"st, falls man $a_i=i$, $i=1,\dots,4$ setzt und f"ur die Variet"at $V$ die geeignet
in den $\CP^{1,2,3,4}$ abgebildete Modulkurve $\X1{N}$ nimmt.\hfill\break
Wir notieren dazu noch folgendes einfaches

\LN{2.1.2} Seien $V$ und $W$ projektive Variet"aten im $\CP^{a_1,\dots,a_n}$
mit festgew"ahltem Koordinatensystem.
Dann gilt bei fester Basisfolge $\ker \varphi_{V\cup W}=\ker \varphi_V \cap
\ker\varphi_W $.
\SE

{\bf Beweis:} Mit obigen Bezeichnungen gilt $\qquad\ker\varphi_{V\cup W}=
 \lambda_{a_1,\dots,a_n}^{-1}(I(V\cup W))= \lambda_{a_1,\dots,a_n}^{-1}(I(V)\cap I(W))=
 \lambda_{a_1,\dots,a_n}^{-1}(I(V))\cap\lambda_{a_1,\dots,a_n}^{-1}(I(W))=
\ker\varphi_{V}\cap\ker\varphi_{W}$.
\BOX

Da $\OUNQ\cong\OUQ$ und die Chernzahlen einer $N$-\MGF deren rationale
Kobordismusklasse sowohl in $\OUNQ$ als auch $\OUQ$ festlegen, k"onnen
wir Geschlechter zu $\OUNQ$ und $\OUQ$ miteinander identifizieren.

Die charakteristische Potenzreihe $Q(x)$ f"ur ein Geschlecht $\varphi$ ist
durch seine Werte auf $\OSUQ\subset\OUQ$
nur bis auf einen Faktor $e^{\alpha x}$ bestimmt: F"ur eine
$SU$-\MGF $M_n$ gilt
\vskip-10mm
\begin{equation}\label{qsu}
\varphi(M_n)=\prod_{i=1}^{n}Q(x_i)[M_n]=e^{\alpha c_1(M_n)}\prod_{i=1}^{n}Q(x_i)[M_n]=
\prod_{i=1}^{n}e^{\alpha x_i}Q(x_i)[M_n],
\end{equation}
\vskip-3mm
da $c_1(M_n)=0$.
Man gelangt zu einem \GE f"ur $\OUQ$, falls der Wert von $\varphi $ auf der Basismannigfaltigkeit
$X_1$ vorgegeben wird. Dadurch ist dann die obige Konstante $\alpha$ aus dem
Exponentialfaktor eindeutig festgelegt.

\nopagebreak
Der Kern eines Geschlechtes (das, falls notwendig, wie im 2.~Beispiel erkl"art, homogen
geschrieben wird) ist ein homogenes Ideal (d.h.~er wird von homogenen Elementen
\nopagebreak erzeugt).
Umgekehrt definiert ein homogenes
Ideal ein Geschlecht, n"amlich \nopagebreak die Projektion \nopagebreak auf die Quotientenalgebra.
\pagebreak[3]


\section{Das universelle komplexe elliptische Geschlecht}

Das universelle elliptische Geschlecht wird als das Geschlecht definiert,
welches zu einer Potenzreihe geh"ort, die die L"osung einer bestimmten
Differentialgleichung $2$-ter Ordnung ist. Dazu folgendes

\LN{2.2.1} Sei $Q(x)=\frac{x}{f(x)}=1+a_1x+a_2x^2+\dots $ und $h(x):=\frac
{f'(x)}{f(x)}$. Sei weiter \index{$S(y)$}$S(y)=y^4+q_1 y^3 +q_2 y^2 +q_3 y +q_4 $ ein normiertes
Polynom 4-ten Grades in $y$ mit Koeffizienten $q_1$ bis $q_4$. Dann besitzt die
Differentialgleichung
\begin{equation}\label{dgl}
\index{$h(x)$}
\left(h'(x)\right)^2=S(h(x))
\end{equation}
eine eindeutige L"osung $h(x)\in\Q[q_1,q_2,q_3,q_4][[x]][x^{-1}]$. Diese
bestimmt wiederum eindeutig die Potenzreihe $Q(x)
\in \Q[q_1,q_2,q_3,q_4][[x]]$. Die Koeffizienten $a_n$ von
$Q(x)$ sind dabei homogene Polynome vom Gewicht $n$ in $q_1$ bis $q_4$, falls
$q_1$ bis $q_4$ die Gewichte $1$ bis $4$ zugeordnet werden.
\SE

{\bf Beweis:} Die Funktion $h(x)$ besitzt bei der Normierung $a_0=1$ eine
Reihenentwicklung
\begin{equation}
h(x)=\frac{1}{x}+c_1+c_2x+c_3x^2+\cdots\ .\nonumber
\end{equation}
Die Differentialgleichung (\ref{dgl}) hat dann die Gestalt
\begin{eqnarray*}
\frac{1}{x^4}(-1+c_2x^2+2c_3x^3+\cdots)^2&=&
\frac{1}{x^4}\bigl[(1+c_1x+c_2x^2+c_3x^3+\cdots)^4\\
&&\qquad +q_1x(1+c_1x+c_2x^2+c_3x^3+\cdots)^3\\
&&\qquad +q_2x^2(1+c_1x+c_2x^2+c_3x^3+\cdots)^2\\
&&\qquad +q_3x^3(1+c_1x+c_2x^2+c_3x^3+\cdots)\\
&&\qquad +q_4x^4\bigr].
\end{eqnarray*}
Koeffizientenvergleich liefert die Beziehungen
\begin{eqnarray*}
0&=&4c_1+q_1\\
-2c_2&=&4c_2+6c_1^2+3q_1c_1+q_2\\
-4c_3&=&4c_3+4c_1^3+12c_1c_2+q_1(3c_2+3c_1^2)+2q_2c_1+q_3\\
&\dots&\\
-2(n-1)c_n+\alpha_n(c_1,\dots,c_{n-1})&=&4c_n+\beta_n(c_1,\dots,c_{n-1})+q_1\,\gamma_{n-1}(c_1,\dots,c_{n-1})\\ \nopagebreak
&&\quad +q_2\,\delta_{n-2}(c_1,\dots,c_{n-2})+q_3\,\varepsilon_{n-3}(c_1,\dots,c_{n-3})\\
&&\quad +q_4\,\delta_{n,4}.
\end{eqnarray*}\pagebreak[3]
In der letzten Gleichung sind $\alpha_n$, $\beta_n$, $\gamma_{n-1}$, $\delta_{n-2}$
und $\varepsilon_{n-3}$ gewichtet homogene Polynome in den $c_i$ vom Gewicht
$n$, $n$, $n-1$, $n-2$ bzw.~$n-3$, und $\delta_{n,4}$ ist das Kronecker-Symbol.
Man kann also mit $c_1=-q_1/4$ beginnend
die Koeffizienten $c_n$ von $h(x)$ als homogene Polynome vom Gewicht $n$ in
den $q_i$ bestimmen.

Sei $\frac{1}{Q(x)}=1+b_1x+b_2x^2+b_3x^3+\cdots$ die zu $Q(x)$ reziproke
Potenzreihe. Die Koeffizienten $b_n$ gen"ugen den Rekursionsgleichungen
\begin{equation}\label{rin}
b_n=-a_n+\phi_n(a_1,\dots,a_{n-1}),
\end{equation}
wobei $\phi_n$ ein gewichtet homogenes Polynom in $a_1$ bis $a_{n-1}$ vom Gewicht $n$
ist. F"ur die Gleichung
$h(x)=\frac{f'(x)}{f(x)}=\frac{Q(x)}{x}\cdot\bigl(\frac{x}{Q(x)}\bigr)'$
erh"alt man
$$\frac{1}{x}(1+c_1x+c_2x^2+c_3x^3+\cdots)=\frac{1}{x}(1+a_1x+a_2x^2+a_3x^3+\cdots)
(1+2b_1x+3b_2x^2+4b_3x^3+\cdots).$$
Koeffizientenvergleich liefert
\begin{eqnarray*}
c_1&=&a_1+2b_1=-a_1\\
c_2&=&a_2+2a_1b_1+3b_2=-2a_2+a_1^2\\
&\dots&\\
c_n&=&a_n+(n+1)b_n+\psi_n(a_1,\dots,a_{n-1},b_1,\dots,b_{n-1}).
\end{eqnarray*}
Unter Verwendung der Rekursionsgleichungen (\ref{rin}) folgt
\begin{equation}
a_n=-\frac{1}{n}c_n+\pi_n(a_1,\dots,a_{n-1}).\nonumber
\end{equation}
Hierbei sind $\psi_n$ und $\pi_n$ gewichtet homogene Polynome. Induktiv
sieht man wieder, da"s $a_n$ ein gewichtet homogenes Polynom vom Gewicht $n$ in
den $c_i$ und damit auch in den $q_i$ ist.\BOX

\DN{2.2.2} Das universelle komplexe elliptische Geschlecht \index{$\GU$}$\GU:\OUQ \to
\Q[A,B,C,D]$ ist das komplexe Geschlecht zu der Potenzreihe $Q(x)$, die zur
L"osung der Differentialgleichung (\ref{dgl}) geh"ort. F"ur die Unbestimmten
$q_1$ bis $q_4$
werden die folgenden homogenen Polynome in $A$, $B$, $C$ und $D$ gesetzt, wobei
die Unbestimmten $A$ bis $D$ mit den Gewichten $1$ bis $4$ versehen seien
(Die Wahl der Polynome erkl"art sich aus dem n"achsten Satz):
\begin{equation}\label{Q1234}\index{$q_1$,$q_2$,$q_3$,$q_4$}
\qquad q_1=2 A,\qquad q_2=\frac{3}{2}A^2-\frac{1}{4}B,
\qquad q_3=\frac{1}{2}A^3-\frac{1}{4}AB+4 C,\hfill \end{equation}
$$\displaylines{\qquad\qquad q_4=\frac{1}{16}A^4-\frac{1}{16}A^2 B+2 A C+\frac{1}{64}B^2-2 D.} $$

Sind $q_1$ bis $q_4$ gegeben, kann man umgekehrt $A$ bis $D$ bestimmen:
\begin{equation}\label{ABCD}\index{$A$,$B$,$C$,$D$}
\qquad A=\frac{1}{2}q_1,\qquad B=\frac{3}{2}q_1^2-4 q_2, \qquad
C=\frac{1}{32}q_1^3-\frac{1}{8}q_1 q_2+\frac{1}{4}q_3,\hfill \end{equation}
$$\displaylines{\qquad\qquad D=\frac{3}{128}q_1^4-\frac{1}{8}q_1^2 q_2
+\frac{1}{8}q_1 q_3+\frac{1}{8} q_2^2-\frac{1}{2} q_4.}$$

Da die Koeffizienten von $Q(x)$ nach Lemma 2.2.1 homogen in den $q_1$ bis $q_4$
sind und damit auch in $A$, $B$, $C$ und $D$, ist der Wert von $\GU$ auf einer
komplex $n$-dimensionalen \MGF ein homogenes Polynom vom Gewicht $n$ in
$A$, $B$, $C$ und $D$.
F"ur die ersten Koeffizienten der Potenzreihe $Q(x)$ erh"alt man\footnote{Die
Berechnungen wurden mit Hilfe des Computerprogramms \sc REDUCE \rm durchgef"uhrt.}
\begin{eqnarray}
 a_1&=&\frac{1}{2}A,\label{ai}\\
 a_2&=&\frac{1}{2^4\cdot 3}(6A^2-B),\nonumber\\
 a_3&=&\frac{1}{2^5\cdot 3}(2A^3-AB+16C),\nonumber\\
 a_4&=&\frac{1}{2^9\cdot 3^2 \cdot 5}(60A^4-60A^2B+1920AC+7B^2-1152D),\nonumber\\
 a_5&=&\frac{1}{2^{10}\cdot 3^2 \cdot 5}(12A^5-20A^3B+960A^2C+7AB^2-1152AD+32CB).\nonumber
\end{eqnarray}

Die Polynome der zugeh"origen multiplikativen Folge sind
\vskip -0.4cm
\begin{equation}\label{folge}\end{equation}
\vskip -1.2cm
$$ \displaylines{\quad K_1=\frac{1}{2}Ac_1,\hfill}$$
$$ \displaylines{\quad K_2=\frac{1}{2^4\cdot 3}(2B c_2+(6A^2-B)c_1^2),\hfill}$$
$$ \displaylines{\quad K_3=\frac{1}{2^5\cdot 3}(48Cc_3+(2AB-48C)c_2c_1+(2A^3-AB+16C)c_1^3),\hfill}$$
$$ \displaylines{\quad K_4=\frac{1}{2^9\cdot 3^2 \cdot 5}((-8B^2+4608D)c_4+
(5760AC+8B^2-4608D)c_3c_1+\hfill\cr
\qquad\qquad\quad +(24B^2-2304D)c_2^2+ (120A^2B-5760AC-28B^2+4608D)c_1^2c_2+\hfill\cr
\qquad\qquad\quad +(60A^4-60A^2B+1920AC+7B^2-1152D)c_1^4),\hfill}$$
$$ \displaylines{\quad K_5=\frac{1}{2^{10}\cdot 3^2 \cdot 5}
(960BCc_5+(-8AB^2+4608AD-960BC)c_4c_1+\hfill\cr
\qquad\qquad\quad +(8AB^2+2880A^2C-4608AD+480BC)c_3c_1^2+(24AB^2-2304AD)c_2^2c_1+\hfill\cr
\qquad\qquad\quad +(40A^3B-2880A^2C-28AB^2+4608AD-160BC)c_2c_1^3+\hfill \cr
\qquad\qquad\quad +(12A^5-20A^3B+960A^2C+7AB^2-1152AD+32BC)c_1^5).\hfill}$$

Die spezielle Wahl f"ur die Polynome $q_1$ bis $q_4$ in $A$ bis $D$ erkl"art sich aus
den Werten von $\GU$ auf den in $\OU$ ausgezeichneten \MGFS $W_1$ bis $W_4$.

\SN{2.2.3} Auf den ersten f"unf \MGFS der in Abschnitt 1.4 konstruierten Basisfolge
$W_1$, $W_2$, $W_3$, $W_4$, $W_5$, $\dots$
f"ur $\OUQ$ nimmt das universelle elliptische Geschlecht die folgenden Werte an:
$\GU(W_1)=A$,\enspace $\GU(W_2)=B$,\enspace $\GU(W_3)=C$,\enspace $\GU(W_4)=D$ und $\GU(W_5)=0$.\SE

\bf Beweis: \rm Die Chernzahlen von $W_1$ bis $W_5$ wurden in Abschnitt 1.4 berechnet.
Einsetzen in die Gleichungen (\ref{folge}) ergibt die Behauptung.\BOX

Wir k"onnten hier durch Rechnen mit formalen Potenzreihen zu zeigen versuchen, da"s
$\GU(W_n)=0$ f"ur alle $n\ge5$ ist, d.h.~$\ker\GU=\langle W_5,W_6,\dots \rangle$.
Dies wird uns aber sp"ater im Zusammenhang mit dem Beweis der Inklusion
$\JSU\subset\ker\GU$ gelingen. F"ur \hbox{$SU$-\MGFS{}} k"onnen wir jetzt schon die folgende
Aussage beweisen:

\SN{2.2.4} Das universelle elliptische \GE einer $SU$-\MGF ist ein homogenes
Polynom, das nicht von der Unbestimmten $A$ abh"angt.\SE

{\bf Beweis:} Das Polynom $S(y)$ aus der Differentialgleichung (\ref{dgl})
hat in den Unbestimmten $A$ bis $D$ geschrieben die Gestalt
\begin{equation}\label{sy}
S(y)=\left(y+\frac{A}{2}\right)^4-\frac{1}{4}B\,\left(y+\frac{A}{2}\right)^2+4C\,\left(y+\frac{A}{2}\right)+\frac{1}{64}B^2-2D.
\end{equation}
Ist $\tilde h(x)$ die L"osung der Differentialgleichung (\ref{dgl}) f"ur $A=0$,
so ist daher $h(x)=\tilde h(x)-\frac{A}{2}$ die allgemeine L"osung.
F"ur die zugeh"origen charakteristischen Potenzreihen erh"alt man
die Beziehung $Q(x)=e^{(A/2)x}\tilde Q(x)$, denn aus $f(x)=e^{-(A/2)x}\tilde f(x)$
folgt $h(x)=\frac{f'(x)}{f(x)}=\tilde h(x)-\frac{A}{2}$, und die Potenzreihe
$Q(x)$ bestimmt sich nach Lemma 2.2.1 eindeutig aus $h(x)$. Wie im letzten
Abschnitt bemerkt (Gleichung (\ref{qsu})), ist der Wert der Geschlechter zu den
Potenzreihen $\tilde Q(x)$ und $e^{(A/2)x}\tilde Q(x)$ auf $SU$-\MGFS der
gleiche. Da $\tilde Q(x)$ zu $S(y)$ mit $A=0$ geh"ort, d.h.~nicht von $A$
abh"angt, folgt die Behauptung. \BOX


Setzen wir f"ur $A$, $B$, $C$ und $D$ komplexe Zahlen, so
k"onnen wir bei von Null verschiedener Diskriminante von $S(y)$ die
L"osung $h(x)$ der Differentialgleichung (\ref{dgl}) explizit angeben
(\cite{jung}, Abschnitt 1.1):
Sei $\wp(x)$ die Weierstra"ssche $\wp$-Funktion\index{$\wp$} (als Laurentreihe
im Nullpunkt geschrieben) zum Gitter $L$ mit den Gitterkonstanten
\begin{equation}\index{$g_2$,$g_3$}
g_2=\frac{1}{48}B^2-2D\quad\hbox{und}\quad g_3=-\frac{1}{1728}B^3+\frac{1}{12}BD-C^2.\nonumber
\end{equation}
Dieses Gitter ist nicht entartet, da wir $\hbox{Diskr }S(y)\not=0$ voraussetzten.
Die Differentialgleichung (\ref{dgl}) hat dann die L"osung
\begin{equation}\label{lsg}\index{$h(x)$}
h(x)=-\frac{1}{2}\,\frac{\wp'(x)+\wp'(z)}{\wp(x)-\wp(z)}-\frac{1}{2}A.
\end{equation}
Dabei ist $z$ ein von $0$ verschiedener Punkt auf der elliptischen Kurve
$\C/L$, der im Weierstra"smodell der Kurve die Koordinaten
$(\wp(z),\wp'(z))=(\frac{1}{24}B,C)$ besitzt.  Weiter wird in \cite{jung} (Lemma 1.1.7)
gezeigt, da"s die Menge der Paare $(L,z)$\index{$(L,z)$} mit $z\in\C/L$, $z\not=0$ und die
Polynome $S(y)$ mit $A=0$ und $\hbox{Diskr }S\not=0$ sich eineindeutig entsprechen.


Wenn $\C[A,B,C,D]$ die graduierte Koordinatenalgebra des gewichtet projektiven
Raumes $\CP^{1,2,3,4}$ ist, k"onnen wir nach der Konstruktion aus dem
letzten Abschnitt das universelle elliptische Geschlecht auch als das Geschlecht
$\GU:\OUC\to \C[A,B,C,D]$ ansehen, welches zur gewichtet projektiven Variet"at
$\CP^{1,2,3,4}$ und einer Basisfolge $W_1$, $W_2$, $W_3$, $W_4$, $W_5$, $V_6$,
$\dots$ geh"ort. Dabei seien $V_n$, $n\ge 6$ beliebige \MGFS aus $\ker
\GU\mid_{\Omega^U_n\otimes\C}$ mit $s_n(V_n)\not= 0$. Solche existieren stets,
denn sei $V_n'$ eine beliebige \MGF aus $\Omega_n^U\otimes\C$ mit $s_n(V_n')\not=0$,
so ist $\GU(V_n')$ ein Polynom $P(A,B,C,D)$ vom Gewicht $n$ in $A$ bis $D$.
F"ur $V_n:=V_n'-P(W_1,W_2,W_3,W_4)$ gilt dann $\GU(V_n)=0$ und $s_n(V_n)=s_n(V_n')\not=0$,
da $P(W_1,W_2,W_3,W_4)$ f"ur $n\ge5$ eine Summe aus lauter zerlegbaren \MGFS
ist.

\section{Die Geschlechter $\varphi_N$ und die Modulkurven $C_N$}

Wir f"uhren in diesem Abschnitt die komplexen elliptischen Geschlechter $\GN$ der Stufe $N$ als
Spezialf"alle des universellen komplexen elliptischen Geschlechtes $\GU$
ein. Sie k"onnen auch als die Geschlechter aufgefa"st werden, welche zu
den Modulkurven $\X1{N}$ geh"oren, die geeignet in den $\CP^{1,2,3,4}$ abgebildet sind.
Die zu den Spitzen von $\G1{N}$ geh"origen Geschlechter sind die altbekannten
Geschlechter $\chi_y(X)$ und $\chi(X,K^{k/N})$ (vgl.~\cite{hihab}).
Beweise zu den folgenden Tatsachen "uber die elliptischen Geschlechter der Stufe
$N$ finden sich in \cite{himod,hiell,jung}.

Sei $L$ ein Gitter in $\C$ und $\alpha \in \C/L$ ein von Null verschiedener
$N$-Teilungspunkt der
zugeh"origen elliptischen Kurve, d.h.~$N\cdot\alpha=0$. Es gibt genau eine
bzgl.~$L$ elliptische Funktion $k(x)$ mit Divisor $N\cdot(0)-N\cdot(\alpha)$ und der
Normierung $k(x)=x^N+O(x^{N+1})$ der Taylorentwicklung im Ursprung. Die
Funktion $f(x)=\root N \of {k(x)}$ ist wohldefiniert, falls gefordert wird,
da"s $f(x)=x+O(x^2)$. Sie ist elliptisch bzgl.~einem Untergitter $\tilde
L\subset L$, dessen Index gleich der Ordnung von $\alpha\in \C/L$ ist. Wir f"uhren
auf
$$\cal L \rm _N :=\{(L,\alpha)\mid L \in \C \mbox{ Gitter, } \alpha
\in \C/L \mbox{ $N$-Teilungspunkt }\not=0\} $$
die folgende "Aquivalenzrelation ein: $(L,\alpha)\sim(L',\alpha') :
\Leftrightarrow $ Es gibt ein $\mu\in\C\setminus\{0\}$ mit $(\mu L,\mu\alpha)
=(L',\alpha')$. Weiter sei $\cal L \rm _N^p \subset  \cal L \rm _N$ die
Menge der Paare $(L,\alpha)$, f"ur die $\alpha$ primitiver $N$-Teilungspunkt
ist, d.h.~$\alpha\in\C/L$ hat Ordnung $N$. Die Funktion
$Q(x):=\frac{x}{f(x)}$ transformiert sich beim "Ubergang zu "aquivalenten
Paaren wie folgt: Ist $Q_{L,\alpha}(x)$ die zum Paar $(L,\alpha)$ geh"orende
charakteristische Potenzreihe, so ist $Q_{\mu L,\mu\alpha}(x)=
Q_{L,\alpha}(\mu^{-1}x)$ die zum Paar $(\mu L,\mu\alpha)$ geh"orende
Potenzreihe (Normierung und Divisor f"ur $f(x)$ stimmen "uber\-ein), d.h.~die
Koeffizienten $a_k$ von $Q(x)$ sind homogene Gitterfunktionen vom Gewicht $-k$
bzgl.~$\cal L \rm _N$.

Auf der oberen Halbebene \H der komplexen Zahlen operiert die
Kongruenzuntergruppe\index{$\Gamma_1(N)$}
$$ \Gamma_1(N):=\lbrace\left(a\enspace b \atop c\enspace d\right)\in SL_2(\Z) \mid
c\equiv 0 \pmod{N},\enspace a\equiv d\equiv 1 \pmod{N}\rbrace$$
von $SL_2(\Z)$ verm"oge $\sigma:SL_2(\Z)\times\H\to\H$,
$(\left(a\enspace b \atop c\enspace d \right),\tau)
\mapsto \frac{a\tau+b}{c\tau+d}$. Der Bahnenraum sei mit $\H/\Gamma_1(N)$ bezeichnet.
F"ur $\cal L \rm _N^p$ gilt die folgende  Charakterisierung:

Die Zuordnung $\tau\mapsto (2\pi i(\Z\tau+\Z),\frac{2\pi i}{N})$ definiert eine
Bijektion zwischen $\H/\Gamma_1(N)$ und $\cal L \rm ^p_N/\simeq$.

Wegen $\cal L \rm _N=\DCUP_{\scriptstyle n\mid N \atop \scriptstyle n>1}\cal L \rm _n^p$
folgt, da"s $\DCUP_{\scriptstyle n\mid N \atop \scriptstyle n>1}\H/\Gamma_1(n)\to
\cal L \rm _N/\!\simeq$, $(\tau,n)\mapsto (2\pi i(\Z\tau+\Z),\frac{2\pi i}{n})$
bijektiv ist.

Die Koeffizienten $a_k$ der Potenzreihe $Q(x)$ zum Paar $(2\pi i(\Z\tau+\Z),
\frac{2\pi i}{N})$, aufgefa"st als Funktionen in $\tau\in\H$, transformieren sich
deshalb wie Modulformen vom Gewicht $k$ zur Modulgruppe $\Gamma_1(N)$.
Durch explizite Konstruktion von $f(x)$ zeigt man: $a_k(\tau)$ ist tats"achlich
Modulform zu $\Gamma_1(N)$, d.h.~$a_k(\tau)$ ist insbesondere holomorph in den
Spitzen, l"a"st sich also holomorph auf die kompaktifizierte Modulkurve $\X1{N}$
fortsetzen.

Um zu einer expliziten Darstellung von $f$ zu gelangen, betrachte f"ur
$q=e^{2\pi i\tau}$ die ganze Funktion
\begin{eqnarray}\label{Phi}\index{$\Phi(\tau,x)$}
\Phi(\tau,x)=(1-e^{-x})\prod_{n=1}^{\infty}(1-q^ne^{-x})(1-q^ne^x)/(1-q^n)^2,
\end{eqnarray}
welche Nullstellen der Ordnung $1$ in den Gitterpunkten von $L=2\pi i(\Z\tau+\Z)$
besitzt. Sie ist bis auf einen Faktor der Gestalt $e^{a x^2+b x}$ die Weierstra"ssche
Sigma-Funktion zu $L$. W"ahlt man f"ur $\alpha$ den $N$-Teilungspunkt $\frac{2\pi i}{N}$,
so erh"alt man die Darstellung
\begin{eqnarray}\label{fN}
f(x)=\frac{\Phi(\tau,x)\Phi(\tau,-\frac{2 \pi i}{N})}{\Phi(\tau,x-\frac{2\pi i}{N})}.
\end{eqnarray}
F"ur Kapitel 3 sei hier noch das genaue Transformationsverhalten von $f(x)$
zum $N$-Teilungspunkt $\alpha=\frac{2\pi i}{N}$ angegeben:
\begin{equation}\label{trans}
  f(x+2\pi i)=f(x),\qquad f(x+2\pi i\tau)=e^{-\frac{2\pi i}{N}} f(x).
\end{equation}

\DN{2.3.1} Das komplexe elliptische Geschlecht \index{$\varphi_N$}$\GN$ der Stufe $N$ ist das \GE
zur Potenzreihe $Q(x)=\frac{x}{f(x)}$ zum Gitter $L=2\pi i(\Z\tau+\Z)$ und
primitivem $N$-Teilungspunkt $\frac{2 \pi i}{N}$.

Da die Koeffizienten $a_k$ von $Q(x)$ Modulformen vom Gewicht $k$ sind, ist der
Wert von $\GN$ auf einer komplex $n$-dimensionalen $U$-\MGF $X_n$ eine
Modulform vom Gewicht $n$ zur Modulgruppe $\Gamma_1(N)$.

Die Funktion $f(x)$ gen"ugt nun zwei Differentialgleichungen.
Eine Differentialgleichung $2$-ter Ordnung ist von der gleichen Gestalt
wie die Differentialgleichung (\ref{dgl}) in Lemma 2.2.1:
\begin{equation}\label{dg1}
\left(\frac{f'}{f}\right)'^{\quad 2}=S\left(\frac{f'}{f}\right),
\end{equation}
wobei $S(y)=y^4+q_1y^3+q_2y^2+q_3y+q_4$ ein normiertes Polynom $4$-ten Grades
ist, dessen Koeffizienten $q_i$ Modulformen vom Gewicht $i$ zu $\G1{N}$ sind.
Da die Potenzreihe $f(x)$ zu $\GN$ der Differentialgleichung (\ref{dgl}) mit
speziellen $q_i$ gen"ugt, l"a"st sich $\GN$ "uber $\GU$ faktorisieren.
Die zweite Differentialgleichung ist von $1$-ter Ordnung und besitzt die Gestalt
\begin{equation}\label{dg2}\index{$P_N(y)$}
\frac{1}{f^N}+d_{2N}f^N=P_N\left(\frac{f'}{f}\right).
\end{equation}
Hier ist $P_N(y)=y^N+d_1y^{N-1}+\cdots+d_{N-1}y+d_N$ ein normiertes Polynom
$N$-ten Grades. Die~$d_i$ sind Modulformen vom Gewicht $i$ zu $\Gamma_1(N)$. Das
Polynom $P_N(y)$ ist ein sogenanntes Zolotarev-Polynom. Es hat die beiden
zus"atzlichen Eigenschaften $d_{N-1}=0$ und $P_N(y)^2=4d_{2N}$ f"ur die $y\not=0$
mit $P_N'(y)=0$.
Zwischen $S(y)$ und $P_N(y)$ besteht die Beziehung
\begin{equation}\label{bez}
S(y)P_N'(y)^2=N^2y^2(P_N(y)^2-4d_{2N}).
\end{equation}

Diese Differentialgleichung liefert zwei (von $N$ abh"angende) gewichtet homogene
Bedingungen $R_{N-1}$ und $R_{N+1}$\index{$R_{N-1}$, $R_{N+1}$} an die Koeffizienten
$q_1$ bis $q_4$ von $S(y)$. Die Polynome $R_{N-1}$ und $R_{N+1}$ definieren
eine Kurve $C_N$\index{$C_N$} im
gewichtet projektiven Raum $\CP^{1,2,3,4}$
mit den homogenen Koordinaten $q_1$ bis $q_4$, wobei $q_i$ mit dem Gewicht $i$
versehen sei.
Nun wird umgekehrt jeder Punkt der Modulkurve $C_N$ durch ein (evtl.~entartetes)
Gitter und einen $N$-Teilungspunkt repr"asentiert. Genauer gilt der folgende

\SN{2.3.2 (s.~\cite{jung}, Kapitel 3)}\index{$\Phi$} Die Abbildung
$$ \Phi:\DCUP_{\scriptstyle n\mid N \atop \scriptstyle n>1}\X1{n}\to\CP^{1,2,3,4},$$
die f"ur jedes $n\mid N$, $n>1$
dem Paar $(\tau,n)$ die Koeffizienten $q_1$
bis $q_4$ des Polynomes $S(y)$ aus der Differentialgleichung (\ref{dg1}) zuordnet
(wobei $f(x)^n$ die elliptische Funktion zu $L=2\pi i(\Z\tau+\Z)$ mit Divisor $n\cdot(0)-n\cdot(\frac{2\pi i}{n})$ und
Normierung $f(x)=x+O(x^2)$ ist), ist nach der Fortsetzung auf die
Spitzen surjektiv auf die Modulkurve $C_N$.
Sie ist bis auf Werte in den Spitzen auch injektiv. Die $\Phi(\X1{n})$
bilden gerade die Zerlegung von $C_N$ in die irreduziblen Komponenten.
Das Bild $\Phi(\X1{n})$ ist unabh"angig vom Vielfachen $N$ von $n$.
\SE

{\bf Beispiel $N=2$:}~F"ur $N=2$ erhalten wir das elliptische Geschlecht der Stufe $2$.\index{$\varphi_2$}
Das Polynom $S(y)$ hat die Gestalt
\begin{equation}
S(y)=y^4+4\delta y^2+4(\delta^2-\epsilon),
\end{equation}
die Differentialgleichung (\ref{dg2}) lautet
\begin{equation}
\frac{1}{f^2}+\epsilon f^2=(f'/f)^2+2\delta,\hbox{\quad oder\quad}
f'^2=\epsilon f^4-2\delta f^2+1.
\end{equation}
Setzen wir $q=e^{2\pi i\tau}$, so haben die Modulformen \index{$\delta,\epsilon$}
$\delta$ und $\epsilon$ vom Gewicht $2$ und $4$ die Form
\begin{eqnarray}
\delta  &=&\frac{1}{4}+6\,\sum_{n=1}^{\infty}\Bigl(\sum_{{\scriptstyle d \mid n}\atop{\scriptstyle d\equiv 1\,(2)}}d\Bigr)q^n,\\
\epsilon&=&\frac{1}{16}\,\prod_{n=1}^{\infty}\left(\frac{1-q^n}{1+q^n}\right)^8.
\end{eqnarray}
Um $\varphi_2$ als Spezialfall von $\GU$ zu schreiben, m"ussen wir im
$A$,$B$,$C$,$D$-Koordinatensystem (vgl.~(\ref{ABCD}))
\begin{equation}
A=0,\qquad B=-16\delta,\qquad C=0\qquad\hbox{und}\qquad D=2\epsilon\nonumber
\end{equation}
setzen. Da die Koeffizienten $a_n$ von $Q(x)$ nach Lemma 2.2.1 Polynome vom
Gewicht $n$ in $q_1$ bis $q_4$ sind und im Stufe-$2$-Fall $q_1=q_3=0$ gilt,
ist die Potenzreihe $Q(x)$ gerade, d.h.~das Geschlecht $\varphi_2$ ist schon f"ur
$SO$-\MGFS definiert. Setzen wir $\epsilon=0$, $\delta=-\frac{1}{8}$
(bzw.~$A=0$, $B=2$, $C=D=0$), so erhalten wir das $\hat A$-Geschlecht;
\index{$\hat A$}f"ur $\epsilon=\delta^2=1$ (bzw.~$A=0$, $B=-16$, $C=0$, $D=2$)
erhalten wir die Signatur.\index{$\sign$}

Wir wollen
nun die Geschlechter $\varphi_n$ f"ur $n\mid N$, $n>1$ zusammenfassen:

\DN{2.3.3}\index{$\tilde\varphi_N$}
Sei $W_1$, $W_2$, $W_3$, $W_4$, $V_5$, $V_6$, $\dots$ eine Basisfolge
von $\OUC$ mit $\langle V_5,V_6,\dots\rangle =\ker \GU$ und $W_1$ bis $W_4$ die \MGFS aus
Abschnitt 1.4. Es bezeichne dann $\tilde\varphi_N$ das \GE zur Variet"at $C_N=
V(\langle R_{N-1},R_{N+1}\rangle)\subset\CP^{1,2,3,4}$ wie in Abschnitt 2.1
definiert. Dabei ist das $A$,$B$,$C$,$D$-Koordinatensystem f"ur den
$\CP^{1,2,3,4}$ zu verwenden.
Die beiden Gleichungen $R_{N-1}$ und $R_{N+1}$
sind die homogenen Bedingungen vom Gewicht $N-1$ und $N+1$ an die
$q_1$ bis $q_4$ bzw.~$A$ bis $D$, die sich aus (\ref{bez}) ergeben.

Die Definition h"angt nat"urlich nicht von den gew"ahlten
Basismannigfaltigkeiten $V_5$, $V_6$,~$\dots$ aus $\ker\GU$ ab.

Im $A$,$B$,$C$,$D$-Koordinatensystem f"ur den $\CP^{1,2,3,4}$ besteht also
die folgende Situation
$$\begin{array}{ccccc}\typeout{zwei Pfeile einfuegen !!!}
\OUC &           & \mapr{\GU} & & \C[A,B,C,D]                   \\
     & \tilde\GN &            & & \downarrow\pi_N               \\
     &           &            & & K(C_N)\cong\C[A,B,C,D]/I(C_N).
\end{array}$$\typeout{Pfeil !!!!}

Der Kern von $\pi_N$ ist das Verschwindungsideal $I(C_N)$ der Kurve $C_N$
und nach dem Hilbertschen Nullstellensatz daher das Radikal von $\langle
R_{N-1},R_{N+1} \rangle$.

Den genauen Zusammenhang zwischen $\GN$ und $\tilde\GN$ beschreibt

\LN{2.3.4} $\ker\tilde\GN=\bigcap_{\scriptstyle n\mid N\atop\scriptstyle 1<n}\ker\GN $\SE

{\bf Beweis:} Wegen Lemma 2.1.1 und Satz 2.3.2 gilt $\ker\tilde\GN=\ker
\varphi_{C_N}=\break \bigcap_{\scriptstyle n\mid N\atop\scriptstyle 1<n}\ker\varphi_{\Phi(\X1{n})}$. Nach
Definition ist $X \in \ker \varphi_{\Phi(\X1{n})} \Leftrightarrow
\GU(X)\in I(\Phi(\X1{n}))$. Die Punkte $P \in \Phi(\X1{n})$ lassen sich durch
$\tau\in\X1{n}$ parametrisieren: $P(\tau)=(A(\tau):B(\tau):C(\tau):D(\tau))$. Da
$\GU(X)(P(\tau))=\varphi_n(X)(\tau)$, ist $\GU(X)\in I(\Phi(\X1{n}))
\Leftrightarrow\varphi_n(X)=0$.\BOX

Ist $N$ eine Primzahl, so ist die Kurve $C_N$ irreduzibel, und es gilt $\ker\tilde
\GN=\ker\GN$.\typeout{genauer: Bildraeume verschieden}

\SN{2.3.5} F"ur das Verschwindungsideal der Kurve $C_N$ gilt $$I(C_N)=
\rad(\langle R_{N-1},R_{N+1}\rangle)=\langle R_{N-1},R_{N+1}\rangle.$$
\SE

{\bf Beweis:} 

Zu zeigen ist, da"s das von den beiden homogenen Polynomen $R_{N-1}$ und $R_{N+1}$
erzeugte Ideal $\AN=\langle R_{N-1},R_{N+1} \rangle$ gleich seinem Radikal
$\rad(\AN)$ ist. Wir werden dies auf die in \cite{jung}, S.~45--47 bewiesene
entsprechende Eigenschaft des Ideals $\bar\AN:=\AN\cap\C[B,C,D]=\langle
\res_A(R_{N-1},R_{N+1})\rangle$ zur"uckf"uhren, dessen Verschwindungsvariet"at
$V(\bar\AN)\subset\CP^{2,3,4}$ die Projektion $\bar C_N$ der Kurve $C_N\subset
\CP^{1,2,3,4}$ auf die gewichtet projektive Ebene $\{A=0\}=\CP^{2,3,4}$ ist.
(Mit $\res_A(R_{N-1},R_{N+1})$ sei die Resultante zwischen $R_{N-1}$ und $R_{N+1}$
bez"uglich der Variablen $A$ bezeichnet.)

Sei $\AN=\cQ_1\cap\dots\cap \cQ_r$ eine reduzierte Prim"arzerlegung von $\AN$.
Dabei m"oge das Prim"arideal $\cQ_i$ zum Primideal $\cP_i$ geh"oren, und
es es kann angenommen werden, da"s die Ideale $\cQ_i$ und $\cP_i$ alle
homogen sind (s.~\cite{kunz}, S.~197, Lemma 4.1). Jedes $\cP_i$-prim"are
Ideal $\cQ_i$ besitzt eine L"ange $l_i$. Dies ist die maximale L"ange $l$
einer echt aufsteigenden Kette von $\cP_i$-prim"aren Idealen
$$\cQ_i=\cR_1\subset\cR_2\subset\dots\subset\cR_l=\cP_i.$$
Das Prim"arideal $\cQ_i$ ist genau dann prim, wenn seine L"ange $l_i=1$ ist
(s.~\cite{gro1}, S.~134).

Das Ideal $\AN$ ist ein vollst"andiger Durchschnitt (s.~\cite{kunz} S.~140),
denn $\AN=\langle R_{N-1},R_{N+1} \rangle$ hat eine $2$-gliedrige Basis
und die H"ohe oder Rang ist ebenfalls $2$. Letzteres, da $V(\AN)=C_N$ nur
aus $1$-dimensionalen Komponenten besteht, also $\AN$ die homomogene
Dimension $1$ bzw.~die Krulldimension $\kdim \AN =2$ hat, und damit
$\rang \AN=\kdim \C[A,B,C,D]-\kdim \AN=4-2=2$ ist. F"ur vollst"andige
Durchschnitte gilt nun der Ungemischtheitssatz, der besagt, da"s das Ideal
$\AN$ ungemischt ist, d.h.~alle Prim"arideale $\cQ_i$ den Rang $2$ haben.
Sie entsprechen eineindeutig den irreduziblen  Komponenten von $C_N$;
es gibt also keine eingebetteten Komponenten (s.~\cite{kunz}, S.~193
(Polynomringe sind Cohen-Macaulay-Ringe) oder \cite{gro2}, S.~185 (vollst"andige
Durchschnitte hei"sen dort Hauptklassenideale)).

Es gilt $\rad(\AN)=\cP_1\cap\dots\cap\cP_r$, und die Primideale $\cP_i$
sind die Verschwindungsideale der einzelnen irreduziblen Komponenten der
Kurve $C_N$. Wir haben also $l_i=1$ f"ur $i=1,\dots,r$ zu zeigen.

Projiziert man die Kurve $C_N=V(\AN)$ verm"oge $\pi:\CP^{1,2,3,4}\to\CP^{2,3,4}$,
$(A:B:C:D)\mapsto (B:C:D)$ in die gewichtet projektive Ebene $\CP^{2,3,4}$,
so erh"alt man eine Kurve $\bar C_N$. F"ur ihr Verschwindungsideal $I(\bar C_N)$
wurde in \cite{jung}, S.~45
$$I(\bar C_N)=\langle \res_A(R_{N-1},R_{N+1})\rangle$$
bewiesen. (Man beachte, da"s das hier verwendete $A,B,C,D$-Koordinatensystem und das in
\cite{jung} verwendete so zusammenh"angen, da"s die Projektionsabbildungen die
gleichen sind.)

Um die Ideale $I(C_N)$ und $I(\bar C_N)$ miteinander zu vergleichen,
ben"otigen wir den Grad $h_0(I)$
eines homogenen Ideals $I$ in einem graduierten Polynomring $\C[x_1,\dots,
x_k]$, $\grad x_i=$ \hbox{$d_i\in\N$}. Die "ubliche Definition "uber den Leitkoeffizienten
des Hilbertpolynomes von $I$, die man bei der gew"ohnlichen Graduierung
$d_1=\dots=d_k=1$ verwendet, ist bei beliebiger Graduierung nicht verwendbar,
da die Hilbertfunktion $f(n)$ auch f"ur gro"se $n$ im allgemeinen kein Polynom
mehr sein wird.

Mit der Bezeichnung $M^{(n)}$ f"ur die Elemente vom Gewicht $n$ eines graduierten
$\C[x_1,\dots,x_k]$-Moduls $M$ (also $M=\bigoplus_{n\ge0}M^{(n)}$, $p\cdot
M^{(n)}\subset M^{(n+l)}$ f"ur $p\in \C[x_1,\dots,x_k]$ homogen vom Gewicht $l$)
definiere f"ur ein homogenes Ideal $I$
$$P_I(t)=\sum_{n=0}^{\infty}\dimc(\C[x_1,\dots,x_k]/I)^{(n)}\cdot t^n.$$
So gilt zum Beispiel
$$P_{\langle 0\rangle}(t)=\prod_{i=1}^{k}\frac{1}{1-t^{d_i}}.$$
Ganz allgemein ist die Poincar\'ereihe $P_I(t)$ nach \cite{atmac}, S.~117
eine rationale Funktion der Gestalt
$$P_I(t)=\frac{Q(t)}{\prod_{i=1}^{k}1-t^{d_i}}$$
mit einem Polynom
$Q(t)$. Die Ordnung der Polstelle bei $t=1$ ist gerade die Krulldimension
von $I$ (und um $1$ gr"o"ser als die Dimension von $V(I)$).
Schreiben wir $P_I(t)$ in der Form
$$ P_I(t)=\frac{1}{(1-t)^{\kdim I}}\cdot \frac{\tilde Q(t)}{\prod_{i=1}^{k}
(\sum_{l=0}^{d_i-1}t^l)},$$
mit einem Polynom $\tilde Q(t)=Q(t)/(1-t)^{k-\kdim I}$,
so definieren wir den Grad $h_0(I)$ des Ideals $I$ durch
$$h_0(I)=\tilde Q(t)\mid_{t=1}.$$
Diese Definition stimmt wegen $\sum_{n=0}^{\infty}h_0(I)\cdot{n+(\kdim I)-1
\choose (\kdim I)-1)}\cdot t^n=\frac{1}{(1-t)^{\kdim I}}\cdot h_0(I)$ f"ur
$d_1=\cdots=d_k=1$ mit derjenigen Definition, die das Hilbertpolynom
verwendet, "uberein. Es gilt stets
\begin{equation}\label{hnull}
h_0(I)>0,
\end{equation}
denn die linke Seite von
$$P_I(t)\cdot\prod_{i=1}^k\bigl(\sum_{l=0}^{d_i-1}t^l\bigr) =
                        \frac{1}{(1-t)^{\kdim I}}(h_0+h_1(1-t)+\cdots) $$
ist eine Potenzreihe in $t$ mit lauter nichtnegativen Koeffizienten.

Wir ben"otigen nun einige Eigenschaften der Poincar\'ereihe und des Grades,
die im Falle des standardgraduierten $\C[x_1,\dots,x_k]$ wohlbekannt sind.
(vgl.~\cite{vogel}, Kap.~1 Abschnit C, S.~43 ff.).

Sei $I\subset\C[x_1,\dots,x_k]$ ein homogenes Ideal und $\varphi\in\C[x_1,\dots,x_k]$ ein
homogenes Polynom vom Grad $r$. Dann gilt (vgl.~\cite{vogel}, S.~43 f.)
\begin{eqnarray}\label{poidi}\label{poidiii}
P_{I+\langle \varphi\rangle}(t)&=&P_I(t)-P_{(I:\varphi)}(t)\cdot t^r \\
&=&P_I(t)\cdot(1-t^r)\qquad\qquad \hbox{wenn $(I:\varphi)=I$},\nonumber
\end{eqnarray}
Zum Beweis verwende $\dimc (I+\langle\varphi\rangle)
^{(n)}=\dimc I^{(n)}+\dimc (\langle\varphi\rangle)^{(n)}-\dimc (I\cap
\langle\varphi\rangle)^{(n)}$, sowie $(I\cap\langle\varphi\rangle)=(I:\varphi)
\cdot\langle\varphi\rangle$ und $\dimc ((I:\varphi)\cdot\langle\varphi\rangle)
^{(n)}=\dimc (I:\varphi)^{(n-r)}$.

Seien $\varphi_1,\dots,\varphi_s\in\C[x_1,\dots,x_k]$ Formen mit den zugeh"origen
Graden $r_1,\dots,r_s$. Wenn $(\langle \varphi_1,\dots,\varphi_{i-1}\rangle
:\varphi_i)=\langle\varphi_1,\dots,\varphi_i\rangle$ f"ur alle $1\le i\le s$ gilt,
so hat man mit $I=\langle\varphi_1,\dots,\varphi_s\rangle$
\begin{equation}\label{vd}
P_I(t)=\frac{\prod_{j=1}^s 1-t^{r_j}}{\prod_{i=1}^k 1-t^{d_i}},\quad\hbox{ also }
\quad h_0(I)=r_1\cdot\dots\cdot r_s.
\end{equation}
Beweis von (\ref{vd}) mittels vollst"andiger Induktion:
Es ist $P_{\langle \varphi_1\rangle}(t)=\frac{1-t^{r_1}}{\prod_{i=1}^k 1-t^{d_i}}$
nach (\ref{poidi}). Sei (\ref{vd}) f"ur $s=i-1$ schon bewiesen. Dann ergibt sich
mit (\ref{poidiii})
\begin{eqnarray*}
P_{\langle \varphi_1,\cdots,\varphi_i\rangle}(t)&=&P_{\langle \varphi_1,\cdots,\varphi_{i-1}\rangle}(t)\cdot
(1-t^{r_i})\\
&=&\frac{\prod_{j=1}^{i-1} 1-t^{r_j}}{\prod_{i=1}^k 1-t^{d_i}}\cdot(1-t^{r_i})
=\frac{\prod_{j=1}^{i} 1-t^{r_j}}{\prod_{i=1}^k 1-t^{d_i}}.
\end{eqnarray*}

Aus (\ref{poidiii}) erh"alt man sofort
(vgl.~\cite{vogel}, S.~47 (1.36)):
\begin{equation}\label{hi}
\qquad\ \ \hbox{wenn $\kdim(I,\varphi)=\kdim(I)=\kdim(I:\varphi)$, so\ \,}
h_0(I,\varphi)=h_0(I)-h_0((I:\varphi)),
\end{equation}
\begin{equation}\label{hii}
\quad\hbox{wenn $\kdim(I,\varphi)=\kdim(I)>\kdim(I:\varphi)$, so\ \,}
h_0(I,\varphi)=h_0(I).\qquad\quad\ \ \
\end{equation}

Sei $\cP\subset\C[x_1,\dots,x_k]$ ein homogenes Primideal und $\cQ$ ein
homogenes $\cP$-prim"ares Ideal der L"ange $l(\cQ)$. Dann gilt
\begin{equation}\label{gradf}
h_0(\cQ)=l(\cQ)\cdot h_0(\cP),
\end{equation}
denn der Beweis von (\ref{gradf}) in \cite{vogel}, S.~47 f.~verwendet nur die
beiden Eigenschaften (\ref{hi}) und (\ref{hii}) des Grades $h_0$, so da"s er
auch in unserer gewichtet homogenen Situation g"ultig bleibt.

Sei $I\subset\C[x_1,\dots,x_k]$ ein homogenes Ideal mit der Prim"arzerlegung
$I=\cQ_1\cap\dots\cap\cQ_r$, wobei $\cQ_i$ zum Primideal $\cP_i$ geh"ort.
Ebenfalls wie in \cite{vogel}, S.~49 erh"alt man
\begin{equation}\label{grad}
h_0(I)=\sum_{\cQ} l(\cQ)\cdot h_0(\cP),
\end{equation}
wobei $\cQ$ alle $\cP$-prim"aren Komponenten von $I$ mit $\kdim \cQ=\kdim I$
durchl"auft. Insbesondere folgt aus (\ref{grad})
\begin{equation}\label{ungl1}
h_0(I)\ge h_0(\rad(I)).
\end{equation}
Da unser Ideal $\AN\subset\C[A,B,C,D]$ ungemischt ist, folgt umgekehrt
$\AN=\rad(\AN)$, wenn wir $h_0(\AN)=h_0(\rad(\AN))$ zeigen k"onnen.

Sei $I\subset\C[A,B,C,D]$ ein homogenes Ideal mit $R_{N-1}\in I$, sei
$\bar I=I\cap\C[B,C,D]$ das Eliminationsideal. Definiere f"ur $l \ge 1$ das
homogene Ideal $\bar I_l$ als das Ideal, das von den Anfangskoeffizienten
$\psi_l$ der homogenen Polynome
\begin{equation}\label{il}
\psi_l\cdot A^l+\psi_{l-1}\cdot A^{l-1}+\cdots+\psi_1\cdot A+\psi_0\,\in\,I,
\quad\psi_i\in\C[B,C,D]
\end{equation}
gebildet wird. Wir erhalten die aufsteigende Idealkette
$$I\,\subset\,\bar I_1\,\subset\,\bar I_2\,\subset\,\dots\,\subset\,\bar I_s\,
=\,\bar I_{s+1}\,=\,\dots\quad\hbox{in $\C[B,C,D]$.}$$
F"ur unser spezielles Ideal $I$ gilt
\begin{equation}\label{AA}
\bar I_{N-1}=\bar I_{N}=\dots =\C[B,C,D].
\end{equation}
Beweis von (\ref{AA}): Das nach Voraussetzung in $I$ enthaltende Polynom
$R_{N-1}$ hat nach \cite{jung}, S.~42 die Gestalt $R_{N-1}=\alpha\cdot
A^{N-1}+\cdots$ mit $\alpha\in\C$, $\alpha\not=0$ ($q_1^{N-1}=2^{N-1}\cdot
A^{N-1}$). Ist $\psi\in\C[B,C,D]$ eine Form, so liegt wegen $\frac{1}{\alpha}
\cdot\psi\cdot R_{N-1}=\psi\cdot A^{N-1}+\cdots\in I$ die Form $\psi$ in
$\bar I_{N-1}$, d.h.~$\bar I_{N-1}=\C[B,C,D]$.

Z"ahlt man die "uber $\C$ linear unabh"angigen Polynome vom Grad $n$ in $I$ ab,
erh"alt man f"ur $n\ge N-2$ unter Ausn"utzung der Darstellung (\ref{il}) und
mit (\ref{AA})
$$\dimc I^{(n)}=\dimc \bar I^{(n)}+\sum_{k=1}^{N-2}\dimc \bar I_k^{(n-k)}+
\sum_{k=N-1}^n \dimc \C[B,C,D]^{(n-k)}.$$
Unter Verwendung von $\dimc\C[A,B,C,D]^{(n)}=\sum_{k=0}^n\dim \C[B,C,D]^{(n-k)}$
folgt
$$ \dimc (\C[A,B,C,D]/I)^{(n)}=\dimc(\C[B,C,D]/\bar I)^{(n)}+\sum_{k=1}^{N-2}
\dimc (\C[B,C,D]/\bar I_k)^{(n-k)},$$
also
$$P_I(t)=P_{\bar I}(t)+\sum_{k=1}^{N-2}P_{\bar I_k}(t)\cdot t^k +S(t)$$
mit einem Polynom $S(t)$.\footnote{Beispiel $N=3$: Sei $I=\langle R_{N-1},
R_{N+1}\rangle=\langle B+\frac{3}{4}A^2,D+\frac{1}{2}AC\rangle$. Dann
$\bar I=\langle C^2 B+3D^2 \rangle$, $\bar I_1=\langle C,D\rangle$, $\bar I_2=
\bar I_3=\dots=\langle 1\rangle$. F"ur die Poincar\'ereihen gilt
$$\frac{(1-t^2)(1-t^4)}{(1-t)(1-t^2)(1-t^3)(1-t^4)}=\frac{1-t^8}{(1-t^2)
(1-t^3)(1-t^4)}+\frac{(1-t^3)(1-t^4)}{(1-t^2)(1-t^3)(1-t^4)}\cdot t.$$}
Zusammen mit (\ref{hnull}) ergibt sich f"ur ein $l \le N-2$
\begin{equation}\label{ungl2}
h_0(I)=h_0(\bar I)+\sum_{k=1}^{l}h_0(\bar I_k)\ge h_0(\bar I).
\end{equation}

Ist jetzt $I=I(C_N)=\rad(\AN)\ni R_{N-1}$, so gilt $\bar I=\rad(\AN)\cap
\C[B,C,D]=\langle \res_A(R_{N-1},R_{N+1}) \rangle$ f"ur das Eliminationsideal.
Denn es ist $\res_A(R_{N-1},R_{N+1})\in\bar I$, und umgekehrt folgt
aus $f\in I(C_N)\cap\C[B,C,D]$, da"s $f$ in $I(\bar C_N)=\langle\res_A(R_{N-1},
R_{N+1})\rangle$ liegt. Wegen $\kdim \AN=2$ gilt $(\langle R_{N-1}\rangle:\langle
R_{N+1}\rangle)=\langle R_{N-1} \rangle$ und mit (\ref{vd}) daher
$$P_{\AN}(t)=\frac{(1-t^{N-1})(1-t^{N+1})}{(1-t)(1-t^2)(1-t^3)(1-t^4)},\quad
\hbox{also $h_0(\AN)=N^2-1$.}$$
Da $\overline{\rad (\AN)}=\langle\res_A(R_{N-1},R_{N+1})\rangle$ homogen vom
Grad $N^2-1$ ist, hat man
$$P_{\overline{\rad (\AN)}}(t)=\frac{1-t^{N^2-1}}{(1-t^2)(1-t^3)(1-t^4)},\quad
\hbox{also $h_0(\overline{\rad (\AN)})=N^2-1$.}$$
Zusammen mit (\ref{ungl1}) und (\ref{ungl2}) erhalten wir schlie"slich
$$N^2-1=h_0(\AN)\ge h_0(\rad(\AN))\ge h_0(\overline{\rad (\AN)})=N^2-1,$$
d.h.~$h_0(\AN)= h_0(\rad(\AN))$, w.z.b.w.\BOX

Kommen wir nun zu einer genaueren Untersuchung der Geschlechter, die zu den
Spitzen der Modulkurven geh"oren.

F"ur die Anzahl $cu(N)$ der Spitzen der Modulkurve $\X1{N}$ gilt (vgl.~\cite{himod}, S.~174):
$$ cu(2)=2,\quad cu(4)=3,\quad cu(N)=\frac{1}{2}\sum_{n\mid N}\varphi(n)
\varphi(N/n)\hbox{\enspace sonst,}$$
dabei ist $\varphi$ die Eulersche $\varphi$-Funktion.

Einer Modulform $f$ zu $\Gamma_1(N)$ kann man die Werte in den Spitzen zuordnen,
die bei ungeradem Gewicht allerdings nur bis auf ein Vorzeichen wohldefiniert
sind.
\typeout{Anzahl Spitzen und Probleme mit dem Vorzeichen}
Die Abbildung $\Phi :\X1{N}\to\CP^{1,2,3,4},\tau\mapsto
(q_1(\tau):\dots:q_4(\tau))$ aus Satz 2.3.2 ist wegen $(q_1:q_2:q_3:q_4)\sim
(-q_1:q_2:-q_3:q_4)$ auch in den Spitzen wohldefiniert. "Uber die Werte von $\Phi$ in den
Spitzen von $\G1{N}$ gibt der folgende Satz Auskunft.

\SN{2.3.6 (s.~\cite{himod}, S.~124)} Das Polynom $S(z)=z^4+q_1z^3+q_2z^2+q_3z+q_4$
hat in den Spitzen von $\G1{N}$ die folgende Gestalt
\begin{eqnarray*}
\hbox{\quad (i)\enspace }S(z)&=&(z-k/N)^2(z+(N-k)/N)^2, \hbox{ mit $0<k<N$}\\
\hbox{\enspace oder (ii)\enspace}S(z)&=&z^2(z^2+2\frac{1-y}{1+y}z+1),
\hbox{mit $-y=e^{2\pi i l/N}$ und ggT$(l,N)=1$,}\\
\noalign{\vskip 2mm\hbox{bzw.~im $A,B,C,D$-Koordinatensystem geschrieben}\vskip 3mm}
\hbox{\quad (i')\enspace }S(z)&=&\bigl((z+\frac{A}{2})^2-\frac{B}{8}\bigr)^2,\\
\noalign{\hbox{mit $A=2(1/2-k/N)$, $0<k<N$, $B=2$ und $C=D=0$}\vskip 2mm}
\hbox{oder (ii')\enspace } S(z)&=&\left(z+\frac{A}{2}\right)^4-\frac{1}{4}B\,\left(z+\frac{A}{2}\right)^2
+4C\,\left(z+\frac{A}{2}\right)+\frac{1}{64}B^2-2D,\\
\noalign{\vskip1mm\hbox{mit $A=\frac{1-y}{1+y}$, $B=\frac{2(y^2-10y+1)}{(1+y)^2}$,
$C=\frac{y(y-1)}{(1+y)^3}$ und $D=\frac{y(-y^2+4y-1)}{(1+y)^4}$,}}
\noalign{\vskip1mm\hbox{$-y=e^{2\pi i l/N}$ und ggT$(l,N)=1$.\hfill}}
\end{eqnarray*}\SE
\vskip -7mm

Die Werte der Spitzen vom Typ (i) liegen also alle auf der gewichtet projektiven Geraden
\hbox{$\{C=D=0\}=\CP^{1,2}$}. "Uber die zu den Spitzen geh"origen komplexen
Geschlechter --- Spezialf"alle von $\GN$ bzw.~$\GU$ --- informiert

\SN{2.3.7} In den Spitzen nimmt das elliptische Geschlecht $\GN$ einer
komplex $n$-\nobreak dimensionalen $U$-\MGF~$X_n$ die folgenden Werte an:\index{$\chi_y$}\index{$\chi(\,.\,,K^{k/N})$}
\begin{eqnarray*}
\hbox{\quad(i) }&& \chi(X_n,K^{k/N})\hbox{\enspace mit $0<k<N$,}\\
\noalign{\hbox{dem \GE, das zur Potenzreihe $Q(x)=\frac{x}{1-e^{-x}}e^{-(k/N)x}$ geh"ort}\vskip 3mm}
\hbox{\enspace oder (ii) }&& \chi_y(X_n)/(1+y)^n\
\hbox{mit $-y=e^{\frac{2\pi i l}{N}}$ und ggT$(l,N)=1$,}\\
\noalign{\hbox{dem \GE, das zur Potenzreihe $Q(x)=\frac{x}{1-e^{-x}}\frac{1+ye^{-x}}{1+y}$
geh"ort.}}
\end{eqnarray*}\SE
\vskip -8mm

\DN{2.3.8}\index{$\tilde A$} Sei $\tilde A$ das Geschlecht, das zur Potenzreihe
$$ Q(X)=e^{\frac{A}{2}x}\frac{\sqrt{B/2}\cdot(x/2)}{\sinh\left(
\sqrt{B/2}\cdot(x/2)\right)}$$
geh"ort.\footnote{
$\tilde A(X)$ ist im wesentlichen das Hilbertpolynom von $X$.}
\typeout{Zusammenhang A-tilde mit Hibertpolynom vgl. theproof}

Die Potenzreihe zu $\tilde A$ gen"ugt der Differentialgleichung (\ref{dgl}) mit $C=D=0$
und $\tilde A$ l"a"st sich daher "uber $\GU$ faktorisieren. Man kann $\tilde A$ also
auch als das Geschlecht zur Kurve $\CP^{1,2}$ im \hbox{$A,B$}-Koordinatensystem
und zu der wie in Definition 2.3.3 gew"ahlten
Basisfolge $W_1$, $W_2$, $W_3$, $W_4$, $V_5$, $\dots$ auffassen. Die zugeh"orige
multiplikative Folge von Polynomen erh"alt man aus der von $\GU$, indem man
$C=D=0$ setzt.

Die Geschlechter, die zu den Spitzen von Typ (i) geh"oren, wollen wir nun
zusammenfassen. Die Bilder unter $\Phi$ der Spitzen vom Typ (i) auf der projektiven
Geraden $\CP^{1,2}$ im $A$, $B$-Koordinatensystem sind nach Satz 2.3.6. (i') die
Punkte $P_{N,k}:=(2(1/2-k/N)t:2t^2)$, $0<k<N$, wobei die Punkte $P_{N,k}$ und
$P_{N,N-k}$ zusammenfallen (setze $t=-t$). Sei $M_N=\bigcup_{k=1}^{[\frac{N}{2}
]}P_{N,k}$ die Vereinigung der Punkte $P_{N,k}$. Die homogenen Polynome $T_{N-1}$,
die das Verschwindungsideal der Variet"at $M_N$ im $\CP^{1,2}$ erzeugen,
haben offenbar die folgende Gestalt:\index{$T_{N-1}$}
$$T_{N-1}=\cases{\prod_{k=1}^{[\frac{N}{2}]}\left(A^2-2(1/2-k/N)^2B\right), &
falls $N$ ungerade, \cr A\cdot\prod_{k=1}^{[\frac{N}{2}]-1}\left(A^2-2(1/2-k/N)^2B\right), &
falls $N$ gerade.\cr}$$

\DN{2.3.9}\index{$\tilde A_N$} Sei $\tilde A_N$ das Geschlecht, das zur Variet"at $M_N\subset\CP
^{1,2}$ im $A$, $B$-Koordinatensystem
und der Basisfolge $W_1$, $W_2$, $W_3$, $W_4$, $V_5$, $\dots$ geh"ort.

\LN{2.3.10} Es gilt
$$ \ker \tilde A_N=\bigcap_{k=1}^{N-1}\ker\chi(X,K^{k/N}).$$\SE

{\bf Beweis:} Wegen Lemma 2.1.1 gilt $\ker\tilde A_N=\bigcap_{k=1}^{N-1}
\ker\varphi_{P_{N,k}}$. Nach Definition bedeutet $X\in \ker\varphi_{P_{N,k}}$
da"s $\tilde A(X)$ im Ideal $\langle A^2-2(1/2-k/N)^2B\rangle$ ($k\not=N/2$)
bzw.~in $\langle A\rangle$ ($k=N/2$) liegt.
Dies ist "aquivalent zu $X \in\ker\chi(X,K^{k/N})=\ker
\chi(X,K^{(N-k)/N})$, denn f"ur $A=2(1/2-k/N)t$ und $B=2t^2$ ist $\tilde A$
das Geschlecht zu $Q(X)=e^{(1/2-k/N)tx}\frac{tx/2}{\sinh(tx/2)}=
e^{-(k/N)tx}\frac{tx}{1-e^{-tx}}$, der homogen geschriebenen charakteristischen
Potenzreihe von $\chi(X,K^{k/N})=(-1)^{\dim_{\C}X}\chi(X,K^{(N-k)/N})$.\BOX

Das Toddsche- oder arithmetische Geschlecht \index{$Td$}$\hbox{Td}(X)$ erh"alt man aus
dem universellen Geschlecht $\GU(X)$ f"ur $C=D=0$ und $A=1$, $B=2$.



%


\section{Berechnung der Ideale $I_*^{N,1}$, $J_*^N$, $I_*^N$
und $I_*^{SU,t}$, $J_*^{SU}$, $I_*^{SU}$}

Mit den Ergebnissen aus dem letzten Abschnitt werden wir zeigen, da"s
die in Kapitel 1 konstruierten \MGFS ausreichen, die Ideale zu erzeugen.
Umgekehrt erhalten wir eine Charakterisierung der komplexen elliptischen
Geschlechter.

Wichtiges Hilfsmittel ist

\SN{2.4.1 (Hirzebruch \cite{hiell}, S.~58)} Sei $X$ eine $N$-\MGF mit
\SO. Ist der Typ $t$ der Operation $\not\equiv 0\!\! \pmod{N}$, dann
gilt $\GN(X)=0$.\hfil\break
Das elliptische Geschlecht $\GN$ der Stufe $N$ ist strikt multiplikativ
in $U$-Faserb"undeln (d.h.~Totalraum $E$, Basis $B$ und Faser $F$ sind
$U$-\MGFS mit vertr"aglicher $U$-Struktur (s.~\cite{bohi}, II (21.8))) mit
$N$-\MGFS $F$ als Faser und kompakter zusammenh"angender Liegruppe $G$ von
$U$-Automorphismen von $F$ als Strukturgruppe. Also gilt insbesondere
$\GN(E)=\GN(F)\cdot\GN(B)$.\SE

In \cite{hiell} wird der Satz nur f"ur komplexe \MGF formuliert, in der
Einleitung wird aber darauf hingewiesen, da"s er auch f"ur $U$-\MGFS
und $S^1$-Operationen, die die stabil fastkomplexe Struktur respektieren,
richtig bleibt.

{\bf Korollar 2.4.2:}\sl \ Es gelten die beiden Inklusionen
$$\hbox{1.\quad} I_*^{N,1}\subset\ker \tilde\GN \hbox{\quad und\quad 2.\quad}
J_*^N\subset \ker \tilde\GN.$$\SE

{\bf Beweis:} Zu 1. \ Eine $N$-\MGF $X$ mit \SO vom Typ $t$ ist f"ur
$n\mid N$ auch eine $n$-\MGF mit \SO vom Typ $t\equiv 1\pmod{n}$. F"ur
$t\equiv 1\pmod{N}$ folgt daher aus Satz 2.4.1 $X\in \bigcap_{{n\mid N}\atop
{n>1}}\ker \varphi_n$, wegen Lemma 2.3.4 $X\in \ker\tilde\GN$.

Zu 2. \ Sei $X\in J_*^N$, also $X=\widetilde\CP(E\oplus F)$ getwistet projektives
B"undel mit Basis $B$ und Faser $\widetilde\CP_{p,q}$, $N\mid p-q$. Der getwistet
projektive Raum $\widetilde\CP_{p,q}$ l"a"st als ge\-twistet projektives B"undel zu den
trivialen B"undeln $E'=\C^p$ und $F'=\C^q$ "uber einem Punkt nach Satz 1.3.6
\SOS von beliebigem Typ zu. Satz 2.4.1 liefert $\varphi_n(X)=
\varphi_n(X)\cdot\varphi_n(\widetilde\CP_{p,q})=\varphi_n(X)\cdot 0=0$ f"ur alle $n\mid N$,
$n>1$, da die Faser auch eine $n$-\MGF ist. Wie in 1. folgt $X\in\ker\tilde
\GN$. (Mit Hilfe des Residuensatzes kann man f"ur $n\mid p-q$ unter Verwendung von (\ref{trans})
$\varphi_n(\widetilde\CP_{p,q})=0$ auch direkt nachrechnen.)
\BOX

Damit k"onnen wir nun zeigen:

\SN{2.4.3}\index{$I_*^{SU,t}$}\index{$J_*^{SU}$} Seien $W_5, W_6, \dots$ die getwistet projektiven B"undel aus
Abschnitt 1.4. Dann gilt f"ur $t\not=0:$
$$ I_*^{SU,t}=J_*^{SU}=\langle W_5, W_6, \dots \rangle = \ker \GU \mid_{\OSUC}.$$
\SE

{\bf Beweis:} Wir zeigen
\begin{eqnarray*}
1.&&\langle W_5, W_6, \dots \rangle\subset I_*^{SU,t},\quad
\langle W_5, W_6, \dots \rangle\subset J_*^{SU},\\
2.&& I_*^{SU,t}\subset\ker\GU\mid_{\OSUC},\quad J_*^{SU}\subset\ker\GU\mid_{\OSUC},\\
3.&& \ker\GU\mid_{\OSUC}\subset\langle W_5, W_6, \dots \rangle.
\end{eqnarray*}
Zu 1. Die \MGFS $W_5, W_6, \dots$ wurden in Abschnit 1.4 gerade so konstruiert
(Lemma 1.4.7).

Zu 2. \ Sei $X_n\in I_*^{SU,t}$ eine $n$-dimensionale $SU$-\MGF, sei $N$ eine
ganze Zahl teilerfremd zu $t$ und gr"o"ser als $n+1$. Da eine $SU$-\MGF auch
eine $N$-\MGF ist, gilt $X_n\in I_*^{N,t\,\bmod\,N}=I_*^{N,ggT(t,N)}=I_*^{N,1}$.
Nach Korollar 2.4.2 ist $\tilde\GN(X_n) =\pi_N\circ\GU(X_n)=0$. Mit Satz 2.3.5
folgt $\GU(X_n)\in\ker\pi_N=\langle R_{N-1}, R_{N+1}\rangle$. Wegen $N-1>n$
mu"s $\GU(X_n)=0$ sein, da $R_{N-1}$ und $R_{N+1}$ homogen vom Grad $N-1$
bzw.~$N+1$ sind.
Analog hat man $X_n\in J_*^{SU} \Rightarrow X_n\in J_*^N \Rightarrow
\tilde\GN(X_n)=0 \Rightarrow \GU(X_n)=0$ f"ur $N-1>n$.

Zu 3. \ Nach 1. und 2. gilt $\langle W_5, W_6, \dots \rangle \subset\ker\GU\mid_{\OSUC}$.
Daher l"a"st sich $\GU$ "uber $\OSUC/\langle W_5, W_6, \dots \rangle\iso\C[W_2,W_3,W_4]$
faktorisieren: $\GU=\tilde\GU\circ\pi$. Wegen Satz 2.2.3 ist die Abbildung
$\tilde \GU:\C[W_2,W_3,W_4]\to\C[B,C,D]$ ein Isomorphismus, d.h.~insbesondere
ist $\ker \GU\mid_{\OSUC}\subset\ker\pi=\langle W_5, W_6, \dots \rangle$.
\BOX


Man kann daher die \MGFS $W_5, W_6, \dots $ f"ur die in
Definition 2.3.3 verwendeten \MGFS $V_5, V_6, \dots $ nehmen.

\SN{2.4.4}\index{$I_*^{N,1}$}\index{$J_*^N$} Seien $W_5, W_6, \dots$ die getwistet projektiven B"undel aus Abschnitt
1.4. Dann gilt
$$ I_*^{N,1}=J_*^{N}=\langle \CP_{N-1}, \widetilde\CP_{N+1,1},W_5, W_6, \dots \rangle = \ker\tilde\GN.$$
\SE

{\bf Beweis:} Die Inklusionen $\langle\CP_{N-1}, \widetilde\CP_{N+1,1}, W_5, W_6, \dots\rangle
\subset I_*^{N,1}\subset\ker \tilde\GN$ und \break $\langle\CP_{N-1}, \widetilde\CP_{N+1,1},
 W_5, W_6, \dots\rangle \subset J_*^N\subset\ker \tilde\GN$ sind nach den
Konstruktion in Abschnitt 1.4
und Korrollar 2.4.2 klar; zu zeigen bleibt $\ker \tilde\GN\subset
\langle\CP_{N-1}, \widetilde\CP_{N+1,1}, W_5, W_6, \dots\rangle$.\hfil\break
Nach Definition ist $\tilde\GN=\pi_N\circ\GU$ und nach Satz 2.3.5
gilt $\ker \pi_N =\langle R_{N-1},R_{N+1}\rangle$. Der Kern von $\GU$ auf ganz $\OUC$
ist wegen $\GU(W_1)=A$ und Punkt 3 des Beweises des letzten Satzes
gerade $\langle W_5, W_6, \dots \rangle$. Es gen"ugt daher, die Inklusion
\begin{equation}\label{inkl}
\langle \GU(\CP_{N-1}),\GU(\widetilde\CP_{N+1,1})\rangle\supset
\langle R_{N-1},R_{N+1}\rangle=\ker \pi_N
\end{equation}
von Idealen in $\C[A,B,C,D]$ zu zeigen.

Da $\GU(\CP_{N-1})$ in $\langle R_{N-1},R_{N+1}\rangle$
liegt und das Gewicht $N-1$ hat, gibt es ein $\alpha\in\C$ mit
\begin{equation}\label{rel1}
\GU(\CP_{N-1})=\alpha\cdot R_{N-1}.
\end{equation}
Wir m"ussen $\alpha\not=0$ zeigen.
Das $\chi_y$-Geschlecht ist, wie wir im vergangenen Abschnitt gesehen haben
(Satz 2.3.6 und 2.3.7 (ii)), ein Spezialfall des universellen elliptischen Geschlechtes,
l"a"st sich also "uber dieses
faktorisieren. W"are $\GU(\CP_{N-1})=0$, so w"urde daher das $\chi_y$-Geschlecht
von $\CP_{N-1}$ verschwinden. Es ist aber $\chi_y(\CP_{N-1})=\frac{1-(-y)^N}{1-(-y)}\not=0$
(vgl.~\cite{hihab}, S.~15) --- ein Widerspruch.

Ebenso liegt das Polynom $\GU(\widetilde\CP_{N+1,1})$ vom Grad $N+1$
in $\langle R_{N-1},R_{N+1}\rangle$, so da"s es sich in der Form
\begin{equation}\label{rel2s}
\GU(\widetilde\CP_{N+1,1})=(c\cdot A^2+ c'\cdot B)\cdot R_{N-1}+c''\cdot R_{N+1}
\end{equation}
mit $c,c',c''\in\C$ darstellen l"a"st. Da  $\GU(\CP_1)=A$,
$\GU(\CP_2)=\frac{3}{8}A^2-\frac{1}{16}B$ und $\GU(\CP_{N-1})=\alpha\cdot R_{N-1}$
mit $\alpha\not=0$ ist, k"onnen wir f"ur (\ref{rel2s}) auch
\begin{equation}\label{rel2}
\GU(\widetilde\CP_{N+1,1})=(\beta\cdot\GU^2(\CP_1)+\gamma\cdot\GU(\CP_2))\cdot\GU(\CP_{N-1})+
\delta\cdot R_{N+1}
\end{equation}
mit $\beta,\gamma,\delta\in\C$ schreiben.
Wir zeigen wieder $\delta\not=0$. Dazu
nehmen wir das Gegenteil an und versuchen $\beta$ und $\gamma$ zu berechnen.
Wenn wir die Gleichung (\ref{rel2}) f"ur das $\chi_y$-Geschlecht spezialisieren, erhalten
wir
$$(\beta(1-y)^2+\gamma(1-y+y^2))(\frac{1-(-y)^N}{1-(-y)})=\chi_y(\widetilde\CP_{N+1,1})
=\frac{(-y)^1-(-y)^{N+1}}{1-(-y)}.$$
F"ur das letzte Gleichheitszeichen vgl.~\cite{hata}, S.~717. Es folgt
$\beta (1-y)^2+\gamma(1-y+y^2)=-y$ und damit $\beta=1$ und $\gamma=-1$.
\hfil\break
Setzen wir $B=C=D=0$, erhalten wir einen anderen Spezialfall $\varphi_A$ von
$\GU$. Die zugeh"orige charakteristische Potenzreihe ist dann
$Q(X)=e^{\frac{A}{2} x}$ (s.~Def.~2.3.8). F"ur $A=2$ ergibt sich
$\varphi_A(\CP_n)=(e^g)^{n+1}[\CP_n]=$ Koeff. von $g^n$ in  $e^{(n+1)g}=
\frac{(n+1)^n}{n!}$, sowie $\varphi_A(\widetilde\CP_{N+1,1})=(e^g)^{N+1}\cdot(e^{-g})^1
[\widetilde\CP_{N+1,1}]=(-1)\cdot$ Koeff. von $g^{N+1}$ in $e^{Ng}=-\frac{N^{N+1}}{(N+1)!}$. Wenden
wir dies auf (\ref{rel2}) f"ur $\delta=0$ an, so gilt $((\frac{2}{1!})^2-\frac{3^2}{2!})\frac{N^{N-1}}
{(N-1)!}=-\frac{N^{N+1}}{(N+1)!}\Leftrightarrow \frac{1}{2}=\frac{N}{N+1}$.
Widerspruch f"ur $N\ge 2$.

Also ist $\delta\not=0$, und die
Gleichungen (\ref{rel1}) und (\ref{rel2}) ergeben zusammen die Behauptung.\BOX

Die Ideale $I_*^{N,t}$ mit $t\mid N$, $t\not=N$ gen"ugen der Inklusion
$$ I_*^{N,t}\subset\bigcap_{{\scriptstyle n\mid N \atop \scriptstyle n>1}
\atop \scriptstyle n\notmid t}\ker \varphi_n, $$
wie sich analog zu Korollar 2.4.2 aus Satz 2.4.1 ergibt.

{\bf Vermutung:}\index{$I_*^{N,t}$}{ \sl F"ur alle Paare $(N,t)$ mit $N\ge 2$, $t\mid N$, $t\not=N$
gilt
$$ I_*^{N,t}=\bigcap_{{\scriptstyle n\mid N \atop \scriptstyle n>1}
\atop \scriptstyle n\notmid t}\ker \varphi_n. $$}

Bei festem $N$ w"are der Idealverband der $I_*^{N,t}$ damit isomorph zum
Teilerverband von $N$.
Die \MGFS, die f"ur $t>1$ zur Erzeugung der Ideale $I_*^{N,t}$ fehlen, k"onnen
wegen des letzten Satzes keine getwistet projektive B"undel sein.

Um die Ideale $I_*^N$ und $I_*^{SU}$, die von zusammenh"angenden $N$- bzw.~$SU$-\MGFS
mit effektiver \SO erzeugt werden (keine Einschr"ankung an den Typ),
zu berechnen, ben"otigen wir noch den folgenden

\SN{2.4.5 (Hattori \cite{hat} und Kri\v cever \cite{kri})}
Sei $X$ eine $N$-\MGF mit effektiver \SO. Dann gilt
$\chi(X,K^{k/N})=0$ f"ur $0<k<N$.
\SE

{\bf Korollar 2.4.6:\enspace}\sl Es gilt die Inklusion
$I_*^N\subset\ker \tilde A_N$.\SE

{\bf Beweis:} Lemma 2.3.10.\BOX

\SN{2.4.7}\index{$I_*^{SU,0}$}\index{$I_*^{SU}$} Seien $W_3, W_4, W_5, \dots$ die \MGFS aus Abschnitt 1.4. Dann gilt
$$ I_*^{SU,0}=I_*^{SU}=\langle W_3, W_4, W_5, \dots \rangle = \ker \hat A \mid_\OSUC.$$
\SE

{\bf Beweis:}
Wir zeigen
\begin{eqnarray*}&1.&\quad \langle W_3, W_4, W_5, \dots \rangle\subset I_*^{SU,0},\cr
&2.&\quad I_*^{SU,0}\subset I_*^{SU},\cr
&3.&\quad I_*^{SU}\subset \ker\hat A\mid_\OSUC,\cr
&4.&\quad \ker\hat A\mid_\OSUC\subset\langle W_3, W_4, W_5 \dots \rangle.\qquad\qquad\qquad\qquad\end{eqnarray*}

Zu 1. \ Nach Lemma 1.4.3 und Lemma 1.4.4 besitzen $W_3$ und $W_4$ effektive
$S^1$-Operationen. Diese $S^1$-Operationen m"ussen den Typ $0$ haben, denn
$W_3$ und $W_4$ k"onnen f"ur $t\not=0$ nicht in $I^{SU,t}$ liegen, da nach
Satz 2.4.3 $I_*^{SU,t}=\langle W_5,W_6,\dots\rangle$ gilt. Die \MGFS
$W_5$, $W_6$, $\dots$ besitzen nach Lemma 1.4.7 effektive
$S^1$-Operation vom Typ $0$. \hfil\break
Zu 2. \ Nach Definition. (Die Umkehrung $I_*^{SU}\subset I_*^{SU,0}$ ist --- anders
als im Fall der $N$-\MGFS{} --- nicht unmittelbar ersichtlich.)\hfil\break
Zu 3. \ Folgt aus Satz 2.4.5, da $X\in I_*^{SU} \Rightarrow X\in I_*^2$ und
$\ker \chi(X,K^{1/2}) =\ker \hat A$. Statt Satz 2.4.5 kann man auch den
Satz von Atiyah-Hirzebruch \cite{athis1} verwenden, denn jede $SU$-\MGF ist eine
Spin-\MGF.\hfil\break
Zu 4. \ Das homogen geschriebene $\hat A$-Geschlecht ist $\GU$ spezialisiert
zu $A=C=D=0$, was man z.B. aus Definition 2.3.8 sehen kann, indem man $A=0$ setzt.
Nach 1.~bis 3.~gilt $\langle W_3, W_4, W_5 \dots \rangle\subset\hat A\mid_{\OSUC}$.
Daher l"a"st sich $\hat A$ "uber $\OSUC/\langle W_3, W_4, W_5 \dots \rangle=\C[W_2]$
faktorisieren. Da $\hat A(W_2)=B$, folgt
$\hat A\mid_{\OSUC}\subset\langle W_3, W_4, W_5 \dots \rangle$.
\BOX

\SN{2.4.8}\index{$I_*^N$} Seien $W_3, W_4, W_5 \dots$ die \MGFS aus Abschnitt 1.4. Dann gilt
$$ I_*^{N}=\langle \CP_{N-1}, W_3, W_4, W_5, \dots \rangle = \ker\tilde A_N.$$
\SE

{\bf Beweis:}
Die Inklusion $\langle\CP_{N-1}, W_3 ,W_4, W_5, \dots\rangle
\subset I_*^{N}\subset\ker \tilde A_N$ ist nach der Konstruktion in
Abschnitt 1.4 und Korrollar 2.4.6 klar, zu zeigen ist $\ker \tilde A_N\subset
\langle\CP_{N-1}, W_3, W_4, W_5, \dots\rangle$.\hfil\break
Nach Definition ist $\tilde A_N=\mu_N\circ \tilde A$, wobei $\mu_N$ die
Projektion $\mu_N:\C[A,B]\to\C[A,B]/\langle T_{N-1} \rangle$ ist und
$\tilde A$ das universelle elliptische Geschlecht f"ur $C=D=0$, also
$\ker \tilde A= \langle W_3, W_4, W_5, \dots\rangle$.
Wir m"ussen somit $\tilde A(\CP_{N-1})=\alpha\cdot T_{N-1}$ mit $\alpha\not=0$ zeigen.

Das Geschlecht $\tilde A$ geh"ort zur Potenzreihe
$Q(x)=e^{(A/2)x}\frac{\sqrt{B/2}(x/2)}{\sinh(\sqrt{B/2}(x/2))}=x/f(x)$.
Daher ist $\tilde A(\CP_{N-1})=\hbox{Koeff. von $y^{N-1}$ in }g(y)=f^{-1}(y).$
Mit Schatten-Analysis kann man die Behauptung direkt folgern
(vgl.~\cite{kona} und \cite{roman}, S.~67).

Einfacher ist die folgende "Uberlegung, analog zum Beweis von Satz 2.4.4.
Da \break $\tilde A_N(\CP_{N-1})=0$, ist $\tilde A(\CP_{N-1})=\alpha\cdot T_{N-1}$.
Setzen wir $B=0$, folgt wie dort $\tilde A(\CP_{N-1})=\frac{(\frac{A}{2}N)^{N-1}}{(N-1)!}\not=0$
und daher $\alpha\not=0$. \BOX

Das Geschlecht $\tilde A_N$ ist also das einzige rational kobordismustheoretische
Hindernis f"ur effektive \SOS auf $N$-\MGFS.\footnote{Es ist bekannt (vgl.~\cite{law}
Theorem 11.15), da"s das $\tilde A_N$-Geschlecht einer komplexen \NM
verschwindet, falls diese eine Ricci-positive K"ahler Metrik zul"a"st. Gibt es
au"ser $\tilde A_N$ weitere rational kobordismustheoretische Hindernisse ?
Im Falle von orientier\-ten Spin-\MGFS ist das $\hat A$-Geschlecht vermutlich
das einzige derartige Hindernis f"ur Ricci-postive Metriken. Denn vollst"andige
Durchschnitte mit $c_1>0$ lassen nach Yau's Beweis der Calabivermutung
Ricci-positive Metriken zu (vgl.~\cite{law}, S.~298) und sie erzeugen
wahrscheinlich den Kern von $\hat A$ in $\Omega_*^{Spin}\otimes \Q$.}

Wir haben mit Satz 2.4.4 eine geometrische Beschreibung des elliptischen
Geschlechtes der Stufe $N$ gegeben. Umgekehrt hat das rein geometrisch
definierte Ideal $I_*^{N,1}$ etwas mit Modulkurven zu tun:
$$\OUNC / I_*^{N,1}\iso\C[W_1,W_2,\dots]/\langle \CP_{N-1},\widetilde\CP_{N+1,1},W_4,W_5,\dots
\rangle \iso $$
$$ \C[A,B,C,D]/\langle R_{N-1},R_{N+1} \rangle = K(C_N).$$

\section{Die Starrheit des universellen komplexen elliptischen
Geschlechtes bei $S^1$-Operationen auf $SU$-Mannigfaltigkeiten}

In diesem Abschnitt werden wir das universelle komplexe elliptische \GE
einer $SU$-\MGF als Potenzreihe schreiben mit Koeffizienten, die Indizes
des getwisteten Dolbeault-Komplexes sind. Mit Hilfe dieser Darstellung k"onnen
wir die Starrheit von $\GU$ bei \SOS auf $SU$-\MGFS aus dem entsprechenden
Resultat f"ur $\GN$ bei \SOS auf \NMS herleiten. "Aquivalent dazu ist die
Multiplikativit"at von $\GU$ in $U$-Faserb"undel mit $SU$-\MGFS als Faser.
Wir zeigen, da"s sich jedes komplexe Geschlecht  mit dieser Eigenschaft "uber
$\GU$ faktorisieren l"a"st. Schlie"slich geben wir noch eine Charakterisierung
des $\chi_y$-Geschlechtes an.

Das universelle komplexe elliptische \GE einer $SU$-\MGF $X_d$ ist nach
Satz 2.2.4 ein Polynom, das nicht nicht von $A$ abh"angt, so da"s wir im folgenden
$A=0$ setzen d"urfen, solange wir uns auf $SU$-\MGFS beschr"anken. Weiter
hatten wir gesehen (\ref{lsg}), da"s f"ur feste komplexe Zahlen $B$, $C$ und $D$ die
Differentialgleichung (\ref{dgl}) die explizite L"osung
\begin{eqnarray}\label{hx}
h(x) & = & -\frac{1}{2}\frac{\wp_L'(x)+\wp_L'(z)}{\wp_L(x)-\wp_L(z)}
\end{eqnarray}
besitzt, falls die Diskriminante $\Delta=g_2^3-27 g_3^2=-\frac{1}{32}B^3C^2+
\frac{9}{2}BC^2D+\frac{1}{16}B^2D^2-27C^4-8D^3$ des Polynoms $S(y)=y^4-
\frac{1}{4}B y^2+ 4 C y +\frac{1}{64}B^2-2 D$ von Null verschieden ist.
In diesem Falle gibt es genau ein Gitter $L$ in $\C$ und einen Punkt
$z\in \C/L$, $z \not=0$, so da"s $(B,C,D)=(24 \wp_L(z),\wp_L'(z),6\wp_L^2(z)-(1/2)g_2(L))$,
und umgekehrt geh"ort zu jedem solchen Paar $(L,z)$ ein Tripel $(B,C,D)$.

Da ein Polynom in $B$, $C$ und $D$ schon durch seine Werte auf einer nichtleeren
offenen Teilmenge des $\C^3$ festgelegt ist, verlieren wir keine Information. wenn
wir $\GU(X_d)$ nur f"ur solche komplexen Zahlentripel $(B,C,D)$ betrachten, die
zu einem Paar $(L,z)$ geh"oren.

Ersetzen wir das Paar $(L,z)$ durch $(\lambda L,\lambda z)$ f"ur ein
$\lambda \in \C^*$, so geht $(B,C,D)$ in $(\lambda^{-2} B,\lambda^{-3} C,\lambda^{-4} D)$
"uber, die charakteristische Potenzreihe $Q(x)$ von $\GU$ in $Q(x/\lambda)$ und
daher $\GU(X_d)$ in $\lambda^{-d}\cdot \GU(X_d)$ (vgl.~\cite{jung}, S.~6).
Beschr"anken wir uns daher auf normierte Gitter der Gestalt $L=2\pi i(\Z\tau+\Z)$,
$\tau \in \H$, so verlieren wir nur einen uninteressanten Homogenit"atsfaktor,
wenn wir $\GU(X_d)$ als Funktion von $\tau$ und $z$ auffassen. (Man erh"alt
ihn wieder zur"uck, wenn man $\GU$ wieder homogen schreibt.)

Die Weierstra"ssche Funktion $\wp(\tau,z)$, ihre Ableitung $\wp'(\tau,z)$\index{$\wp$}\index{$g_2$,$g_3$}
bez"uglich $z$ und die Gitterkonstante $g_2(\tau)$ sind als Funktionen
auf $\H\times \C$ meromorphe Jacobiformen vom Gewicht $2$, $3$ und $4$
(und Index $0$) zur vollen Modulgruppe $PSL_2(\Z)$ (vgl.~\cite{eizag} und
\cite{himod} Anhang I). Abgesehen von gewissen Regularit"atsbedingungen
ist eine meromorphe Funktion $\Phi:\H\times\C\to\CP_1$ eine Jacobiform vom
Gewicht $k$, falls sie die beiden folgenden Transformationseigenschaften
besitzt:
\begin{eqnarray*}
\qquad & (1) & \quad \Phi(\tau,x+\omega)=\Phi(\tau,x)
\hbox{ f"ur alle $\omega \in 2\pi i (\Z\tau+\Z)$,}\\
&&\qquad\qquad\hbox{d.h.~$\Phi(\tau,.)$ ist bez"uglich dem Gitter
        $2\pi i (\Z\tau+\Z)$ elliptisch}\\
\noalign{\vskip3mm}
\qquad & (2) & \quad \Phi(\frac{a\tau+b}{c\tau+d},\frac{x}{c\tau+d})(c\tau+d)^{-k}
=\Phi(\tau,x)\hbox{ f"ur alle $\left(\begin{array}{cc}a&b\\c&d\end{array}\right)
\in PSL_2(\Z)$.}
\end{eqnarray*}

Sei $\JC$\index{$\JC$} die bez"uglich des Gewichts graduierte $\C$-Algebra
der meromorphen Jacobiformen $\Phi(\tau,x)$ zu $PSL_2(\Z)$, die in der zweiten
Variable nur Pole im Gitter $2\pi i(\Z\tau+\Z)$ besitzen. Da $\wp$, $\wp'$ und
$g_2$ in $\JC$ liegen\footnote{Vermutlich gilt sogar $\JC\iso\C[\wp,\wp',g_2]$.},
k"onnen wir das universelle komplexe elliptische \GE auch als den graduierten
$\C$-Algebrenhomomorphismus
$$\GU:\OSUC\to\JC$$
ansehen.

In Abschnitt 2.3 hatten wir f"ur festes $q=e^{2\pi i\tau}$, $\tau\in\H$ die
bzgl.~$x$ ganze Funktion\index{$\Phi(\tau,x)$}
$$\Phi(\tau,x)=(1-e^{-x})\prod_{n=1}^{\infty}(1-q^ne^{-x})(1-q^ne^x)/(1-q^n)^2$$
definiert.

\LN{2.5.1} Die charakteristische Potenzreihe $Q(x)$ des als Funktion von
$(\tau,z)\in\H\times\C$ aufgefa"sten komplexen elliptischen Geschlechtes
$\GU$ von $SU$-\MGFS besitzt die Produktentwicklung
\begin{eqnarray*}
Q(x)&=&e^{kx}\frac{x\cdot\Phi(\tau,x-z)}{\Phi(\tau,x)\Phi(\tau,-z)}\qquad\qquad\\
&=&e^{kx}\cdot
\frac{x}{1-e^{-x}}(1+ye^{-x})\prod_{n=1}^{\infty}\frac{1+yq^ne^{-x}}
{1-q^ne^{-x}}\cdot\frac{1+y^{-1}q^ne^{x}}{1-q^ne^{x}}\,\cdot\\
&&\qquad\qquad\qquad\qquad\qquad\qquad
\cdot\,\left((1+y)\prod_{n=1}^{\infty}\frac{(1+yq^n)(1+y^{-1}q^n)}
{(1-q^n)^2}\right)^{-1}.
\end{eqnarray*}
Hierbei sei $q=e^{2\pi i \tau}$ sowie $y=-e^z$ gesetzt und $k$ eine von
$\tau$ und $z$ abh"angige Konstante. \SE

Solange wir uns auf $SU$-\MGFS beschr"anken, k"onnen wir den Faktor $e^{kx}$ in
der Produktentwicklung von $Q(x)$ auch weglassen.

{\bf Beweis:} Die in $\C$ meromorphe Funktion $Q(x)=e^{kx}
\frac{x\cdot\Phi(\tau,x-z)}{\Phi(\tau,x)\Phi(\tau,-z)}$ ist gerade so normiert,
da"s $1$ der konstante Term der Taylorentwicklung im Nullpunkt ist. Wegen
Lemma 2.2.1 und (\ref{hx}) m"ussen wir also
\begin{eqnarray}\label{rl}
&&h(x)=\frac{d}{dx}\log f(x)=-k(\tau,z)+\frac{d}{dx}\Phi(\tau,x)-\frac{d}{dx}\Phi(\tau,x-z)\qquad\qquad\qquad\\
\noalign{\vskip2mm}
&&\qquad\qquad\qquad\qquad\qquad\qquad\qquad\qquad\qquad\qquad\qquad\qquad
   =-\frac{1}{2}\frac{\wp'(\tau,x)+\wp'(\tau,z)}{\wp(\tau,x)-\wp(\tau,z)}\nonumber
\end{eqnarray}
zeigen. Dazu zeigen wir, da"s die Ableitung der rechten und der linken Seite
gleich ist, also $-k(\tau,z)$ eine Integrationskonstante ist.

Die rechte Seite von (\ref{rl}) ist eine bzgl.~dem Gitter $L=2\pi i(\Z\tau+\Z)$
elliptische Funktion mit einfachen Polen bei $0$ und $z$ mit den Residuen $1$
und $-1$ (s.~\cite{jung}, Lemma 1.1.1). Die Ableitung hat daher die Haupteile
$-\frac{1}{x^2}$ und $\frac{1}{(x-z)^2}$. Genauer lautet die Laurententwicklung
im Nullpunkt $-\frac{1}{x^2}+\wp(\tau,z)+O(x)$ (\cite{jung}, S.~4).

Um die Ableitung der linken Seite von (\ref{rl}) zu untersuchen, verwendet man
die Beziehung $\frac{d^2}{dx^2}\log \Phi(\tau,x)=-\wp(\tau,x)+c(\tau)$. Denn
$\Phi(\tau,x)$ ist bis auf einen Exponentialfaktor $e^{ax^2+bx}$ die
Weierstra"ssche $\sigma$-Funktion, und f"ur diese gilt $\frac{d^2}{dx^2}\log
\sigma(\tau,x)=-\wp(\tau,x)$ (\cite{himod}, S.~157). Die Ableitung der linken
Seite ist somit $-\wp(\tau,x)+\wp(\tau,x-z)$. Da $\wp(\tau,x)$ die
Laurententwicklung $\frac{1}{x^2}+O(x^2)$ besitzt, stimmen die Hauptteile
auf der rechten und linken Seite "uberein, d.h.~sie sind bis auf eine Konstante
gleich. Diese verschwindet aber, da die konstanten Glieder $\wp(\tau,z)$ und
$\wp(\tau,-z)$ der Laurententwicklungen im Nullpunkt gleich sind. \BOX

Sei $E$ ein komplexes Vektorb"undel vom Rang $r$ "uber einer $U$-\MGF $X$.
Definiere dann
$$\Lambda_t(E):=\bigoplus_{k=0}^{\infty}\Lambda^kE\cdot t^k$$
sowie
$$S_t(E):=\bigoplus_{k=0}^{\infty}S^kE\cdot t^k$$
als die Potenzreihe in $t$, deren Koeffizienten die "au"seren bzw.~symmetrischen
Potenzen von $E$ sind. Ist $c(E)=\prod_{i=1}^{r}(1+x_i)\in H^*(X,\Z)$ die Zerlegung
der totalen Chernklasse in formale Wurzeln, so gelten in $H^*(X,\Z)[[t]]$ f"ur
den Cherncharakter die folgenden Formeln (vgl.~\cite{himod}, S.~15 f.)
\begin{equation}\label{ch1}
\ch(\Lambda_t(E)):=\sum_{k=0}^{\infty}\ch(\Lambda^kE)\cdot t^k
=\prod_{i=1}^{r}(1+t\cdot e^{x_i})
\end{equation}
\begin{equation}\label{ch2}
\ch(S_t(E)):=\sum_{k=0}^{\infty}\ch(S^kE)\cdot t^k
=\prod_{i=1}^{r}\frac{1}{1-t\cdot e^{x_i}}.
\end{equation}

Weiter ben"otigen wir noch die beiden Beziehungen
\begin{equation}\label{ch3}
\ch(E\oplus F)=\ch(E)+\ch(F)
\end{equation}
\begin{equation}\label{ch4}
\ch(E\otimes F)=\ch(E)\cdot\ch(F)
\end{equation}
f"ur ein zweites komplexes \VB $F$ "uber $X$.

"Uber der $U$-\MGF $X$ hat man den durch die $U$-Struktur gegebenen
Dolbeault-Komplex, den man mit dem komplexen \VB $E$ "`twisten"' kann.

F"ur den Index des $\chi(X,E)$ des so erhaltenen "`getwisteten"'
Dolbeault-Komplexes liefert der Atiyah-Singer-Indexsatz die Formel
\begin{equation}\label{asi}
\chi(X,E)=\Td(TX)\cdot\ch(E)[X].
\end{equation}
Dabei ist $\Td(X)$ die totale Toddklasse von $X$ und $\ch(E)$ der
Cherncharakter von $E$ (vgl.~\cite{himod}, Kap.~5).

Schlie"slich setzen wir noch
$$\chi_y(X_d,E):=\chi(X_d,\Lambda_y(T^*X_d)\otimes E):=
\sum_{p=0}^{d}\chi(x_d,\Lambda^pT^*X_d\otimes E)\cdot y^p.$$

\vspace{1cm}
\LN{2.5.2} F"ur eine $d$-dimensionale $SU$-\MGF $X_d$ gilt mit den Bezeichnungen
wie in Lemma 2.5.1 und der Abk"urzung $T$ f"ur das stabil fastkomplexe
Tangentialb"undel von $X_d$ die Darstellung
\begin{eqnarray}\label{chiyl}\index{$\CHIYL$}
\CHIYL&=&\chi_y\big(X_d,\bigotimes_{n=1}^{\infty}\Lambda_{yq^n}T^*
\otimes\bigotimes_{n=1}^{\infty}\Lambda_{y^{-1}q^n}T\otimes
\bigotimes_{n=1}^{\infty}S_{q^n}(T\oplus T^*)\big)\cr
\noalign{\vskip1mm}
&=&\GU(X_d)\cdot\Phi(\tau,-z)^d
\end{eqnarray}
Die linke Seite der Gleichung ist dabei so zu verstehen, da"s man das Produkt
$\bigotimes_{n=0}^{\infty}\Lambda_{yq^n}T^*\break
\otimes\bigotimes_{n=1}^{\infty}\Lambda_{y^{-1}q^n}T\otimes
\bigotimes_{n=1}^{\infty}S_{q^n}(T\oplus T^*)$
in eine Reihe $\sum_{n=1}^{\infty}\sum_{|i| \le c_n}R_{n,i}\,y^iq^n$
entwickelt (was stets eindeutig m"oglich ist, da die Koeffizienten von $q^n$
stets Laurentpolynome in $y$ sind), wobei die $R_{n,i}$ zum stabilen
fastkomplexen Tangentialb"undel von $X_d$ assozierte komplexe
Vektorb"undel sind, und man dann gliedweise den Index berechnet. Ebenso entwickele
man die rechte Seite in eine Reihe in $y$, $y^{-1}$ und $q$.\SE

Der Ausdruck (\ref{chiyl}) wurde mit $\CHIYL$ bezeichnet, da er formal als
das "aquivariante $\chi_y$-Geschlecht des freien Schleifenraumes ${\cal L}X_d$
von $X_d$ interpretiert werden kann, auf dem $q \in S^1$ durch Reparametrisierung
der Schleifen in nat"urlicher Weise operiert.

{\bf Beweis:} Mit Hilfe des Atiyah-Singer-Indexsatzes (\ref{asi}) und den
Formeln (\ref{ch1}), (\ref{ch2}),(\ref{ch3}) und (\ref{ch4})
f"ur den Cherncharakter erh"alt man f"ur die linke
Seite von (\ref{chiyl})
\begin{eqnarray}\label{prodchi}
\quad\CHIYL&=&\Td(X_d)\cdot\ch\bigl(\bigotimes_{n=0}^{\infty}\Lambda_{yq^n}T^*
\otimes\bigotimes_{n=1}^{\infty}\Lambda_{y^{-1}q^n}T\otimes
\bigotimes_{n=1}^{\infty}S_{q^n}(T\oplus T^*)\bigr)[X_d] \\
&=&\prod_{i=1}^d
\frac{x_i}{1-e^{-x_i}}(1+ye^{-x_i})\prod_{n=1}^{\infty}\frac{1+yq^ne^{-x_i}}
{1-q^ne^{-x_i}}\cdot\frac{1+y^{-1}q^ne^{x_i}}{1-q^ne^{x_i}}[X_d],\nonumber\cr
\end{eqnarray}
wobei wieder $x_1$, $\dots$, $x_d$ die formalen Wurzeln der Chernklasse von
$X_d$ sind. Obwohl sich die Strukturgruppe von $TX_d$ nicht notwendig nach
$U(d)$ reduzieren l"a"st, reichen $d$ formale Wurzeln, da der sich aus dem Produkt
ergebende Ausdruck in den Chernklassen f"ur $d$ oder mehr Faktoren der gleiche
bleibt.

Nach Lemma 2.5.1 ist (\ref{prodchi}) bis auf den Skalierungsfaktor
$\Phi(\tau,-z)^d$ das komplexe elliptische Geschlecht $\GU(X_d)$ von $X_d$.
\BOX

Da der Index eines elliptischen Komplexes stets eine ganze Zahl ist, hat man als

{\bf Korollar 2.5.3:}\sl~Die Potenzreihe $\CHIYL=\GU(X_d)\cdot\Phi(\tau,-z)^d$
liegt in $\Z[y,y^{-1}][[q]]$.\SE

Ebenso l"a"st sich $\GU(X_d)$ selbst als Element von $\Z[[y]][y^{-1}][[q]]$
schreiben, denn der Skalierungsfaktor $\Phi(\tau,-z)^{-1}=\frac{1}{1+y}
\prod_{n=1}^{\infty}\frac{(1-q^n)^2}{(1+yq^n)(1+y^{-1}q^n)}$ liegt in
$\Z[[y]][y^{-1}][[q]]$.

Bezeichnen wir mit $\JZ$\index{$\JZ$} diejenigen Jacobiformen von
$\JC$, die eine ganzzahlige Fourierreihenentwicklung in $y$ und $q$ besitzen,
so ist $\GU$ ein graduierter Ringhomomorphismus
$$\GU:\OSU\to\JZ.$$

{\bf Beispiel $W_2$:} Das elliptische \GE von $W_2$ besitzt nach \cite{himod},
S.~144 Lemma 3.3 die folgende Reihenentwicklung:
\begin{eqnarray*}
\GU(W_2)&=& 24\,\wp(\tau,z)\\
&=& 24 \left(\frac{-y}{(1+y)^2}+\sum_{n=1}^{\infty}\bigl(\sum_{d\mid n}
d((-y)^d+(-y)^{-d})\bigr)q^n+\frac{1}{12}\bigl(1-24\sum_{n=1}^{\infty}
\sigma_1(n)q^n\bigr)\right).
\end{eqnarray*}
F"ur das unskalierte elliptische Geschlecht folgt
\begin{eqnarray*}
\chi_y(q,{\cal L}W_2)&=&24\,\wp(\tau,z)\cdot\Phi(\tau,-z)^2\qquad\qquad\qquad\qquad\qquad\qquad\qquad\qquad\qquad\qquad\\
&=&\underbrace{2y^2-20y+2}_{=\chi_y(W_2)}
+\underbrace{(1+y)^2(-20y^{-1}-88-20y)}_{
=\chi_y(W_2,T^*\oplus T\oplus T^*y\oplus T y^{-1})}q\\
\noalign{\vskip2mm}
&&+\underbrace{(1+y)^2(2y^{-2}-220y^{-1}-588-220y+2y^2)}_{\scriptsize
 \begin{array}{c}\scriptsize
=\chi_y(W_2,T^*\oplus T\oplus S^2T^*\oplus S^2T\oplus 2\cdot T^*\otimes T\oplus
(T^*\oplus T^*\otimes(T^*\oplus T))y\quad\\ \scriptsize\quad\oplus
(T\oplus T\otimes(T\oplus T^*))y^{-1}\oplus \Lambda
^2T^*y^2\oplus\Lambda^2Ty^{-2})
 \end{array}
}q^2+\,\dots
\end{eqnarray*}
Man sieht an diesem Beispiel auch, da"s die Abbildung $\GU:\OSU\to\JZ$ nicht
surjektiv ist, da $\frac{1}{2}\GU(W_2)\in\JZ$, es aber, wie in Abschnitt 1.4
gezeigt, keine $SU$-\MGF mit den Chernzahlen von $\frac{1}{2}[W_2]$ gibt.

Operiert auf einer $U$-\MGF $X$ eine $S^1$, die die $U$-Struktur respektiert,
so operiert die $S^1$ auch auf den B"undeln $R_{n,i}$ sowie auf den
Kohomologiegruppen des getwisteten Dolbeault-Komplexes. Identifiziert man
die $S^1$ mit $\{\lambda\in\C\mid |\lambda|=1\}$, so ist der Charakter
des $S^1$-Moduls $H^q(X,R_{n,i})$ f"ur $\lambda\in S^1$ eine endliche
Laurentreihe in $\lambda$. Bildet man schlie"slich die alternierende Summe der
$H^q(X,R_{n,i})$, so bekommt man den als "aqui\-varianten Index bezeichneten
virtuellen $S^1$-Charakter $\chi(\lambda,X,R_{n,i})\in R(S^1)\iso\Z[\lambda,
\lambda^{-1}]$.

In Verallgemeinerung der Geschlechter der Stufe $N$ gilt der

\SN{2.5.4} Das universelle komplexe elliptische Geschlecht ist starr bei \SOS
auf $SU$-\MGFS, d.h.~der virtuelle $S^1$-Charakter $\chi(\lambda,X,R_{n,i})\in\C[\lambda
,\lambda^{-1}]$ ist f"ur alle $n$ und $i$ trivial.\footnote{
Dieser Satz wurde auch unabh"angig von Kri\v cever bewiesen (vgl.~\cite{kriell}).
Satz 2.5.2 beantwortet die dort gestellte Frage nach einer
Darstellung von $\GU$ durch Indizes elliptischer Komplexe.}
\SE\pagebreak[3]

{\bf Beweis:} Mit $R_n:=\sum_{|i|\le c_n}R_{n,i}\cdot y^i$ gilt
$$\GU(X_d)\cdot \Phi(\tau,z)^d = \sum_{n=0}^{\infty}\sum_{|i|\le c_n}
\chi(X_d,R_{n,i})y^iq^n = \chi(X_d,R_{n})q^n.$$

Das elliptische \GE $\GU(X_d)$ einer $SU$-\MGF $X_d$ f"ur $z=\frac{2\pi i}{N}$
bzw.~$y=-e^{\frac{2\pi i}{N}}$ ist gerade das elliptische Geschlecht $\GN$ der
Stufe $N$ (s.~S.~\pageref{fN}, Gleichung (\ref{fN})).

Nach \cite{hiell}\typeout{Hirzebruch Starrheitsreferenz}
(vgl.~auch Satz 2.4.1) ist $\GN$ starr bei
\SOS auf \NMS, also insbesondere bei \SOS auf $SU$-\MGFS. Dabei bedeutet
Starrheit, da"s f"ur alle $n\ge 0$ der virtuelle $S^1$-Charakter
$$\chi(\lambda,X_d,R_n)\in \Z[y][\lambda,\lambda^-1],$$
aufgefa"st als Laurentreihe in $\lambda$ mit Koeffizienten in $\Z[y]$,
$y=-e^{\frac{2\pi i}{N}}$ in $\lambda$ konstant ist.

Nun ist ($y$ wieder unbestimmt)
$$\chi(\lambda,X_d,R_n)=\sum_{|i|\le c_n}\chi(\lambda,X_d,R_{n,i})y^i
=\sum_{|i|\le c_n}\left(\sum_j a_{ij}^{(n)}\lambda^j\right)\cdot y^i.$$
Aufgrund der Starrheit von $\GN$ gilt
$$\sum_{|i|\le c_n} a_{ij}^{(n)}\cdot y^i=0\hbox{\qquad f"ur $j\not=0$ und $-y=e^{\frac{2\pi i}{N}}$.}$$
Da dies f"ur alle $N\ge 2$ gilt, folgt
$$a_{ij}^{(n)}=0\hbox{\qquad f"ur $j\not=0$,}$$
also
$$\chi(\lambda,X_d,R_{n,i})=\sum_j a_{ij}^{(n)}\,\lambda^j=a_{i0}^{(n)},$$
d.h.~$\chi(\lambda,X_d,R_{n,i})$ ist konstant in $\lambda$.
\BOX

{\bf Korollar 2.5.5:} \sl Das universelle komplexe elliptische Geschlecht
ist multiplikativ in $U$-Faserb"undeln mit $SU$-\MGFS als Faser und kompakter
zusammenh"angender Liegruppe von $U$-Automorphismen der Faser als Strukturgruppe.
\SE

Das Korollar folgt auch direkt aus dem entsprechenden Resultat (Satz 2.4.1) f"ur
$U$-B"undel mit $N$-\MGFS als Faser:
F"ur ein Faserb"undel $E$ mit Basis $B$ und $SU$-Faser $F$
gilt f"ur alle $N$ nach Satz 2.4.1
$\GN(E)=\GN(B)\cdot\GN(F)$, d.h.~$[E]-[B]\cdot[F]\in\ker\GN$.
Daher ist $[E]-[B]\cdot[F]\in\bigcap_{N\ge 2}\ker\GN=\bigcap_{N\ge 2}\ker\tilde\GN=\ker\GU$,
also $\GU(E)=\GU(B)\cdot\GU(F)$.\typeout{Beweis Korollar und Anmerkung}

Wie der n"achste Satz zeigt, ist $\GU$ durch diese Eigenschaft eindeutig
charakterisiert.

\SN{2.5.6} Sei $\varphi$ ein komplexes Geschlecht, das multiplikativ in
$U$-Faserb"undeln ist mit $SU$-\MGFS als Faser und kompakter zusammenh"angender
Liegruppe von $U$-Automorphismen der Faser als Strukturgruppe. Dann l"a"st sich
$\varphi$ "uber das universelle komplexe elliptische Geschlecht $\GU$
faktorisieren, d.h.~zu $\varphi:\OUC\to \Lambda_*$ existiert ein graduierter
$\C$-Algebrenhomomorphismus $\mu:\C[A,B,C,D]\to\Lambda_*$ mit
$\varphi=\mu\circ\GU$.
\SE

{\bf Beweis:} Es gen"ugt offenbar, $\ker\GU\subset\ker\varphi$ zu zeigen. Wegen
Satz 2.4.3 und $\GU(W_1)=A$ gilt $\ker \GU =\langle W_5,W_6,\dots \rangle$.
Nach Konstruktion sind f"ur $n\ge 2$ die \MGFS $W_{2n+1}$ und $W_{2n+2}$
getwistet projektive B"undel mit Faser $\widetilde{\CP}_{n,n}$. Ist $g$ der Erzeuger
von $H^*(\widetilde{\CP}_{n,n},\Z)\iso H^*(\CP_{2n-1},\Z)\iso\Z[g]/\langle g^{2n} \rangle$,
so ergibt sich f"ur die Chernklasse $c(\widetilde{\CP}_{n,n})=(1+g)^n(1-g)^n=(1-g^2)^n$.
Es verschwinden daher alle Chernzahlen von $\widetilde{\CP}_{n,n}$ und damit auch alle
Geschlechter. F"ur ein multiplikatives Geschlecht $\varphi$ gilt also
$\varphi(W_{2n+1})=\varphi(\hbox{Basis})\cdot\varphi(\widetilde{\CP}_{n,n})=0$ und
ebenso $\varphi(W_{2n+2})=0$, d.h.~$\langle W_5,W_6,\dots\rangle\subset
\ker\varphi$.
\BOX

Der Koeffizient von $q^0$ in der Potenzreihenentwicklung
des unskalierten universellen komplexen elliptischen Geschlechts
$\GU(X_d)\cdot\Phi(\tau,-z)^d$=$\sum_{n=0}^{\infty}\sum_{|i|\le c_n}\chi(X_d,R_{n,i})y^iq^n$
einer $SU$-\MGF
ist nach Lemma 2.5.2 gerade das $\chi_y$-Geschlecht. Von ihm ist bekannt,
da"s es sogar starr bei \SOS auf beliebigen $U$-\MGFS ist
(vgl.~\cite{hata,himod}) und --- was "aquivalent dazu ist ---
multiplikativ in $U$-B"undeln. Der nachfolgende Satz zeigt, da"s diese
Eigenschaft das $\chi_y$-Geschlecht charakterisiert.

\SN{2.5.7} Sei $\varphi$ ein komplexes Geschlecht, das multiplikativ ist in
$U$-Faserb"undeln mit kompakter zusammenh"angender Liegruppe von $U$-Automorphismen
der Faser als Strukturgruppe. Dann l"a"st sich $\varphi$ "uber das
$\chi_y$-Geschlecht faktorisieren.\SE

{\bf Beweis:} Das homogen geschriebene $\chi_y$-Geschlecht erh"alt man aus dem
universellen komplexen elliptischen \GE, das im $q_1$,$q_2$,$q_3$,$q_4$-Koordinatensystem
geschrieben sei, wenn man $q_3=q_4=0$ setzt. Dies ergibt sich wie folgt:

Sei $\phi_{V(\langle q_3,q_4\rangle)}$ das \GE $\pi\circ\GU$, wobei $\pi$ die Projektionsabbildung
$\pi:\C[q_1,q_2,q_3,q_4]\to\C[q_1,q_2,q_3,q_4]/\langle q_3,q_4\rangle\iso\C[q_1,q_2]$ ist.
Dann bestimmen sich die Geschlechter
$\chi_y$ und $\phi_{V(\langle q_3,q_4\rangle)}$ gegenseitig:

Setzt man $q_1=2(1-y)$ und $q_2=(1+y)^2$, so erh"alt man f"ur das Polynom $S(z)$ aus der
Differentialgleichung (\ref{dgl})
$$S(z)=z^4+\,2(1-y)z^3+\,(1+y)^2z^2.$$
Die L"osung $Q(x)$ der Differentialgleichung (\ref{dgl}) ist gerade die charakteristische Potenzreihe
f"ur das $\chi_y$-Geschlecht (s.~\cite{himod}, S.~124, vgl.~auch Satz 2.3.6 und Satz 2.3.7).
Sind umgekehrt f"ur eine $n$-dimensionale $U$-\MGF in dem Polynom
$\phi_{V(\langle q_3,q_4\rangle)}(X_n)=\sum_{i+2j=n}\alpha_i\,q_1^iq_2^j$
die Unbestimmten $q_1$ und $q_2$ durch $2(1-y)$ bzw.~$(1+y)^2$ ersetzt,
so kann man aus diesem Polynom in $y$ die Koeffizienten $\alpha_i$ wieder
induktiv bestimmen und damit auch $\phi_{V(\langle q_3,q_4\rangle)}(X_n)$.

Es gilt also insbesondere
\begin{equation}\label{chiy}
\ker \chi_y=\ker\phi_{V(\langle q_3,q_4\rangle)}.
\end{equation}
Wir k"onnen daher $\chi_y$ auch
als das Geschlecht ansehen, das zur Variet"at $V(\langle q_3,q_4\rangle)\iso\CP^{1,2}
\subset\CP^{1,2,3,4}$ geh"ort.

Sei nun $\varphi:\OUC\to\Lambda_*$ ein Geschlecht, das multiplikativ in $U$-Faserb"undeln
ist (mit zusammenh"angender kompakter Liegruppe von $U$-Automorphismen der Faser
als Strukturgruppe), also inbesondere in solchen Faserb"undeln mit $SU$-\MGFS
als Faser. Nach Satz 2.5.6 l"a"st sich $\varphi$ "uber $\GU$ faktorisieren:
$\varphi=\mu\circ\GU$ mit $\mu:\C[A,B,C,D]\to\Lambda_*$. Damit sich $\varphi$
"uber $\chi_y$ faktorisieren l"a"st, m"ussen wir wegen (\ref{chiy})
$\ker\mu\supset\langle q_3,q_4 \rangle$ zeigen.

F"ur ein $U$-B"undel $\xi$ mit Totalraum $E(\xi)$, Basis $B(\xi)$ und Faser
$F(\xi)$ ist $\varphi(E(\xi))=\varphi(B(\xi))\cdot\varphi(F(\xi))$ "aquivalent
zu $\GU(E(\xi))-\GU(B(\xi))\cdot\GU(F(\xi))\in\ker\mu $.
Sei $\xi_3$ das projektive B"undel $\CP(K\oplus K^2)$ zu dem \VB{}$K\oplus K^2$
"uber dem $\CP_2$; sei $\xi_4$ das projektive B"undel $\CP(K\oplus\epsilon_{\C}\oplus\epsilon_{\C})$
zu dem \VB{}$K\oplus\epsilon_{\C}^2$ "uber dem $\CP_2$.
Hierbei steht $K$ f"ur das Determinantenb"undel des Tangentialb"undels des
$\CP_2$. Setzt man noch
$[Q_i]:=[E(\xi_i)]-[B(\xi_i)]\cdot[F(\xi_i)]$ f"ur $i=3,4$, so erh"alt man nach einer
einfachen Rechnung --- analog zu der in Abschnitt 1.4 durchgef"uhrten ---
\begin{eqnarray*}
\GU([Q_3])&=&\frac{3}{16}(2A^3-AB+16C)=\frac{3}{4}q_3,\cr
\noalign{\vskip3mm}
\GU([Q_4])&=&\frac{9}{512}(36A^4-20A^2B+384AC+B^2-128D)=\frac{9}{16}q_1q_3+\frac{9}{8}q_4,
\end{eqnarray*}
d.h.~$\ker\mu\supset\langle\,\GU([Q_3]),\GU([Q_4])\,\rangle=\langle q_3,q_4\rangle$.
\BOX

Das Geschlecht $\varphi_{V(\langle \Delta\rangle)}$, das zur Diskriminantenfl"ache
$V(\langle \Delta\rangle)\subset\CP^{1,2,3,4}$ geh"ort (und zur Basisfolge $W_1$ $W_2$, $\dots$
und dem $A$,$B$,$C$,$D$-Koordinatensystem; $\Delta=g_2^3-27g_3^2$ Diskriminante
von $S(y)$), ist das getwistete $\chi_y$-Geschlecht
$\chi_y(\,.\,,K^r),$
"uber das sich alle "`klassischen"' Geschlechter wie Signatur, Toddgeschlecht,
$\hat A$-Geschlecht und Eulercharakteristik faktorisieren lassen. Projiziert
man die Fl"ache $V(\langle\Delta\rangle)$ in die Ebene $\{A=0\}=\CP^{2,3,4}$, so
erh"alt man die
Kurve $\bar V(\langle\Delta\rangle)\subset\CP^{2,3,4}$ vom Grad $12$, deren Verschwindungsideal
ebenfalls von $\Delta\in\C[B,C,D]$ erzeugt wird. Diese Projektion entspricht
der Einschr"ankung des Geschlechtes auf $\OSUC$. Es gilt also
$$\chi_y(\,.\,)\mid_{\OSUC}=\chi_y(\,.\,,K^r)\mid_{\OSUC}=\varphi_{\bar V(\langle\Delta\rangle)}\mid_{\OSUC}.$$

Wir beenden diesen Abschnitt mit einer "Ubersicht der wichtigsten in der
Arbeit verwendeten Geschlechter und zugeh"origen Variet"aten:

\epsfig{file=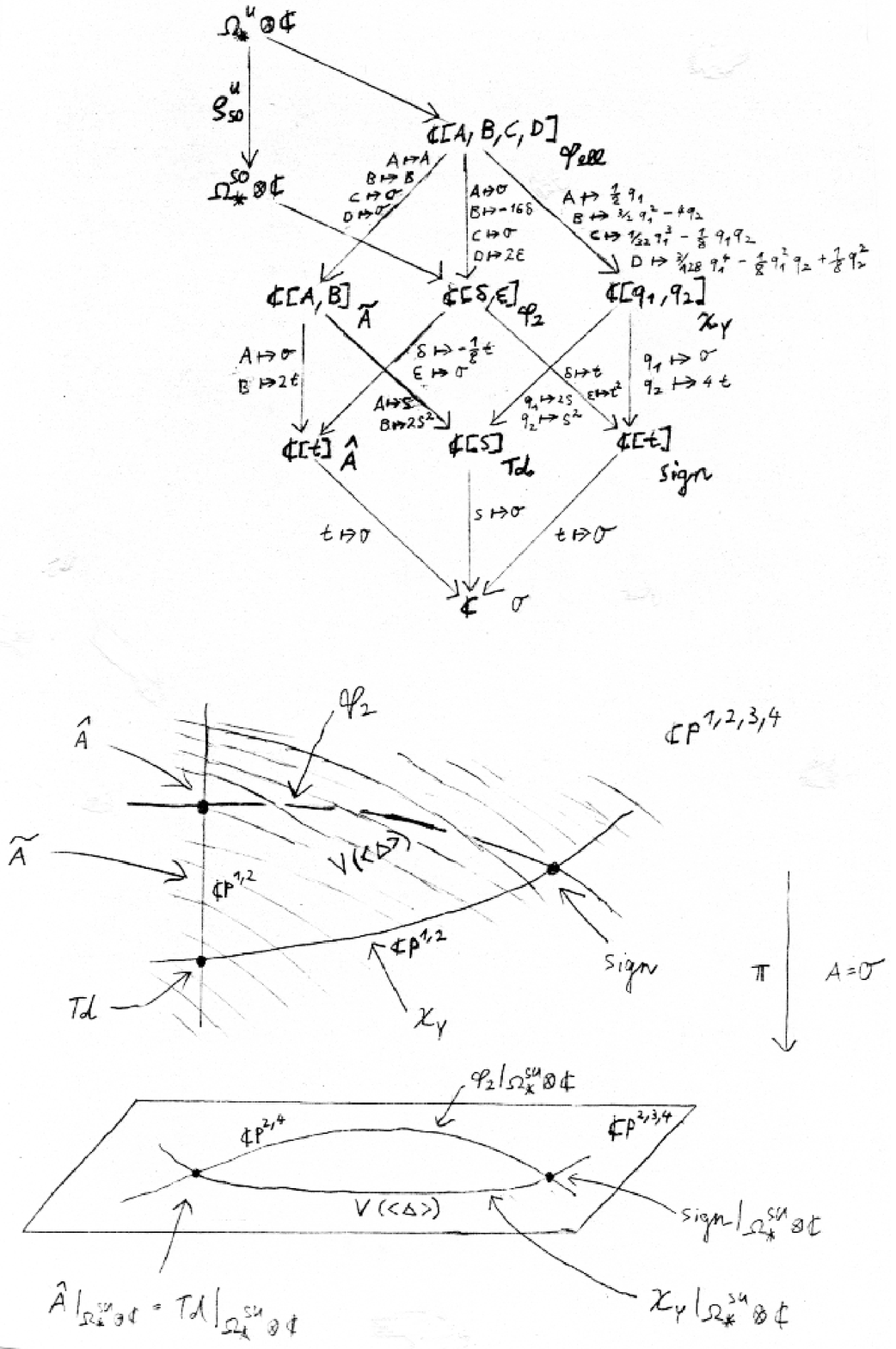,height=200mm,width=150mm}

\chapter{Das Verhalten des elliptischen Geschlechtes beim Aufblasen}

Beim Proze"s des Aufblasens einer kompakten komplexen \MGF $X$
entlang einer komplexen \UMGF $Y$ wird anstelle von $Y$ das
projektive B"undel $ \CP(\nu)$ des Normalenb"undels $ \nu $ von $Y$ in $X$
"`eingeklebt"'. F"ur die so erhaltene \MGF $\tilde X$ ist
bekannt, da"s ihr Toddsches \GE $\Td(\tilde X)$, das zur Potenzreihe $Q(x)=
{\frac{x}{1-e^{-x}}}$ geh"ort, gleich dem von $X$ ist
(vgl.~\cite{hihab}, S.~176). Das Toddsche \GE ist
auch das einzige \GE mit dieser Eigenschaft (vgl.~\cite{himod}, S.~123).
F"ur eine projektiv algebraische \MGF $X$ kann die Invarianz des Toddschen
Geschlechts auch aus der birationalen Invarianz der Hodgezahlen
$ h^{0,p}=h^{p,0}=\dimc\Gamma\left(X,\Lambda^p T^* X \right) $ erhalten
werden, da nach dem Satz von Riemann-Roch-Hirzebruch das arithmetische
\GE $ \chi(X) = \sum_{p \ge 0} (-1)^p h^{0,p} $
gleich dem Toddschen \GE $\Td(X)$ ist.

In diesem Kapitel wird gezeigt, da"s sich das in Kapitel 2 eingef"uhrte
komplexe elliptische \GE $\varphi_N $ der Stufe $N$ beim Aufblasen
nicht "andert, falls $\cod Y \equiv 1 \pmod{N}$ ist. Dazu werden in
Abschnitt 3.1 die ben"otigten Begriffe bereitgestellt, und die bekannte
Beziehung zwischen den Chernklassen von $X$ und $\tilde X$ angegeben. In
Abschnitt 3.2 wird mit Hilfe dieser Formel das elliptische \GE
von $\tilde X$ berechnet.

\section{Das Aufblasen und die Chernklasse}

Ziel dieses Abschnittes ist es, eine Beziehung zwischen der
totalen Chernklasse  der aufgeblasenen \MGF $\tilde X$ und der totalen
Chernklasse von $X$ herzustellen.
Wichtigstes Hilfsmittel ist dabei der Satz von
Riemann-Roch-Grothendieck f"ur Einbe"ttungen von komplexen \MGFS,
wie er von Atiyah-Hirzebruch bewiesen wurde \cite{athigrr}.
Alle in diesem Kapitel auftretenden \MGFS seien als
zusammenh"angend und komplex vorausgesetzt. F"ur eine Erl"auterung der
folgenden Begriffe vergleiche \cite{hihab}, \S~23.

Sei $ K_\omega(M) $ die Grothendieck-Gruppe der koh"arenten analytischen Garben
"uber der kompakten \MGF $M$, sei $ K_\omega'(M)$ die Grothendieck-Gruppe
der holomorphen \VB{} "uber $M$.\index{$K_\omega'(\,.\,)$, $K_\omega(\,.\,)$} Von der kanonischen
Abbildung $j: K_\omega'(M) \to K_\omega(M), [E] \mapsto [\Omega(E)] $, welche
einem \VB die Garbe der Keime von holomorphen Schnitten zuordnet, ist bekannt,
da"s sie, zumindest falls M algebraisch ist,
ein Isomorphismus ist.
Der Cherncharakter eines \VBS $E$ h"angt nur von seiner Klasse
$[E]$ in $ K_\omega'(M) $ ab, ist damit also ein wohldefinierter
Homomorphismus $ ch: K_\omega'(M) \to H^*(M,\Q) $. Die totale Chernklasse
ist ein Homomorphismus in die Gruppe der multiplikativen Einheiten von $H^*(M,\Z)$.
F"ur eine holomorphe Abbildung $h: M \to N $ haben wir die Abbildung
$h^*: K_\omega'(N) \to K_\omega'(M) $, die durch das Zur"uckholen von B"undeln
induziert ist,
sowie die Abbildung $ h_! : K_\omega(M) \to K_\omega(N), [{\cal T}] \mapsto
\sum_{i \ge 0} (-1)^i [R^i h_*({\cal T})] $. Dabei ist $R^i h_*({\cal T}) $
die zur Pr"agarbe $ U \mapsto H^i(h^{-1}(U),{\cal T}) $ assoziierte
Garbe. Zusammen mit dem Zur"uckholen von B"undeln bzw.~dem Pushforward ist
$K_\omega'(.)\iso K_\omega(.)$ ($M$ algebraisch) ein kontravarianter bzw.~kovarianter Funktor,
d.h.~ist $ k: N \to L $ eine weitere Abbildung, so gilt
$ (k \circ h)^* = h^* \circ k^* $ und $ (k \circ h)_! = k_! \circ h_! $.
In der Kohomologie (mit $\Z$- oder $\Q$-Koeffizienten)
gibt es die analoge Abbildung $h^*:H^*(N) \to H^*(M)$
sowie die Umkehr- oder Gysinabbildung $h_*:H^*(M) \to H^*(N)$, die ein
$H^*(N)$-Modul-Homomorphismus ist, falls wir $H^*(M)$ verm"oge $h^*$ als $H^*(N)$-Modul
auffassen, d.h.~f"ur $x\in H^*(N)$,\enspace $y\in H^*(M)$ gilt $h_*(h^*(x)\cdot y)=
x \cdot h_*(y)$.

Ist die Abbildung $h$ eine Einbettung komplexer \MGFS,
so gilt der Satz von
Riemann-Roch-Grothendieck, welcher in der "ublichen Form
geschrieben mit $a \in K_\omega'(M)$ die folgende Gestalt hat:\footnote{
F"ur eine Einbettung $h$ ist die Chernklasse von $h_!(a)\in K_{\omega}(M)$
definiert (s.~\cite{athigrr}).}
\begin{equation}\label{RRG}
ch(h_!(a)) \cdot Td(N) =h_*(ch(a) \cdot Td(M)).
\end{equation}

F"ur die totale Chernklasse selbst gilt die folgende Beziehung:

\SN{3.1.1\enspace (s.~\cite{athigrr})} Sei $h:M \to N$ eine Einbettung von kompakten
komplexen \MGFS, sei $\nu$ das komplexe Normalenb"undel von $M$ in $N$ mit Rang $n$
und $a \in K_\omega'(M)$. Dann gilt
$$ c(h_!a)=1+h_*\bigl(\frac{c(a)*c(\Lambda_{-1}\nu^*)-1}{c_n(\nu)}\bigr).$$
Dabei steht $c(E)*c(F)$ f"ur das universelle Polynom in den Chernklassen von $E$
und $F$ der Chernklasse des Tensorprodukts $E \otimes F$, und es ist $\Lambda_t(\nu^*)
:=\sum_{i=0}^{\infty}\Lambda^i\nu^*\,t^i\in K_{\omega}'(M)[t]$.
\SE

Ist $\nu$ ein Linienb"undel, vereinfacht sich
der Satz 3.1.1 zu
\begin{equation}\label{RRAH}
c(h_!a)=1+h_*\left(\bigl(\frac{1}{v}\bigr)\left(\prod_{i=1}^{\rg a}\bigl(\frac{1+y_i}{1+y_i-v}\bigr)-1
\right)\right),
\end{equation}
wobei $c(a)=\prod_{i=1}^{\rg a}(1+y_i)$ die Zerlegung in die formalen Wurzeln
$y_1,\dots,y_m$ und $v:=c_1(\nu)$ ist. Beachte: der Faktor $1/v$ k"urzt sich
heraus.

Wir kommen nun zur Erl"auterung des Aufblasens.\index{$\tilde X$}
Sei $X$ eine komplex $n$-dimensionale kompakte \MGF und $Y$ eine $k$-dimensionale
kompakte \UMGF. Die Kodimension von $Y$ bezeichnen wir mit $q=n-k$. Sei $U$ eine holomorphe
Koordinatenumgebung eines Punktes $p \in Y \subset X$ mit Koordinaten
$z=(z_1,\dots,z_n)$, wobei die \UMGF $\{ x\in U \mid
z_1(x)=\dots=z_q(x)=0 \}$ gerade $U \cap Y$ sei. Betrachte das
Produkt $U \times \CP_{q-1}$ mit den homogenen Koordinaten $(t_1:\dots:t_q)$
f"ur den $\CP_{q-1}$. In ihm liegt die \UMGF \
$$W=\{(x,t) \in U\times\CP_{q-1}\mid z_i(x)t_j=z_j(x)t_i,\enspace 1\le i<j\le q\}.$$
Die aufgeblasene \MGF ist nun eine komplexe \MGF $\tilde X$ mit Projektionsabbildung
$f:\tilde{X}\to X$,
so da"s f"ur einen Punkt $p\in Y$ und f"ur eine Koordinatenumgebung $U$ um $p$
das Urbild $f^{-1}(U)$ gleich
$W$ ist und der Projektion $f$ die Projektion von $U\times\CP_{q-1}$ auf die erste
Komponente entspricht. F"ur Umgebungen $U$ in $X$, welche $Y$ nicht treffen, sei
$f^{-1}(U)$ biholomorph "aquivalent zu $U$. Die Konstruktion ist unabh"angig von
den gew"ahlten Koordinatenumgebungen, so da"s $\tilde X$ global definiert werden
kann.

Sei $i:Y\to X$ die Einbettungsabbildung, sei $\nu$ das komplex analytische Normalenb"undel von
$Y$ in $X$ und sei $g:\tilde{Y}\to Y$ das assoziierte projektive B"undel mit Faser
$\CP_{q-1}$. Es existiert dann eine Einbettung $j:\tilde{Y}\to \tilde{X}$, so
da"s das folgende Diagramm kommutiert:
\begin{equation}\label{diagramm}
\begin{array}{ccc}
\tilde{Y}   &\mapr{j}&\tilde{X}   \\
\mapd{g}    &        &\mapd{f}    \\
   Y        &\mapr{i}&\enspace X\enspace.
\end{array}
\end{equation}
Das B"undel $\tilde{Y}$ ist eine $k+q-1=n-1$ dimensionale \UMGF von $\tilde{X}$.
Sei $N$ das zugeh"orige eindimensionale Normalenb"undel. Weiterhin betrachten wir
noch
das B"undel $\tilde{Y}^{\Delta}$ entlang der Fasern zu $g:\tilde{Y}\to Y$.

Es gelten offensichtlich die folgenden Beziehungen
\begin{equation}
i^*TX  \cong TY\oplus\nu
\end{equation}
\begin{equation}
T\tilde{Y} \cong g^*TY\oplus\tilde{Y}^{\Delta}
\end{equation}
\begin{equation}
j^*T\tilde{X} \cong T\tilde{Y}\oplus N
\end{equation}
\begin{equation}
g^*i^*TX \cong g^*TY\oplus g^*(\nu) \cong j^*f^*TX,
\end{equation}
dabei sind alle B"undel als differenzierbare B"undel aufzufassen.

Weiterhin ben"otigen wir das folgende Lemma von Porteous \cite{port}, welches
das Tangential\-b"undel von $X$ mit dem von $\tilde{X}$ in Beziehung setzt.

\LN{3.1.2} In $K_\omega(\tilde{X})$ gilt \enspace
 $\Omega(f^*TX)-\Omega(T\tilde{X})=j_!(g^*\nu-N).$
\SE

Sei nun $u \in H^2(\tilde X,\Z)$ die zu $\tilde Y$ Poincar\'e-duale Kohomologieklasse
und $v$ die erste Chernklasse des Normalenb"undels $N$. Es gelten die Beziehungen
\begin{equation}
j^*(u)=v=c_1(N), \qquad j_*(1)=u.
\end{equation}
Weiter hat man f"ur $a,b\in H^*(\tilde Y)$ noch
\begin{equation}\label{produktregel}
j_*(a)\cdot j_*(b)=j_*(a\cdot b\cdot v).
\end{equation}
Wenden wir nun Satz 3.1.1 auf Lemma 3.1.2 an, so erhalten wir
\begin{equation}
c(\tilde X)/f^*c(X)=1+j_*\left(\bigl(\frac{1}{v}\bigr)\left(\prod_{i=1}^{q}\frac{1+x_i-v}
{1+x_i}\cdot\frac{1+v}{1+v-v}-1\right)\right)
\end{equation}
$$
=1+j_*\left(\bigl(\frac{1}{v}\bigr)\left((1+v)\frac{\prod_{i=1}^{q}(1+x_i-v)}
{g^*c(\nu)}-1\right)\right),
$$
wobei $g^*c(\nu)=\prod(1+x_i)$ die Zerlegung in die zur"uckgeholten formalen
Wurzeln des Normalenb"undels $\nu$ ist. Wegen (\ref{diagramm}) gilt
$$j^*f^*c(X)=g^*c(Y)\cdot g^*c(\nu)$$\vskip-1mm
oder
$$ \frac{1}{g^*c(\nu)}=g^*c(Y)\cdot j^*(\frac{1}{f^*c(X)}).$$\vskip1mm
Damit erhalten wir als endg"ultige Formel (vgl.~\cite{port,vdv})

\SN{3.1.3} $$ c(\tilde X )/f^*c(X)=1+j_*\left(\bigl(\frac{1}{v}\bigr)\left(
j^*(\frac{1}{f^*c(X)})\cdot  g^*c(Y)(1+v)\prod_{i=1}^{q}(1+x_i-v)-1\right)\right).$$
\SE

\section{Die Invarianz des elliptischen Geschlechtes}

Der im vergangenen Abschnitt bewiesene Satz 3.1.3 wird
nun zum Beweis des folgenden Satzes verwendet.

\SN{3.2.1}\index{$\varphi_N$} Sei $\tilde X $ die Aufblasung der kompakten komplexen \MGF $X$
entlang der \UMGF $Y$. Falls die komplexe Kodimension von $Y$ der Bedingung
$\cod Y \equiv 1\! \pmod{N}$ gen"ugt, gilt f"ur das komplexe elliptische
\GE der Stufe N die Gleichung $\varphi_N(\tilde X) = \varphi_N(X)$.
\SE

\pagebreak[3]
\bf Bemerkungen: \rm

1. \ F"ur $\cod Y =1 $ ist der Satz trivialerweise richtig, da dann
$\tilde X = X$.

2. \ Aus der Bedingung $ \cod Y\equiv 1 \pmod{N} $ folgt $ \cod Y\equiv 1 \pmod{n} $
f"ur alle $n \mid N$, also gilt dann f"ur alle $ n\mid N $,\enspace $n>1$:
$\varphi_n(\tilde X) = \varphi_n(X)$, d.h.~$\tilde \GN(\tilde X)=\tilde\GN(X)$.
Der Satz bleibt richtig, wenn man
$ \varphi_1$ als das Toddsche \GE definiert.

3. \ Das $\chi_y$-\GE ist, wenn $-y$ eine $N$-te Einheitswurzel
ungleich $1$ ist, bis auf einen Faktor $(1+y)^{\dimc X}$
der Wert von $\varphi_N$ in einer Sorte von Spitzen von $\G1{N}$
(vgl.~Abschnitt 2.3).
F"ur $-y=1$ erh"alt man die Eulercharakteristik, welche keine Invariante beim
Aufblasen ist.

4. \ Die Signatur ist das $\chi_y$-\GE f"ur $-y =-1$. Deren Verhalten
beim Aufblasen ist aber bekannt (vgl.~\cite{fulton}, S.~293). Man hat
$$ \sign(\tilde X) = \cases{
\sign(X), \hfill & \quad \mbox{falls\enspace} $\cod Y$\mbox{\enspace ungerade},\cr
\sign(X)-\sign(Y),& \quad \mbox{falls\enspace} $\cod Y$\mbox{\enspace gerade},
}$$
was im ungeraden Fall mit unserem Satz "ubereinstimmt.

5. \ Eventuell kann man den Proze"s des Aufblasens auch f"ur fastkomplexe oder
stabil fastkomplexe \MGFS definieren. Der Satz 3.1.1 gilt auch
f"ur differenzierbare \MGFS\ (vgl.~\cite{athigrr}), allerdings mu"s man
die Grothendieck-Gruppe $K_\omega(X)$ durch $ K(X) $, die $K$-Theorie von $X$,
ersetzen. Falls Lemma 3.1.2 auch in der $K$-Theorie richtig bleibt, bleibt
auch Satz 3.2.1 richtig. Vielleicht l"a"st sich auf diese Weise auch wieder
eine Charakterisierung der komplexen elliptischen Geschlechter der Stufe $N$
angeben. Das Toddsche \GE kommt allerdings noch hinzu.

Zum Beweis des Satzes ben"otigen wir zwei Lemmata.

\LN{3.2.2} F"ur alle paarweise verschiedene komplexe Zahlen
$x_1$, $\dots$, $x_q$ mit $q\ge 1$ gilt
$$ \sum_{i=1}^{q}\prod_{j \ne i} \frac{x_j}{x_j-x_i}=1 .$$
\SE

\bf Beweis: \rm Betrachte die meromorphe Differentialform
$\frac{dt}{t\cdot(t+x_1)\cdot(t+x_2)\cdot \dots \cdot (t+x_q)} $ auf dem
$\CP_1$. Sie hat Pole in $t=0$, $-x_1$, $-x_2$, $\dots$, $-x_q$ mit den Residuen
$\prod_{i=1}^{q}{\frac{1}{x_i}}$, $\frac{1}{-x_1}\cdot\prod_{j\ne 1}\frac{1}{x_j-x_1}$, $\dots$,
$\frac{1}{-x_q}\cdot\prod_{j\ne q}\frac{1}{x_j-x_q} $. Der Residuensatz liefert
$$ \prod_{i=1}^{q}\frac{1}{x_i}-\sum_{i=1}^{q}\frac{1}{x_i}
\prod_{j \ne i}\frac{1}{x_j-x_i}=0 $$
$$\displaylines{\hfill\Leftrightarrow 1=\sum_{i=1}^{q}\prod_{j \ne i} \frac{x_j}{x_j-x_i}.\hfill\hbox{\BOX}}$$

\LN{3.2.3} Sei $f(x)$ die zum Gitter $L=2 \pi i (\Z\tau+\Z)$ und
o.B.d.A.~dem $N$-Teilungspunkt $\alpha=\frac{2 \pi i}{N}$ geh"orige Funktion
zu dem in Abschnitt 2.2 definierten elliptischen \GE $\GN$ der Stufe~$N$.
F"ur alle modulo dem Gitter $L=2 \pi i (\Z\tau+\Z)$
paarweise in"aquivalente komplexe Zahlen $x_1$, $\dots$, $x_q$ mit $q\equiv 1
\pmod{N}$ gilt dann
$$\sum_{i=1}^{q}\frac{1}{f(x_i)}\prod_{j \ne i} \frac{1}{f(x_j-x_i)}-
\prod_{i=1}^{q}\frac{1}{f(x_i)}=0.$$
\SE

\bf Beweis: \rm Nach Abschnitt 2.3 ist $f$ elliptisch bez"uglich des Gitters
 $\tilde L=2 \pi i (\Z N\tau+\Z)$ und es gilt nach (\ref{trans})
$f(x+2\pi i\tau)=e^{-\frac{2 \pi i}{N}}
f(x)$. $ 1/f $ hat Pole in den Gitterpunkten von $L$ und es ist
$\res_{x=0}1/f=1$.
Betrachte die meromorphe Funktion
$ h(x)=\frac{1}{f(-x)}\cdot\frac{1}{f(x+x_1)}\cdot\dots\cdot\frac{1}{f(x+x_q)}$.
Sie ist elliptisch bez"uglich $L=2 \pi i (\Z\tau+\Z)$, da $q=rN+1$ also
$h(x+2\pi i\tau)=e^{-\frac{2 \pi i}{N}}\cdot \left(e^{-\frac{2 \pi i}{N}}\right)^
{rN+1}\cdot h(x) = h(x)$. Auf der Riemannschen Fl"ache $\C/L$ hat sie Pole
in $0$ und $-x_1$, $\dots$, $-x_q \pmod{L} $ mit den Residuen
$-\prod_{i=1}^{q}{\frac{1}{f(x_i)}}$, $\frac{1}{f(x_1)}\cdot\prod_{j\ne 1}\frac{1}{f(x_j-x_1)}$, $\dots$,
$\frac{1}{f(x_q)}\cdot\prod_{j\ne q}\frac{1}{f(x_j-x_q)}$. Der Residuensatz
liefert
$$\displaylines{\hfill\sum_{i=1}^{q}\frac{1}{f(x_i)}\prod_{j \ne i} \frac{1}{f(x_j-x_i)}-
\prod_{i=1}^{q}\frac{1}{f(x_i)}=0.\hfill\hbox{\BOX}}$$

Da die $x_i$ beliebig waren, gelten die beiden Lemmata auch im Sinne formaler
Potenzreihen.

\bf Beweis von Satz 3.2.1: \rm
Sei $\varphi$ das zu der Potenzreihe $Q(x)=\frac{x}{f(x)}=1+a_1x+a_2x^2+\dots$
geh"orige komplexe \GE und $K_{\varphi}(E):=\sum_{i=0}^{\infty}K_i(c_1(E),
\dots,c_i(E))$ die zugeh"orige multiplikative Folge zum \VB $E$.

Besteht in $H^*(\tilde X)$ die
Gleichung
$$1+c_1+c_2+\dots=1+j_*((\frac{1}{v})(1+b_1+b_2+\dots)-1),$$
wobei $c_i\in H^i(\tilde X)$ und $b_i\in H^i(\tilde Y)$ sei, so folgt aus (\ref{produktregel})
$$K_{\varphi}(c_1,c_2,\dots)=1+j_*((\frac{1}{v})K_{\varphi}(b_1,b_2,\dots)-1).$$

Unter Verwendung der Multiplikativit"at von $K_{\varphi}$ erh"alt man daher aus
Satz 3.1.3 die Gleichung
$$ K_{\varphi}(T\tilde X )/f^*K_{\varphi}(TX)\!=\!1+j_*\left(\!\bigl(\frac{1}{v}\bigr)\left(
j^*(\!\frac{1}{f^*K_{\varphi}(TX)}\!)\cdot g^*K_{\varphi}(TY)Q(v)\prod_{i=1}^{q}Q(x_i-v)-1\right)\!\right)\!.$$
Multiplikation mit $f^*(K_{\varphi}(TX))$ und Auswerten auf $[\tilde X]$ liefert
\begin{eqnarray}\label{phi}
\varphi(\tilde{X})&=&K_{\varphi}(T\tilde{X})[\tilde{X}] = f_*K_{\varphi}(T\tilde{X})[X]\nonumber\\
&=&\!f_*f^*K_{\varphi}(TX)[X]+f_*j_*\left(\!\bigl(\frac{1}{v}\bigr)g^*K_{\varphi}(TY)
   \left(Q(v)\prod_{i=1}^{q}Q(x_i-v)-\prod_{i=1}^{q}Q(x_i)\right)\!\!\right)\![X]\nonumber\\
&=&\varphi(X)+K_\varphi(TY)\cdot g_*\left(\frac{Q(v)}{v}\prod_{i=1}^{q}Q(x_i-v)-
     \frac{1}{v}\prod_{i=1}^{q}Q(x_i)\right)[Y].
\end{eqnarray}
F"ur den Kohomologiering von $\tilde{Y}=\CP(\nu)$ gilt \cite{botu}
$$ H^*(\tilde{Y},\Z) \iso H^*(Y,\Z)[x_1]/\langle (-x_1)^q+c_1(\nu)(-x_1)^{q-1}+\dots+c_n(\nu)\rangle.$$
Dabei ist $x_1=c_1(S)$, und $S$ ist das tautologische Linienb"undel "uber dem
projektiven B"undel $\tilde{Y}$. Dies ist aber gerade unser Normalenb"undel $N$
(vgl.~\cite{botau}, S.~268).
Anderer\-seits spaltet $g^*\nu$ "uber $\tilde{Y}$ in $S$ und Quotientenb"undel,
d.h.~"uber $\tilde{Y}$ gilt $v=x_1$, also $c(g^*(\nu))=\prod_{i=1}^{q}(1+x_i)=
(1+v)\prod_{i=2}^{q}(1+x_i)$. Damit vereinfacht sich (\ref{phi}) zu
\begin{equation}
\varphi(\tilde{X})-\varphi(X)=K_\varphi(TY)\cdot
         g_*\underbrace{\left(\frac{Q(x_1)}{x_1}\left(\prod_{i=2}^{q}Q(x_i-x_1)-
         \prod_{i=2}^{q}Q(x_i)\right)\right)}_{R:=\quad}[Y].
\end{equation}
Um $g_*(R)$ berechnen zu k"onnen, m"ussen wir zur Fahnenmannigfaltigkeit $F(\nu)$
"uber $Y$ "`aufsteigen"', d.h.~zu dem zum $U(q)$-B"undel $\nu$ assoziierten
Faserb"undel mit Faser $U(q)/U(1)^q$. F"ur den Kohomologiering von $F(\nu)$ gilt
\cite{botu}
$$H^*(F(\nu),\Z) \cong H^*(Y,\Z)[x_1,\dots,x_q]/\langle\prod(1+x_i)=c(\nu)\rangle.$$
Betrachte das folgende Diagramm
$$\begin{array}{ccccc}
    & & & &h^*g^*(\nu)=h^*(N) \oplus L_2 \oplus \dots \oplus L_q \\
    & & & &\downarrow                            \\
    & &\qquad\qquad g^*(\nu)=N \oplus Q&  &F(\nu)             \\
    & &\qquad\downarrow & &     \\
\nu & &\qquad\tilde Y &         \\
\downarrow & & & &        \\
Y   & & & &
\end{array}$$
Das zur Fahnenmannigfaltigkeit zur"uckgeholte B"undel $h^*g^*(\nu)$ spaltet in
die Linienb"undel $h^*(N)\oplus L_2 \oplus\dots\oplus L_q$. Wir f"uhren
die Berechnung von $g_*(R)$ auf die
Integration "uber die Faser der Fahnenmannigfaltigkeit zur"uck.
Hierf"ur gibt es die folgende Formel (vgl.~\cite{hihab}, \S~14 und \cite{bohi})
\begin{equation}\label{pi}
g_*h_*(t)=\sum_{\sigma \in S_q}\sign(\sigma)\cdot \sigma(t)/\prod_{i>j}(x_i-x_j).
\end{equation}
Die symmetrische Gruppe $S_q$ operiert dabei auf $t\in\ H^*(F(\nu),\Z)$ verm"oge
Permutation der $x_i$, und $\sign(\sigma)$ bezeichnet das Signum
der Permutation $\sigma$. Weiterhin gilt
\begin{equation}\label{relat}
h_*(\prod_{i>j\ge 2}(x_i-x_j))=(q-1)!\enspace.
\end{equation}

Sei $\varphi$ von nun an das komplexe elliptische \GE der
Stufe $N$. Zum Abschlu"s des Beweises reicht es zu zeigen, da"s $g_*(R)$
verschwindet. Mit (\ref{relat}) erhalten wir
$$\displaylines{(q-1)!\cdot g_*(R)\hfill \cr
\quad =g_*\left(h_*(\prod_{i>j\ge 2}(x_i-x_j))\cdot
\frac{Q(x_1)}{x_1}\left(\prod_{i=2}^{q}Q(x_i-x_1)-
\prod_{i=2}^{q}Q(x_i)\right)\right)\hfill}$$
und wenn wir $h^*(x_i)$ mit $x_i$ identifizieren
$$\displaylines{
\quad =\pi_*\left(\prod_{i>j\ge 2}(x_i-x_j)
\frac{Q(x_1)}{x_1}\left(\prod_{i=2}^{q}Q(x_i-x_1)-
\prod_{i=2}^{q}Q(x_i)\right)\right)\hfill}$$
Anwenden von (\ref{pi})
$$\displaylines{
\quad=\sum_{\sigma \in S_q}\sign(\sigma)\cdot \sigma\left(
\prod_{i>j\ge 2}(x_i-x_j)
\frac{Q(x_1)}{x_1}\left(\prod_{i=2}^{q}Q(x_i-x_1)-\prod_{i=2}^{q}Q(x_i)\right)\right)
/\prod_{i>j}(x_i-x_j)\hfill\cr
\quad=\sum_{\sigma \in S_q} \sigma\left(
\prod_{i>j\ge 2}(x_i-x_j)
\frac{Q(x_1)}{x_1}\left(\prod_{i=2}^{q}Q(x_i-x_1)-\prod_{i=2}^{q}Q(x_i)\right)\right)
/\sigma(\prod_{i>j}(x_i-x_j))\hfill\cr
\quad=\sum_{\sigma \in S_q} \sigma\left(
\frac{Q(x_1)}{x_1}\left(\prod_{i=2}^{q}\frac{Q(x_i-x_1)}{x_i-x_1}
 -\prod_{i=2}^{q}\frac{Q(x_i)}{x_i-x_1}\right)\right)\hfill\cr
\quad=(q-1)!\left(\sum_{i=1}^{q}\frac{1}{f(x_i)}\left(\prod_{j \ne i}\frac{1}{f(x_j-x_i)}-
\prod_{j \ne i}\frac{1}{f(x_j)}\cdot\frac{x_j}{x_j-x_i}\right)\right)\hfill\cr
\quad=(q-1)!\left(\sum_{i=1}^{q}\frac{1}{f(x_i)}\prod_{j \ne i} \frac{1}{f(x_j-x_i)}-
\prod_{i=1}^{q}\frac{1}{f(x_i)}\cdot \sum_{i=1}^{q}\prod_{j \ne i}
\frac{x_j}{x_j-x_i}\right)\hfill}$$
nach Lemma 3.2.1
$$\displaylines{\quad=(q-1)!\left(\sum_{i=1}^{q}\frac{1}{f(x_i)}\prod_{j \ne i} \frac{1}{f(x_j-x_i)}-
\prod_{i=1}^{q}\frac{1}{f(x_i)}\right)\hfill} $$
nach Lemma 3.2.2
$$\displaylines{\quad=0.\hfill\hbox{\BOX}}$$

\twocolumn
\addcontentsline{toc}{chapter}{\protect\numberline{Symbolverzeichnis}}
\vspace*{40pt}
\centerline{\Large\bf Symbolverzeichnis\hfill}
\bigskip
\footnotesize
\begin{list}{}{\leftmargin=0pt\itemsep=0pt\divide\parsep by 3\multiply\parsep by 2}

  \item $(L,z)$, 41
  \item $A$,$B$,$C$,$D$, 39
  \item $C_N$, 43
  \item $F(u,v)$, 34
  \item $I(V)$, 36
  \item $I_*^N$, 18, 56
  \item $I_*^{N,1}$, 54
  \item $I_*^{N,t}$, 18, 55
  \item $I_*^{SU,0}$, 55
  \item $I_*^{SU,t}$, 21, 53
  \item $I_*^{SU}$, 21, 55
  \item $J_*^N$, 20, 54
  \item $J_*^{SU}$, 21, 53
  \item $K(V)$, 36
  \item $K_\omega'(\,.\,)$, $K_\omega(\,.\,)$, 68
  \item $K_n(c_1,c_2,\dots,c_n)$, 34
  \item $N$, 8
  \item $N$-Mannigfaltigkeit, 11
  \item $N$-Struktur, 8, 11
  \item $P_N(y)$, 43
  \item $Q(x)$, 34
  \item $R_{N-1}$, $R_{N+1}$, 43
  \item $S(y)$, 37
  \item $T_{N-1}$, 51
  \item $Td$, 52
  \item $V$, 36
  \item $W_1$, 23
  \item $W_2$, 24
  \item $W_3$, 24
  \item $W_4$, 24
  \item $W_i$, $i\ge 5$, 30
  \item $\CHIYL$, 4, 61
  \item $\CP^{a_1,\dots,a_n}$, 36
  \item $\CP_{N-1}$, $\widetilde\CP_{N+1,1}$, 33
  \item $\GU$, 38
  \item $\Gamma_1(N)$, 42
  \item $\JC$, 58
  \item $\JZ$, 62
  \item $\OSU$, 8
  \item $\OU$, 7
  \item $\OUN$, 12
  \item $\Phi$, 44
  \item $\Phi(\tau,x)$, 42, 58
  \item $\TWIST$, 20
  \item $\chi(\,.\,,K^{k/N})$, 51
  \item $\chi_y$, 51
  \item $\delta,\epsilon$, 44
  \item $\hat A$, 44
  \item $\sign$, 35, 44
  \item $\tilde A$, 51
  \item $\tilde A_N$, 52
  \item $\tilde X$, 69
  \item $\tilde\varphi_N$, 45
  \item $\varphi$, 34
  \item $\varphi_2$, 44
  \item $\varphi_N$, 43, 70
  \item $\varphi_V$, 36
  \item $\varrho^N_*$, 12
  \item $\varrho^U_{SO}$, 35
  \item $\widetilde{\CP}_{p,q}$, 20
  \item $\wp$, 41, 58
  \item $f(x)$, 34
  \item $g(y)$, 34
  \item $g_2$,$g_3$, 41, 58
  \item $h(x)$, 37, 41
  \item $q_1$,$q_2$,$q_3$,$q_4$, 38
  \item $s(X_n)$, 7
  \item $t$, 16

\end{list}

\onecolumn

\addcontentsline{toc}{chapter}{\protect\numberline{Literaturverzeichnis}}

\end{document}